\documentclass[journal]{new-aiaa}
\usepackage[utf8]{inputenc}
\usepackage{textcomp}

\usepackage{graphicx}
\usepackage{amsmath}
\usepackage[version=4]{mhchem}
\usepackage{siunitx}
\usepackage{longtable,tabularx}
\title{Cislunar Space Situational Awareness Constellation Design and Planning with Facility Location Problem
}

\author{Yuri Shimane\footnote{Ph.D. Candidate, Daniel Guggenheim School of Aerospace Engineering.
Student Member AIAA.},
Kento Tomita\footnote{Formerly Ph.D. Candidate, Daniel Guggenheim School of Aerospace Engineering.},
and
Koki Ho\footnote{Dutton-Ducoffe Professor and Associate Professor, Daniel Guggenheim School of Aerospace Engineering. Associate Fellow AIAA. Corresponding Author. Email: kokiho@gatech.edu}
}
\affil{Georgia Institute of Technology, Atlanta, Georgia 30332}

\usepackage{amsmath}
\usepackage{todonotes}
\usepackage{siunitx}
\usepackage{booktabs}
\usepackage{mathtools}
\usepackage{caption}
\usepackage{subcaption}
\usepackage{multirow}
\usepackage{lscape}
\usepackage{algorithm}
\usepackage{algpseudocode}
\usepackage{tikz}
\usetikzlibrary{shapes.geometric, arrows, calc}
\tikzstyle{arrow} = [thick,->,>=stealth]
\newcommand{\Abold}{\boldsymbol{A}}
\newcommand{\bbold}{\boldsymbol{b}}

\newcommand{\rbold}{\boldsymbol{r}}
\newcommand{\vbold}{\boldsymbol{v}}
\newcommand{\fbold}{\boldsymbol{f}}
\newcommand{\xbold}{\boldsymbol{x}}
\newcommand{\zbold}{\boldsymbol{z}}
\newcommand{\suchthat}{\text{s.t.}}
\newcommand{\minorreview}[1]{\textcolor{black}{#1}}
\newcommand{\maximize}[1]{\underset{#1}{\operatorname{maximize}}}

\begin{document}

\maketitle

\begin{abstract} 
Driven by the surmounting interest for dedicated infrastructure in cislunar space, this work considers the satellite constellation design for cislunar space situational awareness (CSSA).
We propose a mixed-integer linear programming (MILP)-based formulation that simultaneously tackles the constellation design and sensor-tasking subproblems surrounding CSSA. 
Our approach generates constellation designs that provide coverage with considerations for the field-of-view of observers. 
We propose a time-expanded $p$-Median problem (TE-$p$-MP) which considers the optimal placement of $p$ space-based observers into discretized locations based on orbital slots along libration point orbits, simultaneously with observer pointing directions across discretized time. 
We further develop a Lagrangian method for the TE-$p$-MP, where a relaxed problem with an analytical solution is derived, and customized heuristics leveraging the orbital structure of candidate observer locations are devised. 
The performance of the proposed formulation is demonstrated with several case studies for CSSA constellations monitoring the cislunar Cone of Shame and a periodic time-varying transit window for low-energy transfers located in the Earth-Moon L2 neck region. 
The proposed problem formulation, along with the Lagrangian method, is demonstrated to enable a fast assessment of near-optimal CSSA constellations, equipping decision-makers with a critical technique for exploring the design trade space. 
\end{abstract}

\section*{Nomenclature}


{\renewcommand\arraystretch{1.0}
\noindent\begin{longtable*}{@{}l @{\quad=\quad} l@{}}
$C_{\rm diff}$               & Lambertian diffuse reflection coefficient \\
$\mathcal{C}_j$     & set of neighbor facility locations of location $j$ \\
CSSA                & Cislunar Space Situational Awareness \\
FLP                 & Facility Location Problem \\
FOV                 & field of view, \SI{}{deg} \\
$f$                 & facility build cost vector \\
LPO                 & Libration Point Orbit \\
$i$                 & pointing direction index \\
$j$                 & facility location index \\
$\mathcal{J}$       & set of candidate facility locations \\
$k$                 & target index \\
$\mathcal{K}$       & set of targets \\
$\ell$              & number of time steps \\
$\boldsymbol{l}$    & line of sight vector \\
$M$                 & visibility tensor \\
$m$                 & number of pointing directions \\ 
$\bar{m}$           & apparent magnitude \\ 
$n$                 & number of candidate facility locations \\
$p$                 & number of facilities to be used \\
$q$                 & number of targets \\
$R_{\rm targ}$      & target radius, \SI{}{m} \\
$\rbold$            & position vector \\ 
$t$                 & time/time index \\
TE-$p$-MP           & time-expanded $p$-median problem \\
$\vbold$            & velocity vector \\
$X_{ijt}$           & binary variable for allocating client $i$ to facility $j$ at time $t$  \\
$Y_{j}$             & binary variable for using facility location $j$ \\
$Z$                 & objective value \\ 
$\eta_{tk}$         & Lagrange multiplier on objective equivalence constraint \\
$\lambda_{jt}$      & Lagrange multiplier on single direction constraint \\
$\phi$              & phase angle, \SI{}{deg} 
\end{longtable*}}

\section{Introduction}
{Cislunar exploration is one of the defining themes in space exploration of this decade. 
Spearheaded by programs such as Artemis and the Commercial Lunar Payload Services (CLIPS), there is a growing number of public and private sector missions to both the lunar surface and various lunar orbits. 
The growing traffic in cislunar space leads to an increasing need for a dedicated cislunar space situational awareness (CSSA) infrastructure - one that can detect and/or monitor the motion of both cooperative and uncooperative assets~\cite{Frueh2021,Vendl2021,Holzinger2021_report,Gordon2024,Wilmer2024,BakerMcEvilly2024}. 
In contrast to Earth orbits, CSSA is complicated by two defining characteristics: the vastness of the volume of interest to be monitored, and the higher nonlinearity of objects in cislunar vicinities. 
Past works have either considered monitoring a general volume of interest~\cite{Bolden2020,Cunio2020,Badura2022,Fahrner2022,Fedeler2022,Frueh2021b,Vendl2021,Visonneau2023}, or focused on monitoring specific translunar and cislunar trajectories of interest~\cite{Dao2020,Fowler2020,Thompson2021,Dahlke2022,Wilmer2022amos,Wilmer2022,Patel2024}. 
The first characteristic renders certain parts of the volume poorly observable from Earth-based observers, thus requiring observations to be conducted from dedicated satellite constellations. 
The second characteristic necessitates more frequent observation and/or information of higher quality to be provided to prediction algorithms to compensate for the higher nonlinearity. 
These defining features of CSSA drive the need for a dedicated constellation infrastructure. 
For a comprehensive review of the CSSA problem, see the review paper by Baker-McEvilly et al.~\cite{BakerMcEvilly2024} and references therein.}

{At its core, the CSSA problem comprises two major problems: constellation design~\cite{Frueh2021,Vendl2021,Holzinger2021_report,Gordon2024,Wilmer2024,BakerMcEvilly2024,Bolden2020,Cunio2020,Badura2022,Frueh2021b,Visonneau2023,Dao2020,Fowler2020,Thompson2021,Dahlke2022,Wilmer2022amos,Wilmer2022,Patel2024} and sensor-tasking~\cite{Fahrner2022,Fedeler2022,Tomita2023}. 
Constellation design consists of finding the optimally distributed network for performing SSA, while sensor-tasking consists of solving for an optimal schedule for each segment within the network to maximize some monitoring performance metric. 
Due to their complexity, each of these problems is commonly tackled individually: constellation design typically assumes volumes of interest can be simultaneously observed at any time as long as it is visible from an observer, while sensor-tasking assumes an already ``fixed'' constellation, and solves for the optimal observation schedule.}

{Future activities in cislunar space remain difficult to predict with certainty, and so is the exact demand for CSSA. 
Observer parameters, such as sensor specifications or the number of assets to be deployed, also remain uncertain. 
To explore the vast trade space for CSSA infrastructure across different combinations of demand and observer parameters, a mathematical framework to efficiently assess the performance of a given constellation is necessitated.}

\subsection{Overview of Existing Approaches for Cislunar Space Situational Awareness Problems}
{Previous works tackling either the constellation design or sensor-tasking problems can be categorized into either nonlinear programming (NLP)-based or linear programming (LP)-based approaches.}

{In an NLP-based approach, a simulation environment is used for parametric studies~\cite{Bolden2020,Vendl2021,Block2022,Dahlke2022,Wilmer2024} or as part of an objective function for algorithms such as evolutionary algorithms (EA)~\cite{Visonneau2023,Clareson2023}, Monte Carlo tree search (MCTS)~\cite{Fedeler2022,Herrmann2022,Klonowski2023}, or machine learning (ML)~\cite{Linares2016,Little2020,Roberts2021,Klonowski2024}, to name a few. 
Due to EA and MCTS handling the objective function as a black box, any arbitrary level of fidelity may be incorporated into the simulation environment without altering the algorithm used. 
Meanwhile, the typically high-fidelity nature of such simulation environments results in objective functions that are computationally expensive, thereby prohibiting extensive parameter trade-off studies. 
Furthermore, the curse of dimensionality renders NLP-based approaches ill-suited for problem formulations involving large numbers of variables, thus necessitating alternative formulations that avoid extensive combinatorial decisions to be made.}

{In an LP-based approach, an LP model, commonly involving binary or integer variables, is constructed and solved~\cite{Shimane2023AMOS,Tomita2023,Patel2023,Patel2024}. 
Nonlinearities surrounding either the constellation design or sensor-tasking, are precomputed to form objectives and/or constraints to the model. 
Mixed-integer LP (MILP) and binary LP (BLP) are well-suited for encoding combinatorial decisions, and a careful formulation can lead to problem formulations that can be efficiently solved by specialized algorithms such as branch and bound (B\&B)~\cite{Land1960,Lawler1966,Balas1968,Morrison2016}, Lagrangian method (LM)~\cite{Geoffrion1974,Guignard2003,Desai2011}, or column generation~\cite{Dantzig1960,Barnhart1998,Lubbecke2005}. 
Additionally, LP algorithms provide an optimality gap, which corresponds to the worst-case degradation of a prematurely terminated suboptimal solution from the predicted optimal solution. 
Such a gap is not readily available in NLP-based approaches such as EA, MCTS, or ML. 
Notable challenges surrounding LP-based approaches include the need for a reasonable discretization scheme to approximate the continuous decision space in observer placement and sensor tasking, and devising an appropriate formulation that is tractable.}

\subsection{Problem Definition and Proposed Approach}
{The constellation design and sensor-tasking subproblems of CSSA are fundamentally coupled. 
The sensor-tasking schedule depends on the architecture, while the actual performance of an architecture is the result of the sensor-tasking schedule. 
To go beyond the approximation of simultaneously monitoring an arbitrary region of space and instead consider observations restricted by the sensor's field-of-view (FOV), one must combine the sensor-tasking of each observer into the constellation design problem.}

{In this work, we propose a logistics-inspired formulation to simultaneously solve the constellation design and sensor-tasking problems. 
Our formulation takes inspiration from the $p$-median problem, a classical formulation for studying the placement problem of facilities in terrestrial logistics. 
We consider the placement of observer spacecraft into orbital ``slots'' along libration point orbits (LPO) and their allocation to a pointing direction along which the observers monitor. 
To accommodate for the time-varying position of observers along their respective orbits, we augment the allocations with a time index, such that the observers are prescribed a specific pointing direction for each time step. 
The resulting formulation is referred to as the time-extended $p$-median problem (TE-$p$-MP). 
To the best of the authors' knowledge, time-expanded $p$-median formulations as applied to space problems are limited to previous iterations of this work~\cite{Shimane2023AMOS,Shimane2023oct_IACLunarFLP}. 
While this work concerns CSSA, an appropriately modified TE-$p$-MP formulation can be employed for various other applications, such as communication relay, planetary defense, or in-orbit servicing, to name a few.}

{While generic B\&B algorithms are known to adequately tackle classical $p$-median problems of moderate size, we also propose a custom LM to compute an approximate solution rapidly. 
LM consists of iteratively solving for an idealized optimal solution to a relaxed problem and a heuristics-based feasible solution, each time updating Lagrange multipliers used to construct the relaxed problem. 
In this paper, an LM with a purely analytical solution to the relaxed problem is derived by decoupling the observer-wise sensor-tasking decision. 
Custom heuristics for the feasible solution are designed to leverage the orbital structure of candidate observer locations. 
We note that the large number of variables, in the order of $10^6$, of the TE-$p$-MP prohibits the use of EAs. 
An NLP-based approach for tackling both constellation design and sensor-tasking problems could consist, for example, of formulating a black-box constellation design problem in the outer loop, with sensor-tasking handled in the inner loop by a pre-trained neural network (NN), or a policy chosen from a predefined set by the outer loop. 
We however opt for the proposed LP-based approach in favor of its ease in constructing a tractable formulation that encodes decisions surrounding the observation schedule. 
With an NLP-based approach, a policy-based approach would result in a limited performance due to the inherently enforced structure by the pre-defined set of policies, while parameterizing the policy to be optimized by an outer loop is non-trivial. 
{Recently, Patel et al.~\cite{Patel2024} also proposed a framework for concurrent CSSA constellation design and sensor-tasking. 
Their approach is a hybrid between NLP and LP, where an outer-layer NLP optimizes the phasing of satellites along pre-defined orbits, and an inner-layer LP assigns discrete targets to observers.
This work differs from~\cite{Patel2024} in the following ways: our formulation optimizes not only for the phasing but also the orbits of observers from a discrete set of options; we also consider target observation by considering a sensor FOV, rather than assuming a one-to-one assignment between observers and targets; finally, our formulation is a single-layer MILP that does not rely on a black-box NLP solver.}
Through the proposed TE-$p$-MP together with the LM, our approach can produce a feasible candidate CSSA constellation design together with a pointing direction schedule for each observer that provides adequate coverage to a given set of demand and observer parameters.}

The remainder of this paper is organized as follows: in Section~\ref{sec:background}, we provide {background information on the facility location problem, the $p$-Median problem, the Lagrangian method, and} the dynamical system within which the TE-$p$-MP is considered. 
This is followed by Section~\ref{sec:demand_definition}, where the static and dynamic SSA demands are defined, along with the observation model adopted in this work. 
Respectively, Sections~\ref{sec:TEpMP} and~\ref{sec:lagrangean_method} are used to introduce the various TE-$p$-MP instances and their Lagrangian relaxation schemes. 
Numerical experiments using the proposed approach are shown in Section~\ref{sec:numerical_results}. Finally, Section~\ref{sec:conclusions} provides conclusions to this work.

\section{Background}
\label{sec:background}
{We first provide background on the facility location problem (FLP), and more specifically the $p$-Median problem. 
These form the basis of the TE-$p$-MP formulated in Section~\ref{sec:TEpMP} for the CSSA application. 
We then briefly introduce of the Lagrangian method for a generic integer LP. 
Finally,} we provide a review of the circular restricted three-body dynamics model adopted to model the motion of observers, along with synodic resonances that exist for LPOs. 

\subsection{Facility Location Problem}
{In Operations Research, the general class of mathematical framework for deciding the location of infrastructure assets and their assignments to clients is known as the facility location problem (FLP). 
FLP is a popular approach for formulating and solving location and allocation problems in terrestrial logistics applications~\cite{Cooper1963,Wesolowsky1973,Hakimi1965,Cornuejols1977,Rosenwein1994DiscreteLT,SimchiLevi2005,flp_WOLF2011,AhmadiJavid2017} as well as space-based applications~\cite{Shimane2023Depot,Shimane2023AMOS,Shimane2023oct_IACLunarFLP}.}

\subsubsection{$p$-Median Problem}
\label{sec:pMP_introduction}
{One variant of the FLP is the $p$-Median problem, where the number of facilities to be deployed is prescribed. 
Let $j = 1,\ldots,n$ index candidate locations where a facility may be built, and $i = 1,\ldots,m$ index the clients to be serviced. 
In the classical $p$-Median problem, let $X \in \mathbb{B}^{m \times n}$ and $Y \in \mathbb{B}^n$ denote 
\begin{align}
    Y_j
    &=
    \begin{cases}
        1 & \text{facility location $j$ is used}
        \\
        0 & \text{otherwise}
    \end{cases}
    \label{eq:pMP_Y_meaning}
    \\
    X_{ij} 
    &=
    \begin{cases}
        1 & \text{client $i$ is allocated to facility location $j$}
        \\
        0 & \text{otherwise}
    \end{cases}
    \label{eq:pMP_X_meaning}
\end{align}
Then, the $p$-Median problem is given by
\begin{subequations}    \label{eq:pMP}
\begin{align}
    \min_{X,Y} \quad& \sum_{j=1}^n f_j Y_j + \sum_{i=1}^m \sum_{j=1}^n c_{ij} X_{ij}
        \label{eq:pMP_objective}
    \\ \suchthat \quad&
    \sum_{j=1}^n X_{ij} = 1 \quad \forall i = 1,\ldots,m
        \label{eq:pMP_service_once}
    \\&
    \sum_{j=1}^n Y_j = p
        \label{eq:pMP_p_facilities}
    \\&
    X_{ij} \leq Y_j \quad \forall i = 1,\ldots,m, \,j = 1,\ldots,n
        \label{eq:pMP_facility_exists}
    \\&
    X_{ij}, Y_j \in \{0, 1\}, \quad \forall i = 1,\ldots,m, \,j = 1,\ldots,n
        \label{eq:pMP_binary_variables}
\end{align}
\end{subequations}
where $f_j$ is the deployment cost of a candidate facility location $j$, and $c_{ij}$ denote the allocation cost of a facility at $j$ to client $i$. 
The objective~\eqref{eq:pMP_objective} is based on the cumulative sum of the deployment costs of the $p$ facilities and the total allocation cost; constraints~\eqref{eq:pMP_service_once} ensure each client is serviced by exactly one facility; constraint~\eqref{eq:pMP_p_facilities} ensures exactly $p$ facility locations are used; constraints~\eqref{eq:pMP_facility_exists} ensure if at least one client is allocated to facility location $j$, then there is a facility at $j$; finally, constraints~\eqref{eq:pMP_binary_variables} ensure $X_{ij}$ and $Y_j$ are binary. 
The TE-$p$-MP formulated in this work is based on~\eqref{eq:pMP}, but with appropriate modifications made (1) to accommodate an objective function based on cumulative target coverage and (2) to include time indices.}
\minorreview{The use of an equality constraint in~\eqref{eq:pMP_p_facilities} instead of an inequality constraint of the form $\sum_{j=1}^n Y_j \leq p$ is reasonable for problem scenarios where adding an additional facility can always improve the objective.
The equality form will be particularly useful when deriving the Lagrangian relaxation in Section~\ref{sec:upperbound_lagrelax}.}

\subsubsection{Including Time in the Facility Location Problem}
{The inclusion of the time dimension into FLP is traditionally coined \textit{multi-period} FLP (MPFLP)~\cite{Wesolowsky1975,Laporte2015,Nickel2019}, where each ``period'' corresponds to the so-called \textit{planning horizon} about which allocation decisions are to be made. 
In the context of space-based applications, the term ``period'' is unfortunately overloaded, as it could also denote the period of an orbit or the synodic period of the Earth-Moon-Sun system. 
To avoid confusion with these definitions, we refer to the proposed formulation with the time dimension as the time-expanded $p$-median problem (TE-$p$-MP). 
The naming convention adopted in this work is also motivated by the fact that even though the MPFLP and the TE-$p$-MP have a mathematical resemblance, the physical meaning of a decision time step is drastically different: following the definition of strategic, tactical, and operational level decisions from Hax and Candea~\cite{hax1984production}, the time-step in the TE-$p$-MP is at the operational level (on the order of hours), while the time-step in the MPFLP is at the strategic or tactical level (on the order of months). 
The MPFLP has been developed in large parts for supply chain and inventory management applications, where demands of commodities typically have predictable, seasonal fluctuations~\cite{Nickel2019}. 
Meanwhile, to this date, the use of MPFLP in space-based applications is limited; the few exceptions include the authors' previous works on the SSA problem~\cite{Shimane2023AMOS} and on cislunar communication relays~\cite{Shimane2023oct_IACLunarFLP}, where facilities have been located to discretized slots along a given set of LPOs, and allocated to observe or provide service to specific targets/clients.}

\subsection{Lagrangian Method for Integer Linear Program}
\label{sec:generic_lagrangian_method}
{
One popular approach to solving complex binary and integer linear programs (ILP/BLP) is the Lagrangian method (LM), which uses Lagrangian multipliers to relax the original problem and iteratively update the multipliers.}

{
Consider a generic ILP given by
\begin{subequations}    \label{eq:generic_ILP}
\begin{align}
    \maximize{\zbold} &\quad \fbold^T \zbold
    \\ \suchthat\quad& \Abold_s \zbold \geq \bbold_s    \label{eq:generic_ILP_simple}
    \\& \Abold_{c} \zbold \geq \bbold_{c} 
        \label{eq:generic_ILP_complicating}
    \\& \zbold \in \mathbb{Z}_{\geq 0}
        \label{eq:generic_ILP_Z}
\end{align}
\end{subequations}
with an optimal solution $Z^*$, where constraints~\eqref{eq:generic_ILP_simple} are \textit{simple}, and constraints~\eqref{eq:generic_ILP_complicating} are \textit{complicating} constraints. 
Determining a given constraint is simple or complicating is problem-specific and requires careful judgment about the problem at hand. 
The fundamental intuition behind LM is to remove complicating constraints by adjoining them using Lagrange multipliers into the objective, in such a way that the original ILP may be decoupled into smaller subproblems that are easier to solve. 
Let $\boldsymbol{\lambda} \geq \boldsymbol{0}$ denote Lagrangian multipliers for each complicating constraint.
The relaxed problem is given by
\begin{equation}    \label{eq:generic_ILP_LagrangianRelaxed}
\begin{aligned}
    \maximize{\zbold} &\quad \fbold^T \zbold
    - \boldsymbol{\lambda}^T \left( \bbold_{c} - \Abold_{c} \zbold \right)
    \\ \suchthat\quad& 
    \text{\eqref{eq:generic_ILP_simple}, \eqref{eq:generic_ILP_Z}}
\end{aligned}
\end{equation}
whose solution is denoted by $Z_{\rm relax}$. 
Since problem~\eqref{eq:generic_ILP_LagrangianRelaxed} is a relaxation of~\eqref{eq:generic_ILP}, the solution to~\eqref{eq:generic_ILP_LagrangianRelaxed} serves as an upper-bound to the original ILP~\eqref{eq:generic_ILP}, which however may not be a feasible solution~\eqref{eq:generic_ILP}. 
Thus, simultaneously, a feasible lower bound solution is reconstructed from $Z_{\rm feas}$ leveraging problem-specific heuristics, such that
\begin{equation}
    Z_{\rm feas} \leq Z^* \leq Z_{\rm relax}
\end{equation}
where $Z_{\rm relax} - Z_{\rm feas}$ represents the optimality gap.
}

{
At each iteration, the LM iteratively updates $\boldsymbol{\lambda}$ to solve for
\begin{align}   \label{eq:generic_ILP_Plambda}
    \min_{\boldsymbol{\lambda}} \quad \text{problem \eqref{eq:generic_ILP_LagrangianRelaxed}}
\end{align}
Since~\eqref{eq:generic_ILP_Plambda} is a piecewise linear concave function with respect to $\boldsymbol{\lambda}$~\cite{SimchiLevi2005}, a commonly adopted approach to solve for iterates of $\boldsymbol{\lambda}$ is the subgradient method~\cite{Nedic2001}.
At a given iteration $h$, $\boldsymbol{\lambda}$ is updated with
\begin{align}   \label{eq:LM_subgradient_multiplier_update}
    \boldsymbol{\lambda}^{(h+1)} &= \boldsymbol{\lambda}^{(h)}
    + s^{(h)} \left(
        \bbold_c - \Abold_c \bar{\zbold}^{(h)}
    \right)
\end{align}
where $s^{(h)}$ is the step-size given by
\begin{equation}
    s^{(h)} = \dfrac{\mu^{(h)} \left(Z_{\mathrm{relax}}^{(h)} - Z_{\mathrm{feas}}^{(h)}\right)}
    {\sum{ \max(0, \bbold_c - \Abold_c \bar{\zbold}^{(h)}) }}
    \label{eq:step_definition}
\end{equation}
In equation~\eqref{eq:step_definition}, $\mu^{(h)}$ is a hyperparameter for scaling the step-size, which is reduced if no reduction of the gap between the lower-bound and upper-bound solutions is observed; initially, $\mu^{(0)} = 2.0$ is typically chosen when employing the subgradient update, with a step reduction factor of 0.5~\cite{SimchiLevi2005}. 
}

\subsection{Dynamical Systems Preliminaries}
\label{sec:dynamical_system_prelim}
\subsubsection{Equations of Motion}
In this work, the circular restricted three-body problem (CR3BP) is used to model the motion of the observer spacecraft. 
{
In the CR3BP, the motion of a spacecraft of negligible mass is studied under the gravitational forces of the Earth and the Moon. 
The two celestial bodies are assumed to be co-rotating about their common barycenter along a circular orbit. 
Following standard notation, we consider the CR3BP dynamics represented in the Earth-Moon rotating frame, defined with its center at the barycenter, the $x$-axis along the Earth-Moon direction, and the $z$-axis aligned with the angular momentum vector of the Earth and the Moon. 
}
Let $\xbold \in \mathbb{R}^6$ denote the state of the spacecraft, composed of $\xbold = [\rbold^T, \vbold^T]^T$, where $\rbold \in \mathbb{R}^3$ denote the position of the spacecraft, $\vbold \triangleq \dot{\rbold} \in \mathbb{R}^3$ denote the velocity of the spacecraft in the rotating frame. 
The CR3BP equations of motion $\fbold(\xbold)$ are given by
\begin{equation}
    \begin{bmatrix}
        \dot{\rbold} \\
        \dot{\vbold}
    \end{bmatrix}
    = 
    \fbold(\xbold)
    =
    \begin{bmatrix}
        \vbold
        \\
        -\dfrac{\mu_1}{\| \rbold_1 \|_2^3} \rbold_1
        -\dfrac{\mu_2}{\| \rbold_2 \|_2^3} \rbold_2 
        - \boldsymbol{\omega} \times \left(\boldsymbol{\omega} \times \rbold\right) 
        - 2\boldsymbol{\omega} \times \vbold 
    \end{bmatrix}
\end{equation}
where $\rbold_1 \in \mathbb{R}^3$ and $\rbold_2 \in \mathbb{R}^3$ are position vectors of the Earth and the Moon with gravitational parameters $\mu_1$ and $\mu_2$, and $\boldsymbol{\omega} = [0, 0, \omega]^T$ is the angular velocity vector of the Earth-Moon rotating frame. 
{
Making use of canonical scales defined in Table~\ref{tab:r3bp_parameters}, $\omega = 1$, $\mu_1 = 1-\mu_2$, and $\rbold_1$ and $\rbold_2$ are given by 
\begin{align}
    \rbold_1 = [x + \mu_2, y, z]^T
    ,\quad
    \rbold_2 = [x - 1+\mu_2, y, z]^T
\end{align}
}
The corresponding state-transition matrix (STM) $\Phi(t,t_0) \in \mathbb{R}^{6\times 6}$, mapping an initial linear perturbation $\delta \xbold \in \mathbb{R}^6$ at time $t_0$ to time $t$, is propagated by the initial value problem
\begin{align}
    \begin{cases}
        \dot{\Phi} (t,t_0) &= \dfrac{\partial \fbold}{\partial \xbold} \Phi(t,t_0)
        \\
        \Phi(t_0,t_0) &= \boldsymbol{I}_6
    \end{cases}
\end{align}

\begin{table}[h]
\centering
\caption{CR3BP parameters}
\begin{tabular}{@{}ll@{}}
\toprule
Parameter                                & Value \\
\midrule
Mass parameter $\mu_2$             & 0.01215058560962404 \\
Length unit $\mathrm{LU}$, \SI{}{km}     & 389703.2648292776 \\
Time unit $\mathrm{TU}$, \SI{}{s}        & 382981.2891290545 \\
Synodic period, \SI{}{day}               & 29.5 \\
\bottomrule
\end{tabular}
\label{tab:r3bp_parameters}
\end{table}

\subsubsection{Libration Point Orbit}
\label{sec:LPO}
A periodic orbit is a path in state space where for some period $T$, $\xbold(t) = \xbold(t+T)$ $\forall t$. In the context of CSSA, periodic orbits about the first and second libration points (L1 and L2), known as libration point orbits (LPOs), are of particular interest due to their motion residing in the cislunar region of interest. 

Multiple known families of LPOs exist in the vicinity of L1 and L2; among the members of each family, a discrete subset of LPOs {have periods approximately equal to an integer ratio of the synodic period $T_{\mathrm{syn}}$, the duration between two consecutive Sun-Earth-Moon alignment}. 
We denote an $M$:$N$ LPO as an orbit with period $T = (N/M) T_{\mathrm{syn}}$. 
In this work, we focus on a subset of such $M$:$N$ resonant LPOs as part of the candidate observer orbits, following previous works reporting their favorable illumination conditions for CSSA applications~\cite{Vendl2021, Visonneau2023}. 
We note that this choice is reminiscent of opting for designing Earth observation constellations with repeating Low-Earth orbits that exhibit repeating ground tracks~\cite{Lee2020}. 
{An actual spacecraft deployed in a cislunar constellation will be subject to full-ephemeris dynamics subject to various perturbations.
Furthermore, due to the Moon's and Earth's eccentricities, $T_{\mathrm{syn}}$ is not constant. 
Thus, in reality, the $M$:$N$ resonance is only an approximation. 
Nevertheless, we assume the spacecraft will be actively controlled via occasional station-keeping maneuvers to remain in the approximate vicinity of the quasi-periodic LPO that geometrically resembles the CR3BP LPO.
In this work, we thus assume the illumination conditions after $M$ revolutions for a given $M$:$N$ resonant LPO are perfectly repeating.
The interested reader is pointed to Klonowski et al.~\cite{Klonowski2023} for an analysis of performance variation due to this drift.}

This work considers the synodic-resonant LPOs in the L1/L2 Lyapunov, Southern and Northern L2 Halo, Southern and Northern Butterfly, Distant Prograde Orbit (DPO), and Distant Retrograde Orbit (DRO) families, with synodic resonances 1:1, 3:2, 2:1, 9:4, 5:2, 3:1, 4:1, and 9:2.
This selection results in 40 LPOs, as shown in Figure~\ref{fig:lpos_considered}. 
The parameters of these LPOs are given in Appendix~\ref{sec:appendix_lpo_conditions}. 
{Each LPO is constructed via a fixed-time variant of the single-shooting differential correction approach~\cite{Howell1984}, using the LPO entry in the JPL Three-Body Periodic Orbit Catalog\footnote{\url{https://ssd.jpl.nasa.gov/tools/periodic_orbits.html}} with the closest period to the sought resonant period as initial guess.}
The choice of which resonant LPOs to include as candidate orbits matters if the relative positions of member satellites within a constellation must be repeated over a $N_{\rm syn}$-multiple of the synodic month, where $N_{\rm syn}$ is the smallest integer such that the product $N_{\rm syn} M/N$ is an integer for all considered $M$:$N$ ratios. 
This translates to whether the sought observer location and pointing schedule are only required to be transient, or if a repeating, steady-state solution is sought. 
With our choice of resonant LPOs, the constellation returns exactly to its initial configuration every 4 synodic months. 
{We note that the TE-$p$-MP formulation developed in later Sections does not explicitly depend on the $N_{\rm syn}$-month repeating pattern of candidate LPOs, and thus one may instead choose to use an arbitrary set of candidate observer LPOs. 
However, the number of discretized time steps necessitated to solve a CSSA constellation and sensor-tasking over an extended duration, for example, lasting multiple years, can be reduced by several folds.}

For each LPO, candidate observer locations, denoted by set $\mathcal{J}$ are defined by discretizing the orbit into $b$ \textit{slots} equally spaced in time, where $b$ denotes the number of slots along a given LPO. We choose $b$ via
\begin{equation}
    b = \operatorname{ceil}\left( \dfrac{P}{\Delta t_b} \right)
    \label{eq:discretize_lpo_slots}
\end{equation}
where $P$ is the period of the LPO, and $\Delta t_b$ is the temporal spacing between each slot. 
\begin{figure}
    \centering
    \includegraphics[width=0.99\linewidth]{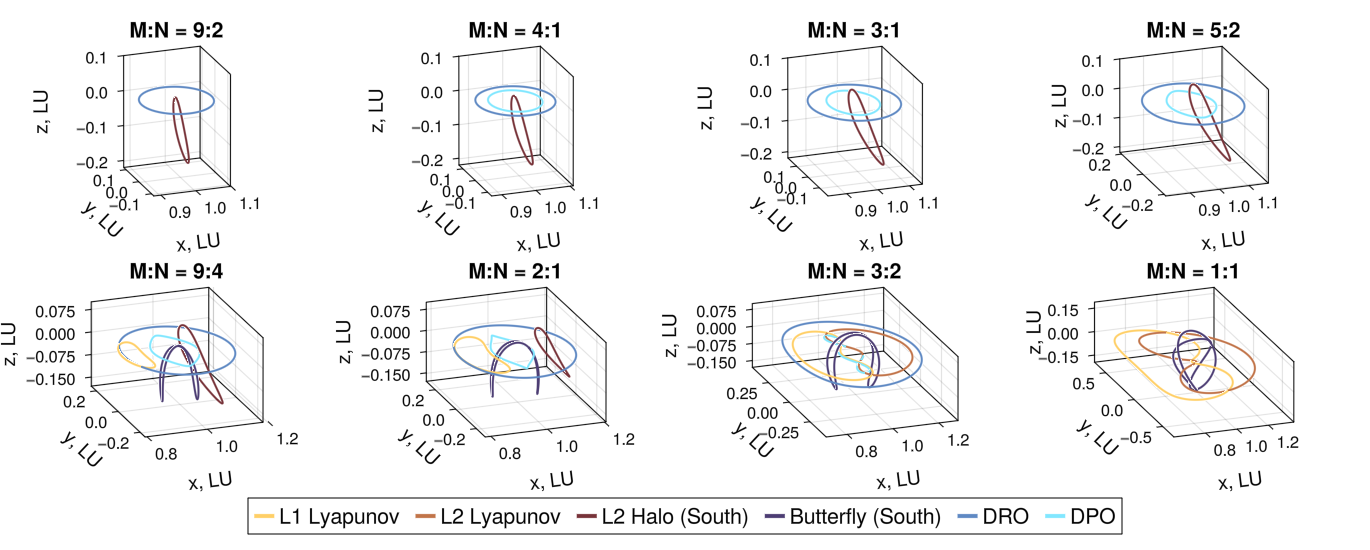}
    \caption{Synodic-resonant LPOs considered as candidate observer orbits, showing only Southern branch for Northern/Southern symmetric LPOs}
    \label{fig:lpos_considered}
\end{figure}

When designing a constellation on LPOs, each member spacecraft's operational cost should also be considered. 
The operational cost is a nontrivial quantity that may include various aspects such as station-keeping cost, navigation requirements, and uplink/downlink frequency; the development of a unifying metric is beyond the scope of this work.
Instead, assuming the primary driving cost comes from station-keeping activities, the linear stability index $\nu$, defined by
\begin{equation}
    \nu = \dfrac{1}{2}
    \left| \lambda_{\max} + \dfrac{1}{\lambda_{\max}}
    \right|
    \label{eq:stability_index}
\end{equation}
where $\lambda_{\max}$ is the largest eigenvalue of $\Phi(T,0)$, provides insight into the station-keeping maneuver $\Delta V$; Folta et al.~\cite{Folta2013} demonstrated a monotonic relationship between $\nu$ and $\Delta V$ on a per-LPO-family basis.
In the context of this work, the operational cost merely serves the purpose of a secondary metric, used only if two different constellations provide the same coverage capability; thus, using $\nu$ as a proxy to represent operational cost is deemed sufficient.

\section{Demand and Visibility Definitions}
\label{sec:demand_definition}
We introduce the demands considered for the CSSA problem. 
Specifically, we first distinguish between \textit{static} and \textit{dynamic} demands; the former consists of coverage requirement that spans the same volume over time, while the latter consists of time-dependent coverage, either to keep custody of known targets or to monitor transfer corridors that are time-dependent. 
In both cases, demands are defined by considering a finite set of targets distributed in space, which approximates a volume of interest for monitoring.
In the static case, the demand for all targets is assumed to be always active, while in the dynamic case, the demand is activated only at time steps when the target exists. 
Consider $k = 1, \ldots, q$ targets and $t = 1,\ldots,\ell$ time steps; we define the demand matrix $D \in \mathbb{B}^{\ell \times q}$ as 
\begin{equation}
    D_{tk} = 
    \begin{cases}
        1 & \text{target $k$ exists at time $t$}
        \\
        0 & \text{otherwise}
    \end{cases}
\end{equation}
Note that in the static case, $D_{tk} = 1$ for all $t$ and $k$, while in the dynamic case, entries of $D_{tk}$ depend on the target motions considered. 
Once the targets are introduced, this Section introduces the visibility model used to determine whether a given target is visible from an observer with a given solar phasing angle. 

\begin{figure}
     \centering
     \begin{subfigure}[b]{0.48\textwidth}
         \centering
         \includegraphics[width=\textwidth]{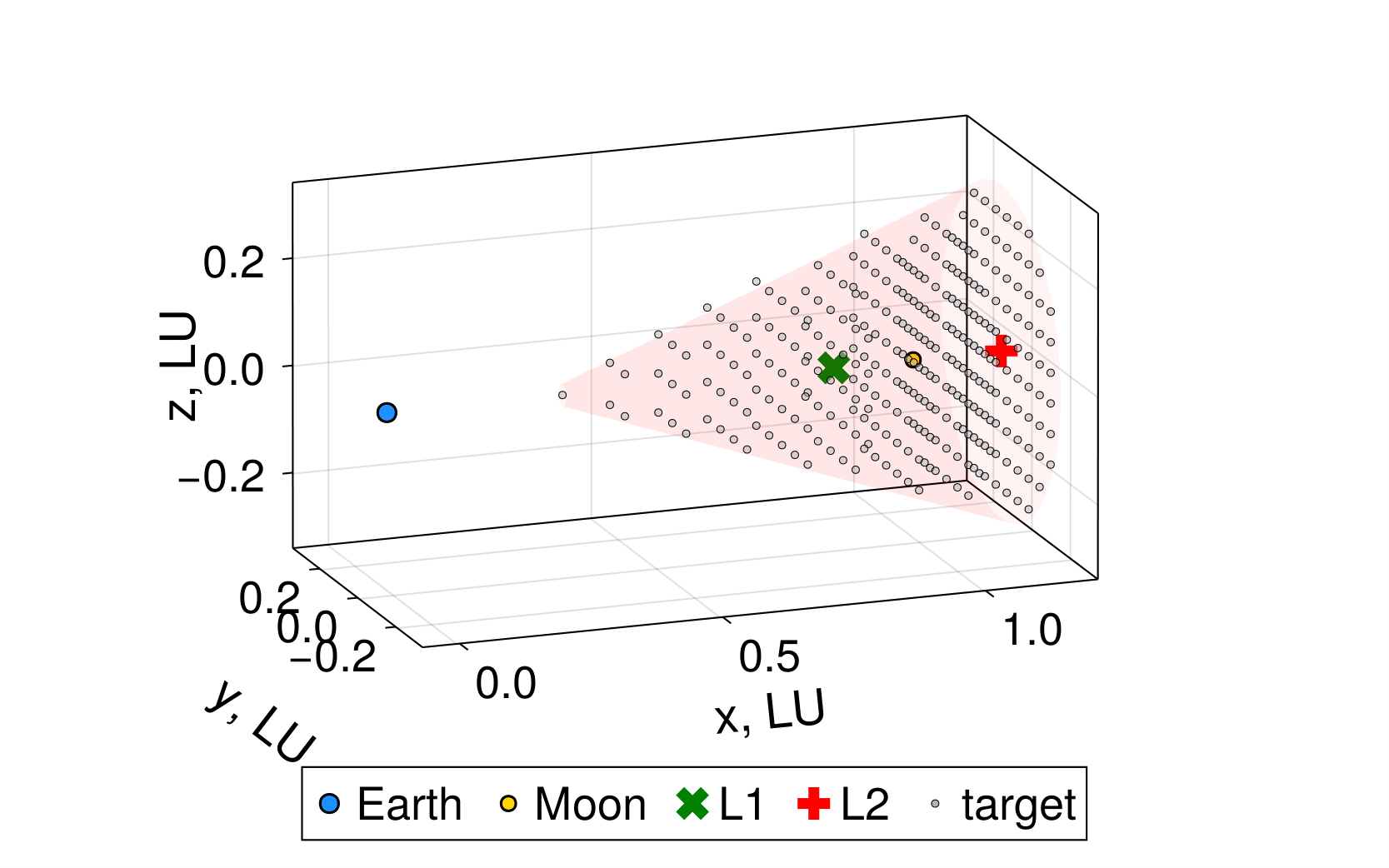}
         \caption{Static target in cone of shame}
         \label{fig:targets_static_cone}
     \end{subfigure}
    \hfill
     \begin{subfigure}[b]{0.48\textwidth}
         \centering
         \includegraphics[width=\textwidth]{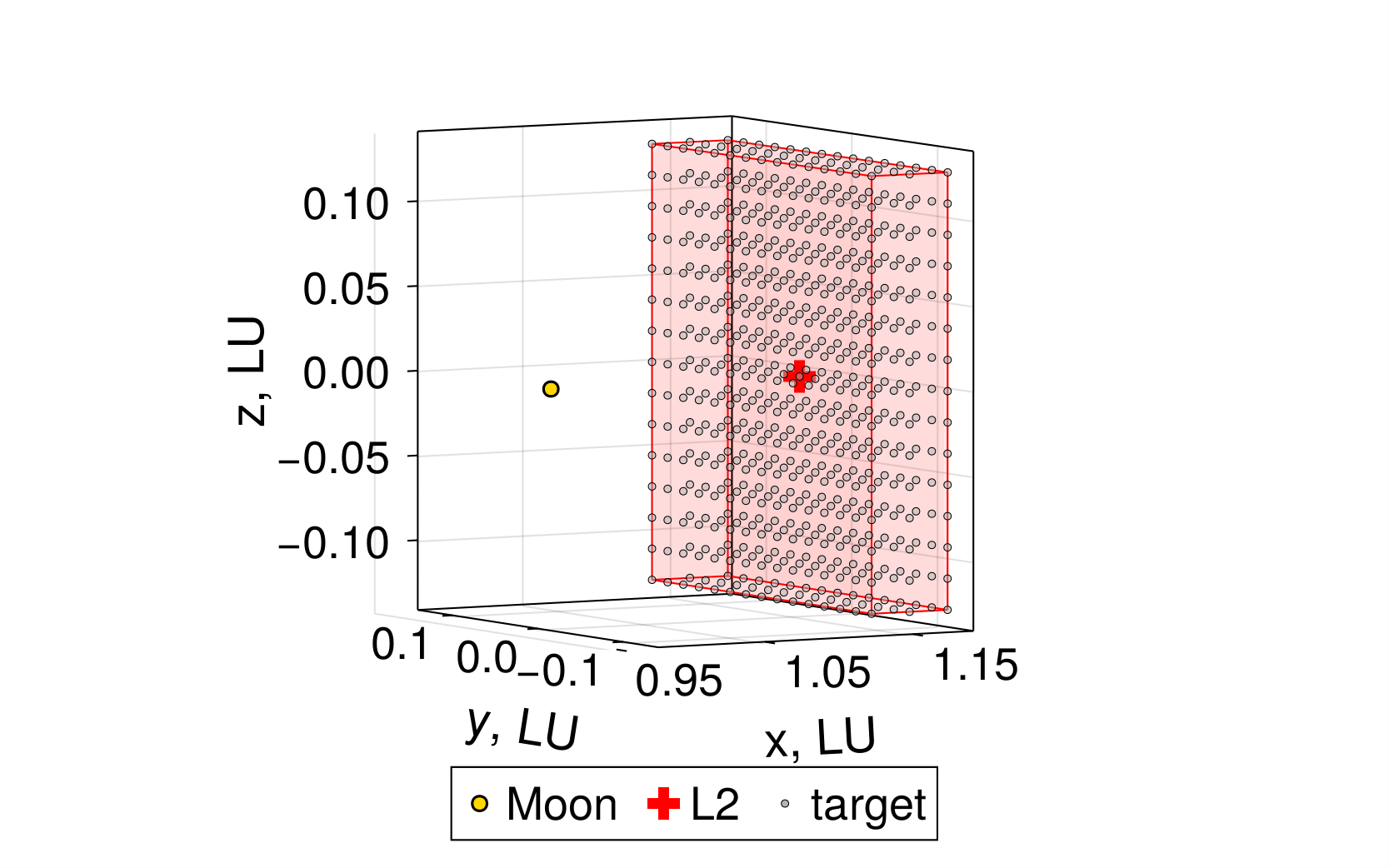}
         \caption{Dynamic target in low-energy transfer window}
         \label{fig:targets_dynamic_let}
     \end{subfigure}
        \caption{Cislunar observation targets for persistent monitoring of volume of interest (red shaded volume), shown in Earth-Moon rotating frame centered at the Earth-Moon barycenter}
        \label{fig:targets}
\end{figure}

\subsection{Static Demand: Cone of Shame}
A static demand is defined by a set of targets distributed in space that remain stationary. While in reality, an object cannot remain at a fixed location in space, each ``target'' in this scenario may be understood as a placeholder for defining a volume of interest. Then, the coverage of such a static demand translates to perpetual coverage of the volume that is spanned by the spatial distribution of these static targets. 

The static demand considered in this work is the so-called \textit{cone of shame}, coined by the Air Force Research Laboratory (AFRL) to designate the lunar exclusion zone, within which tracking objects with Earth-based sensors is particularly challenging~\cite{Roth2020AMOS}.
In this work, we consider a set of target points distributed in a cone between twice the GEO altitude and Earth-Moon L2, with a cone angle of $30^{\circ}$, as shown in Figure~\ref{fig:targets_static_cone}, consisting of 304 targets. 

\subsection{Dynamic Demand: Low-Energy Transfer Transit Window}
A dynamic demand may be defined by a variable number of targets distributed in space for any given instance in time. This formalism may be used to model the observation requirement of a specific trajectory of interest, such as that of the Lunar Gateway, or to monitor a time-varying region of interest. In this present work, we focus on the latter case, specifically looking at the persistent coverage of low-energy transfers (LETs). 

\begin{figure}
    \centering
    \includegraphics[width=0.96\linewidth]{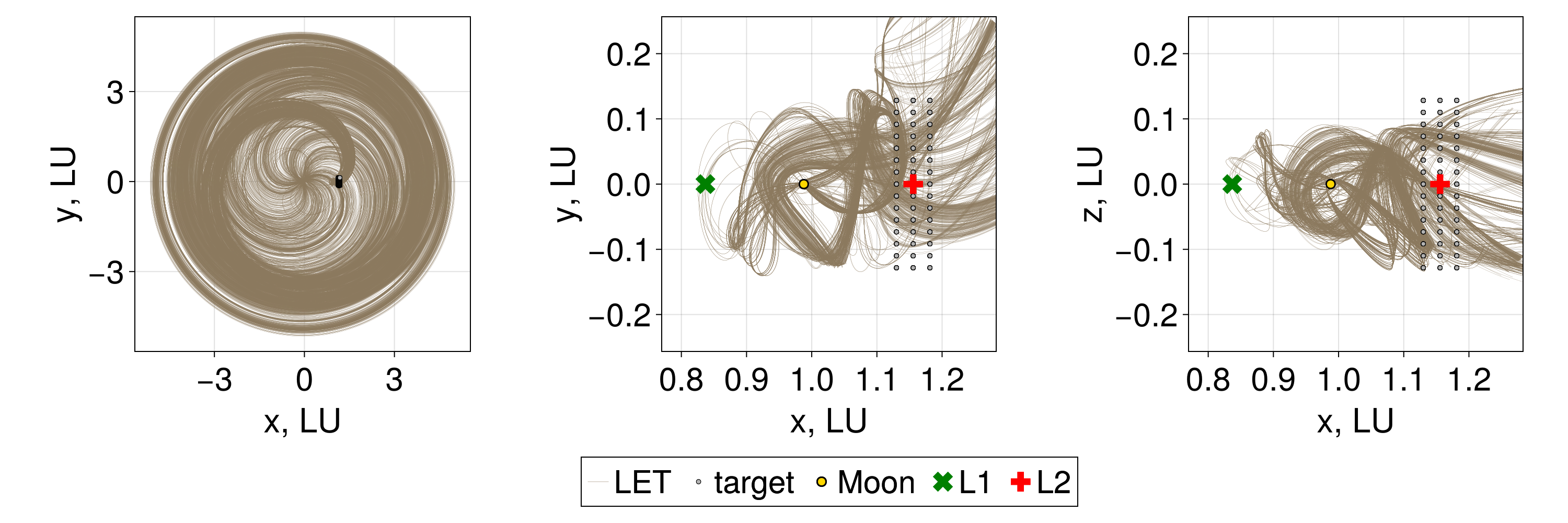}
    \caption{Sample LETs in the Earth-Moon rotating frame centered at the system barycenter. Grey markers indicate transit window monitoring targets.}
    \label{fig:LETs_EMrotating}
\end{figure}

LETs are translunar trajectories that leverage the weak stability boundary of the Earth-Moon-Sun system to reduce the arrival specific energy with respect to the Moon at the cost of longer times of flight~\cite{Koon2001,belbruno2002analytic,jpl-monograph-series12}. 
A sample set of LETs computed with the bi-circular restricted four-body model is shown in Figure~\ref{fig:LETs_EMrotating} in the Earth-Moon rotating frame. 
The left pane in the Figure shows the entire trajectory from perigee to perilune at \SI{100}{km} altitude; the center and right pane shows the portion of each LETs in the lunar vicinity. Due to their high apogee of around $1.5$ million km, LETs are particularly hard to monitor during the majority of their transfer, even with the existence of an in-space SSA network capable of observing areas such as the lunar SOI or the cone of shame. 
One strategic location to monitor LETs is the neck region around Earth-Moon L2, where trajectories pass through before entering the cislunar vicinity. 
The center and right pane in Figure~\ref{fig:LETs_EMrotating} shows this behavior. 
We consider monitoring a volume corresponding to this neck region-based transit window, centered at L2 and spanning \SI{20000}{km} along the $x$-axis and \SI{100000}{km} along the $y$ and $z$-axes of the Earth-Moon rotating frame. 

\begin{figure}
    \centering
    \includegraphics[width=0.28\linewidth]{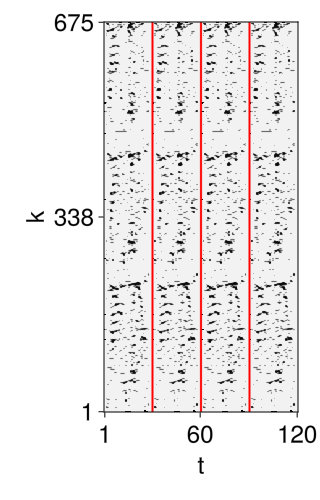}
    \caption{Sparsity of demand matrix $D \in \mathbb{B}^{\ell \times q}$ for LET transit window over 120 time steps corresponding to $4$ synodic months}
    \label{fig:spy_let_demand_matrix}
\end{figure}

The volume of interest is represented by an orthogonal grid of 675 targets, as shown in Figure~\ref{fig:targets_dynamic_let}. 
Since LETs leverage the Sun's tidal effects, they are inherently tied to a specific time of the month when they cross the neck region and arrive in the lunar vicinity. 
Thus, while this orthogonal grid remains static in space, persistent coverage of LETs can be achieved by providing coverage at each time step to portions of the grid where at least one LET transits through.
With a pre-computed database of LETs developed by the authors~\cite{Shimane2023AMOS}, we define $D_{tk}$ such that the existence of the target $k$ at a given time $t$ is based on the existence of a LET that passes near target $k$ at this given time. 
Figure~\ref{fig:spy_let_demand_matrix} shows the sparsity pattern of $D$ over four synodic periods, which repeats every synodic month, as indicated by the red lines marking the end of each month.

\subsection{Metric for Observation Condition}
The observation of a target by an observer must take into account the relative range, illumination conditions as well as the properties of both the sensor and the target. 
We first present the illumination model that computes the apparent magnitude of a target within the FOV of an observer, given the positions of the observer, target, and the Sun, as well as the physical properties of the target. 
Then, we introduce the notion of pointing directions, which considers the FOV of the observer's sensor. 
Finally, the above two conditions are condensed into a boolean visibility flag. 

\subsubsection{Sun Illumination Model}
\label{sec:sun_illumination_model}
\begin{figure}
    \centering
    \includegraphics[width=0.42\linewidth]{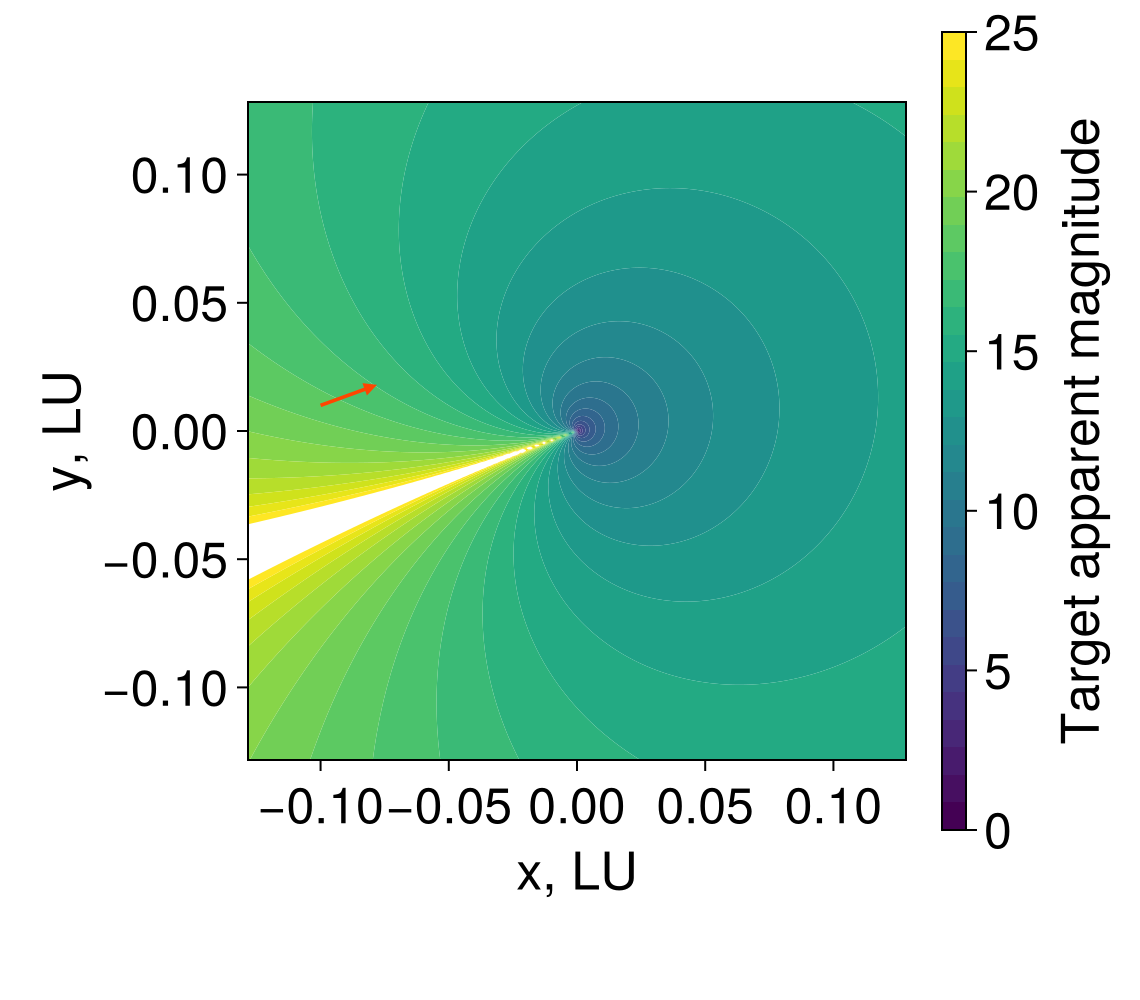}
    \caption{Example visibility magnitude contour with $\phi_S = 200^{\circ}$, with Sun illumination direction aligned with orange arrow}
    \label{fig:visibility_magnitude_contour}
\end{figure}

Let $m_S = -26.74$ denote the apparent magnitude of the Sun, $\rbold_k$ denote the position vector of a target, $\rbold_j$ denote the position vector of an observer, and let $\phi_S$ denote the solar phase angle given by
\begin{equation}
    \phi_S 
    =
    \arccos \left[
        \boldsymbol{l}_{jk}^T \boldsymbol{l}_{Sk}
    \right]
\end{equation}
where $\boldsymbol{l}_{jk}$ is the unit line-of-sight vector from slot $j$ to target $k$, and $\boldsymbol{l}_{Sk}$ is the line-of-sight vector from the Sun to target $k$. 
{
The apparent magnitude of a target object is given by~\cite{Krag1974}
\begin{equation}
    \bar{m}_{\mathrm{target}} (\rbold_j, \rbold_k)
    =
    m_{S} - 2.5 \log _{10}\left[
        \dfrac{R_{\rm targ}^2}{\| \rbold_k - \rbold_j \|^2}
        \left(
            C_{\rm diff} p_{\rm diff}(\phi_S)
            + 
            \dfrac{1}{4} C_{\rm spec}
        \right)
    \right]
    \label{eq:apparent_magnitude}
\end{equation}
where $C_{\rm diff}$ is the diffuse reflectance coefficient, $C_{\rm spec}$ is the diffuse reflectance coefficient, $R_{\rm targ}$ is the target's radius, and $p_{\rm diff}(\phi_S)$ is the diffuse phase function given by
\begin{equation}
    p_{\rm diff}(\phi_S) = 
    \dfrac{2}{3\pi}[\sin (\phi_S)+(\pi - \phi_S) \cos (\phi_S)]
\end{equation}
\minorreview{In this work, we consider targets of $R_{\rm targ} = 2 \,\mathrm{m}$, approximating a medium-sized, robotic lunar spacecraft, with $C_{\rm diff} = 0.2$, and $C_{\rm spec} = 0$.
Values for $C_{\rm diff}$ and $C_{\rm spec}$ are according to~\cite{Vendl2021}, and correspond to a Lambertian sphere model.
An object may exhibit a higher $C_{\rm diff}$, for example within the range $[0.1,1]$ as stated in~\cite{BakerMcEvilly2024}, and the selected value serves as a conservative estimate, which results in dimmer targets.}
}
Figure~\ref{fig:visibility_magnitude_contour} shows the contour of $\bar{m}_{\mathrm{target}}$ centered at the observer with $\phi_S = 200^{\circ}$; the contour demonstrates the strong variation of $\bar{m}_{\mathrm{target}}$ with respect to $\phi_S$, and the weaker variation with respect to the distance $\| \rbold_k - \rbold_j \|$ as the latter increases. 

In this work, the line of sight to the Sun $\boldsymbol{l}_{Sk}$ is assumed to lie on the Earth-Moon orbital plane, rotating clockwise at a constant angular velocity $\omega_S$ about the $z$-axis of the Earth-Moon rotating frame. 
This assumption is based on the small, $\sim 5^{\circ}$ difference between the Sun-Earth and Earth-Moon orbital plane, making a relatively small impact on the illumination conditions. 
Adopting this assumption, the illumination condition is exactly one synodic month-periodic, thereby making a given CSSA architecture and the corresponding sensor-tasking schedule an approximate solution for any epoch. 
Note, however, that this assumption is not a necessary ingredient to the subsequent development of the TE-$p$-MP; if one opts to incorporate $\boldsymbol{l}_{Sk}$ with $z$ components, the TE-$p$-MP solution becomes epoch-specific.

\subsubsection{Pointing Direction}
Instead of allocating observers to look at a specific target, the TE-$p$-MP will optimize for allocations of observers to a particular pointing direction; whether a target is visible to an observer will then depend on whether the target is within the observer's field of view. We consider a discrete set of $m$ pointing directions, indexed by $i = 1, \ldots, m$, which divides the sphere surrounding an observer into $m$ directions, parameterized in terms of azimuth and elevation angles. 

\begin{figure}
    \centering
    \includegraphics[width=0.49\linewidth]{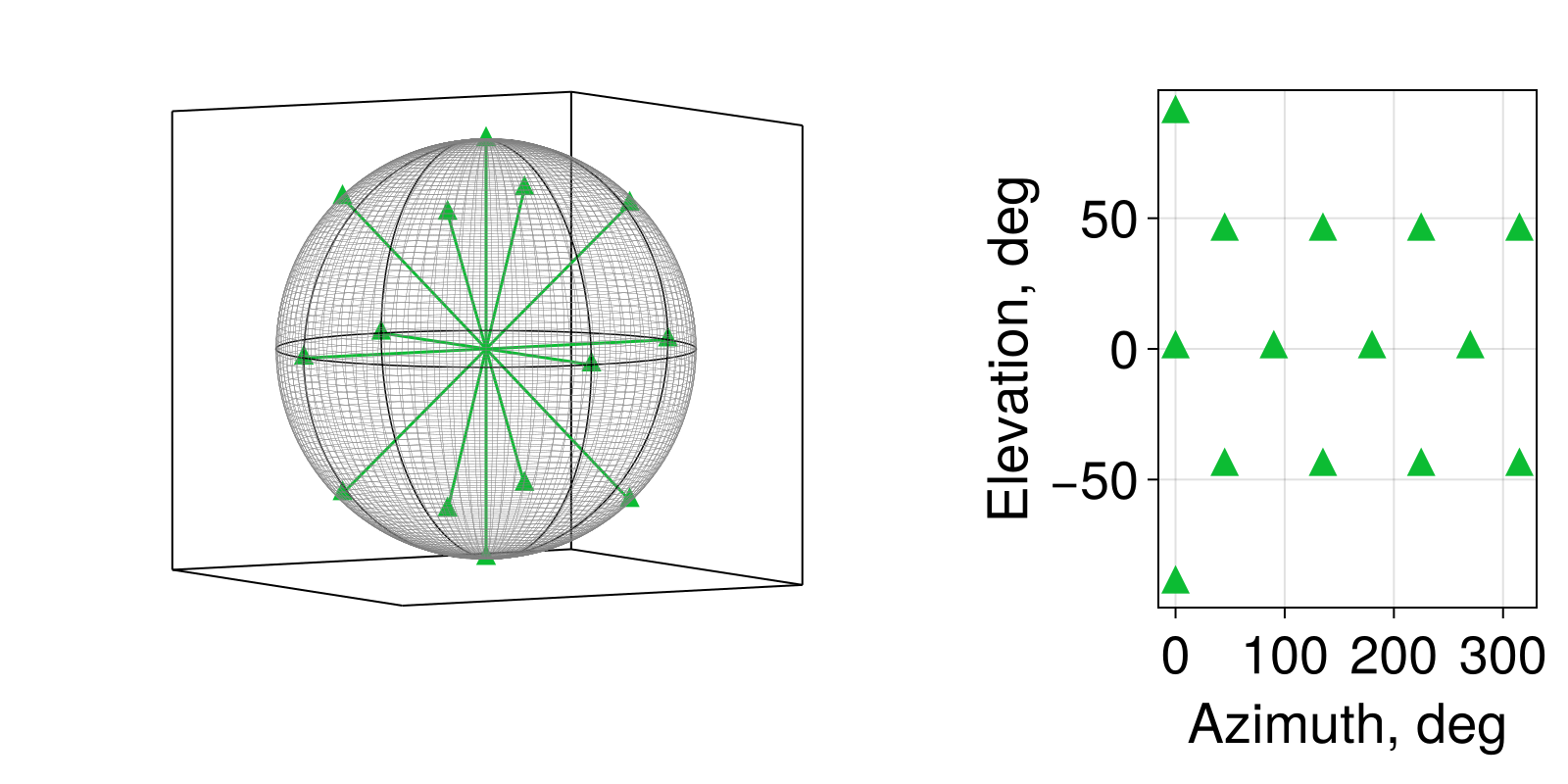}
    \caption{Discretization of pointing directions}
    \label{fig:pointing_directions}
\end{figure}

The problem of dividing a sphere into $m$ equally distributed points is non-trivial, and suboptimal, approximation approaches such as the Fibonacci spiral may be used, especially for large numbers of $m$. 
In this work, because the TE-$p$-MP aims to rapidly explore the impact of different parameters of the SSA architecture, we choose to limit $m$ to a relatively small number with easily interpreted directions. Thus, we define a set of 14 pointing directions, with 6 orthogonal directions, and 8 diagonal directions, as shown in Figure~\ref{fig:pointing_directions}. 
We assume an azimuth and elevation angles of of $0^{\circ}$ points the sensor along the $x$-axis in the Earth-Moon rotating frame. 

\subsubsection{Visibility boolean}
\begin{figure}
    \centering
    \includegraphics[width=0.65\linewidth]{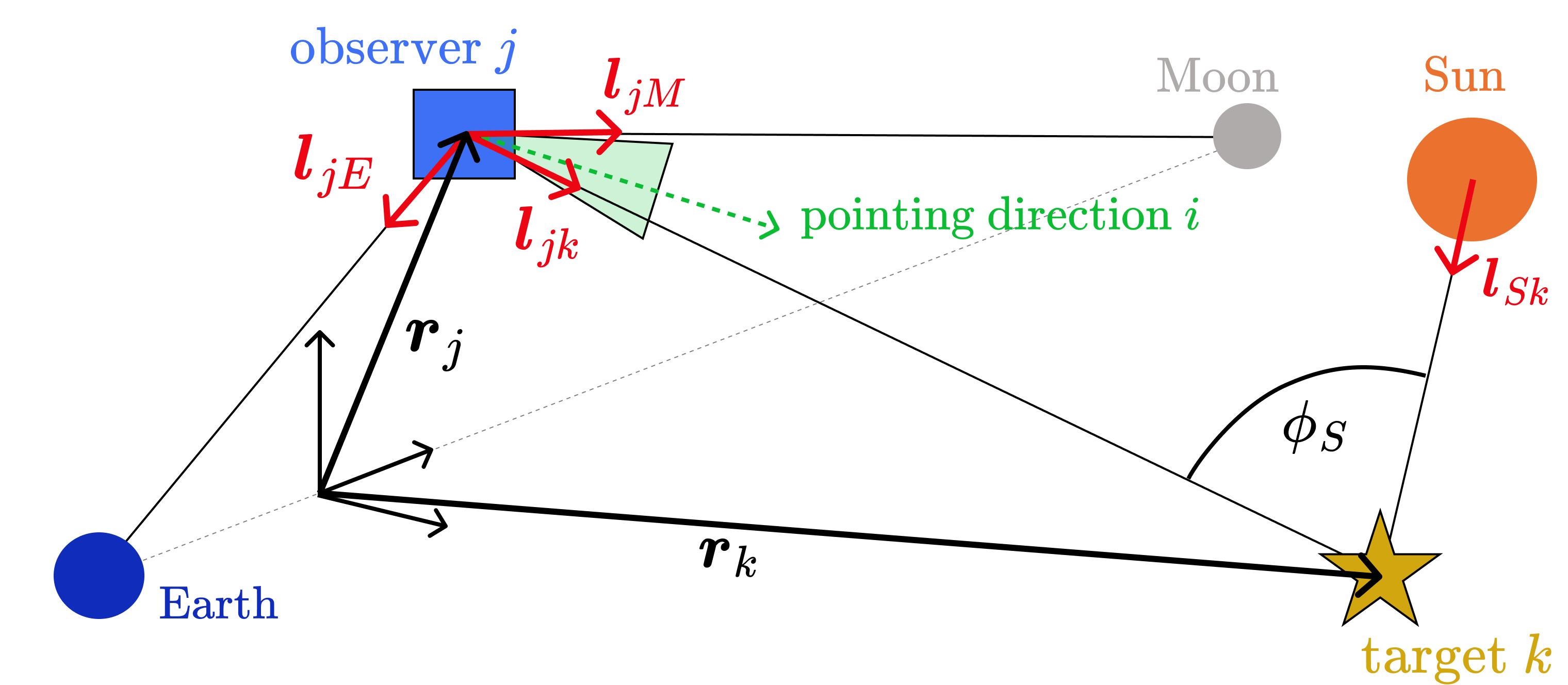}
    \caption{Discretizations used to compute the visibility boolean of a target $k$ from an observer location $j$
    }
    \label{fig:lflp_ijtk_definitions_cr3bp}
\end{figure}

We seek to define a binary parameter that dictates whether an observation of a target is possible and is of sufficient quality. 
We make use of the target's apparent magnitude~\eqref{eq:apparent_magnitude} to define the \textit{visibility boolean}, denoted by $\bar{m}^b \in \mathrm{B}$; the use of apparent magnitude as a metric for assessing whether a target is observable can be found in previous literature, for example in~\cite{Vendl2021,Visonneau2023}. 
Alternatively, one may also incorporate more advanced observation metrics such as the signal-to-noise ratio (SNR)~\cite{Fowler2020,Frueh2021,Badura2022,Klonowski2023} to define the visibility boolean; since the focus of this work is to propose the algorithmic framework for formulating and solving the CSSA constellation design problem, the incorporation of such metrics are deemed beyond the scope of this work. 
Four conditions must be met for an observer to be visible: the first and second conditions ensure the Moon and the Earth are not ahead or behind the target; the third condition ensures the target to lie in the FOV of the observer; the last condition ensures the target's apparent magnitude $\bar{m}_{\mathrm{target}} (\rbold_j, \rbold_k)$ to be below a threshold visible apparent magnitude $\bar{m}_{\mathrm{crit}}$. 
Taking these conditions into account, for a given observer position $\rbold_j(t)$, target position $\rbold_k(t)$, and sensor direction $i(t)$, $\bar{m}^b$ is given by
\begin{equation}
\begin{aligned}
    \bar{m}^b (\rbold_j(t), \rbold_k(t), i(t))
    &=
    \begin{cases}
        0 & \cos^{-1} \left(\boldsymbol{l}_{jM}(t) \cdot \boldsymbol{l}_{jk}(t) \right)
            < \phi_{\mathrm{Moon}}^{\mathrm{crit}}(\rbold_j)
        \\
        0 & \cos^{-1} \left(\boldsymbol{l}_{jE}(t) \cdot \boldsymbol{l}_{jk}(t) \right)
            < \phi_{\mathrm{Earth}}^{\mathrm{crit}}(\rbold_j)
        \\
        0 & \rbold_k(t) \notin \operatorname{FOVCone}_i\left(\rbold_j (t)\right)
        \\
        0 & \bar{m}_{\mathrm{target}} > \bar{m}_{\mathrm{crit}}
        \\
        1 & \text{otherwise}
    \end{cases}
\end{aligned}
\label{eq:visibility_boolean}
\end{equation}
where $\boldsymbol{l}_{jk}$ is the line of sight vector from slot $j$ to target $k$, $\boldsymbol{l}_{jM}$ is the line-of-sight vector from slot $j$ to the Moon, and $\boldsymbol{l}_{jE}$ is the line-of-sight vector from slot $j$ to the Earth. 
In equation~\eqref{eq:visibility_boolean}, {the first condition ensures the angular separation between $\boldsymbol{l}_{jk}$ and $\boldsymbol{l}_{jM}$ is separated by the lunar apparent angular radius $\phi_{\mathrm{Moon}}^{\mathrm{crit}}(\rbold_j)$; similarly, the second condition ensures the angular separation between $\boldsymbol{l}_{jk}$ and $\boldsymbol{l}_{jE}$ is separated by the Earth's apparent angular radius $\phi_{\mathrm{Earth}}^{\mathrm{crit}}(\rbold_j)$; 
}
the third condition ensures the target $\rbold_k$ is within the cone formed by the FOV from $\rbold_j$ pointed direction $i$; {finally, the fourth condition ensures the target appears sufficiently bright to the observer.}
Figure~\ref{fig:lflp_ijtk_definitions_cr3bp} shows the definitions of the quantities defined to compute the visibility boolean. 

\section{Time-Expanded $p$-Median Problem}
\label{sec:TEpMP}
The traditional $p$-Median problem considers the problem of choosing $p$ facilities from $n$ candidate slots to allocate services to $m$ clients. 
While the mathematical formulation developed in this work shares some resemblance to the problem in terrestrial applications, we note that some of the physical representations of variables are different. 

We consider $j = 1,\ldots,n$ candidate observer locations, as is the case in the terrestrial problem. The index $i = 1\ldots,m$ denotes the candidate pointing directions for each observer to point its optical sensor. 
The targets are indexed by $k = 1\ldots,q$. Since the location of the observers, illumination conditions, and (depending on the considered problem) the position of the targets all move in time, an additional time index $t = 1\ldots,\ell$ is introduced. 
Note that, while in terrestrial FLPs, the ``location'' of a facility is its physical position, the ``location'' in space-based FLPs defines a ``slot'' along an orbit in which the observer spacecraft may be placed. Hence, a specific observer ``location'' does not refer to a fixed coordinate in space, but rather this abstract ``slot'' in which the observer may be placed. 
We use the expression ``position'' to denote the physical, Cartesian coordinates of the observer at a given time.

\subsection{Design Variables}
The TE-$p$-MP has two sets of binary variables, $X \in \mathbb{B}^{m \times n \times \ell}$ and $Y \in \mathbb{B}^{n}$. 
In traditional FLPs, the variable $Y$ dictates whether a candidate facility location is used in the solution, while the variable $X$ dictates the allocation of existing facilities to clients. 
In the context of the CSSA problem, this allocation can be understood as an \textit{assignment} of an observer to monitor in a particular direction; furthermore, due to the dynamic nature of space, allocation decisions must be taken at each time step. 
The variables $Y$ takes the same definition as~\eqref{eq:pMP_Y_meaning}, while variables $X$ correspond to
\begin{align}
    X_{ijt} &= 
    \begin{cases}
        1 & \text{observer in location $j$ observes along direction $i$ at time $t$}
        \\
        0 & \text{otherwise}
    \end{cases}
\end{align}

\subsection{Problem Formulation}
We formulate a problem to design a CSSA constellation that maximizes the aggregated coverage of all targets across all time steps.
We collect the visibility boolean $\bar{m}^b (\rbold_j(t), \rbold_k(t), i)$ from equation~\eqref{eq:visibility_boolean} corresponding to observer location $j$, pointing direction $i$, target $k$, at time $t$ into a visibility tensor $M^b \in \mathbb{B}^{m \times n \times \ell \times q}$, 
\begin{equation}
    M^b_{ijtk} = 
    \begin{cases}
        \bar{m}^b (\rbold_j(t), \rbold_k(t), i) & D_{tk} > 0
        \\
        0 & D_{tk} = 0
    \end{cases}
    \label{eq:Mb_ijtk_definition}
\end{equation}
where the condition is based on whether the demand matrix $D$ is non-zero for target $k$ at time $t$.
Then, the TE-$p$-MP is given by
\begin{subequations}
\label{eq:TEpMPAggregate}
\begin{align}
    \maximize{X,Y} \quad& \sum_{t=1}^{\ell} \sum_{k=1}^q 
    \left(\max_{\substack{i=1,\ldots,m,\\j=1,\ldots,n}} M^b_{ijtk} X_{ijt} \right) 
    - \dfrac{1}{\ell} \sum_{j=1}^n f_j Y_j
    \label{eq:TEpMP_objective}
    \\
    \suchthat \quad
    & \sum_{i=1}^m X_{ijt} \leq 1 \quad \forall j,t
    \label{eq:TEpMP_constraint_single_dir}
    \\
    & \sum_{j=1}^n Y_j = p
    \label{eq:TEpMP_constraint_pfacility}
    \\
    & X_{ijt} \leq Y_j \quad \forall i,j,t
    \label{eq:TEpMP_constraint_existence}
    \\
    & X_{ijt}, Y_j \in \{0,1\} \quad \forall i,j,t
    \label{eq:TEpMP_constraint_XYbin}
\end{align}
\end{subequations}
The objective~\eqref{eq:TEpMP_objective} has two terms; {the primary objective is the first term}, which \minorreview{corresponds to the number of observed targets across all time steps.}
Note that the first term in~\eqref{eq:TEpMP_objective} is different from simply taking $\sum_{t=1}^{\ell} \sum_{k=1}^q \sum_{i=1}^{m} \sum_{j=1}^n M^b_{ijtk} X_{ijt}$, which would allow a solution to capitalize on observing the same target by multiple observers rather than observing as many targets as possible. 
\minorreview{The second term is a discount factor that penalizes the use of ``expensive'' observer locations in the case that multiple combinations of facilities yield the same value from the first term.}
The vector $f \in \mathbb{R}^n$, defined by
\begin{equation}
    f_j = 1 - \dfrac{1}{\nu_j + 10}
    \label{eq:f_j_definition}
\end{equation}
models the cost of a given facility location $j$, where $\nu_j$ is the stability index of the LPO, as defined by equation~\eqref{eq:stability_index}. 
{
The definition of $f_j$ with equation~\eqref{eq:f_j_definition} allows the ``cost'' of a given facility location to \minorreview{be distributed around similar values, all under $1$}: for the considered set of LPOs, $\nu_j$ varies between $1$ and $1399.19$, while $f_j$ varies between about $0.9091$ and $0.9993$.}
\minorreview{Directly mapping $f_j = \nu_j$ would be inappropriate for the objective~\eqref{eq:TEpMP_objective}, as the ``cost'' associated with placing a single observer on an unstable LPO, say of $\nu_j = 1399.19$, would be approximately the same as placing $1399$ observers on a stable LPO of $\nu_{j^{\prime}} = 1$, which is clearly unrealistic.
In effect, by defining $f_j \in [0,1)$, we ensure the unstable LPOs, some which are known be advantageous for CSSA~\cite{Vendl2021,Visonneau2023}, are not discouraged if a better target coverage can be achieved.}
Since prohibitively expensive candidate observer locations can be excluded from the candidate observer location set $\mathcal{J}$ to begin with, it is assumed that all candidate LPOs {have acceptable station-keeping cost should they be chosen as part of the constellation.}
\minorreview{To ensure the optimal solution never compromises on the number of observed targets to choose a ``cheaper'' constellation according to costs $f_j$, the second term is multiplied by $1/\ell$. Since $f_j < 1 \, \forall j$, $\sum_{j=1}^n f_j Y_j < p$ due to~\eqref{eq:TEpMP_constraint_pfacility}. Also, for a typical CSSA problem, $p << \ell$, so $(1/\ell) \sum_{j=1}^n f_j Y_j < p/\ell << 1$.
Thus, an additional target observation at any single time-step is always prioritized in an optimal solution to~\eqref{eq:TEpMPAggregate}.}

Constraint~\eqref{eq:TEpMP_constraint_single_dir} ensures each observer location $j$ only looks to up to a single direction $i$ at any given time $t$.
Constraint~\eqref{eq:TEpMP_constraint_pfacility} is the traditional $p$-median constraint, which ensures exactly $p$ observer is placed.
\minorreview{As mentioned in Section~\ref{sec:pMP_introduction}, the use of an exact equality constraint in~\eqref{eq:TEpMP_constraint_pfacility} is due to the assumption that for cislunar SSA where the coverage demand is high, having an additional observer can always improve the objective~\eqref{eq:TEpMP_objective}.}
Constraint~\eqref{eq:TEpMP_constraint_existence} is also a common constraint in the traditional $p$-median problem, enforcing the $j^{\mathrm{th}}$ observer to exist for its allocation to view direction $i$ at time $t$ to be used in the solution. 
Finally, constraint~\eqref{eq:TEpMP_constraint_XYbin} ensures the variables $X$ and $Y$ are binary. 

The max operator in the objective~\eqref{eq:TEpMP_objective} may be interpreted as the infinity norm of the vector formed by concatenating all columns of the matrix $M^b_{_:,:,tk} \odot X_{:,:,t}$, where $\odot$ denotes element-wise multiplication. Then, through the introduction of a slack variable $\theta_{tk} \in [0,1]$ for each combination of $t$ and $k$, problem~\eqref{eq:TEpMPAggregate} may be rewritten as 
\begin{subequations}
\label{eq:TEpMP_slack}
\begin{align}
    \maximize{X,Y,\theta} \quad& \sum_{t=1}^{\ell} \sum_{k=1}^q 
    \theta_{tk} 
    - \dfrac{1}{\ell} \sum_{j=1}^n f_j Y_j
    \label{eq:TEpMP_slack_objective}
    \\
    \suchthat \quad
    & \sum_{i=1}^{m} \sum_{j=1}^n M^b_{ijtk} X_{ijt} \geq \theta_{tk} \quad \forall t,k
    \label{eq:TEpMP_slack_constraint_infnorm}
    \\
    & 0 \leq \theta_{tk} \leq 1 \quad \forall t,k
    \label{eq:TEpMP_slack_constraint_slackbounds}
    \\ 
    & \text{\eqref{eq:TEpMP_constraint_single_dir} $\sim$ \eqref{eq:TEpMP_constraint_XYbin}}
    \nonumber
\end{align}
\end{subequations}
where constraint~\eqref{eq:TEpMP_slack_constraint_infnorm} ensures this modified problem~\eqref{eq:TEpMP_slack} has the same optimal solution as the original problem~\eqref{eq:TEpMPAggregate}. 
Problem~\eqref{eq:TEpMP_slack} is a standard mixed-integer linear program (MILP), which is amenable to standard B\&B solvers such as Gurobi~\cite{gurobi} or CPLEX~\cite{cplex2009v12}. 
While the slack variable $\theta_{tk}$ is a continuous variable in $[0,1]$, the optimal solution to~\eqref{eq:TEpMP_slack} has $\theta_{tk} = 1$ for time $t$ and target $k$ if at least one observer is viewing $k$ at this time, and $\theta_{tk} = 0$ if $\sum_{i=1}^m \sum_{j=1}^n M^b_{ijtk} X_{ijt} = 0$. 
Thus, this slack variable may be interpreted as a binary variable, such that
\begin{equation}
    \theta_{tk}
    = 
    \begin{cases}
        1 & \text{target $k$ is observed by at least one observer at time $t$} \\
        0 & \text{otherwise}
    \end{cases}
\end{equation}
For convenience in comparing obtained solutions later, we define the fraction of observation demand that has been met across all time steps, denoted by $\Theta$, and given by
\begin{equation}
    \Theta \triangleq \frac{\sum_{t=1}^{\ell} \sum_{k=1}^q \theta_{tk}}{\sum_{t=1}^{\ell} \sum_{k=1}^q D_{tk}}
    \label{eq:observation_fraction}
\end{equation}
The formulation~\eqref{eq:TEpMP_slack} enjoys a few favorable properties; firstly, even though the slack variables $\theta$ increase the decision space, it is a continuous variable, mending itself well to standard linear program techniques. 
Secondly, despite the introduction of the time index $t$ and the decomposition of allocation into pointing directions $i$ and targets $k$, the general structure of the traditional $p$-median problem is preserved; this enables efficient use of Lagrangian relaxation, which will be introduced in Section~\ref{sec:lagrangean_method}.

\section{Lagrangian Method for the Time-Expanded $p$-Median Problems}
\label{sec:lagrangean_method}
{We develop a tailored version of the LM for the TE-$p$-MP.
For traditional $p$-Median problems, the \textit{complicating} constraints is typically those that pertain to allocations between $i$ and $j$, resulting in decoupled subproblems corresponding to each candidate facility location $j$.}
LM is particularly advantageous for the $p$-Median problem as the resulting relaxed subproblems can be solved analytically. 
Despite the introduction of the time indices $t$ and slack variables $\theta_{tk}$, we are also able to arrive at a relaxed problem that achieve decomposition into subproblems corresponding to each $j$.
As a result, the upper bound solution~$Z_{\rm relax}$ can be obtained analytically for a given set of Lagrange multipliers. 
We further provide a set of customized heuristics to generate a competitive feasible lower bound solution~$Z_{\rm feas}$. 

{The remainder of this Section is organized as follows: first,} Section~\ref{sec:upperbound_lagrelax} discusses the relaxed problem formulation, along with the approach to decompose the resulting relaxed problem into facility-wise subproblems and to solve each of these subproblems. 
This is followed by Section~\ref{sec:lowerbound_heuristics} which describes the heuristic scheme devised to generate a feasible solution based on the relaxed solution at each iteration. 
An additional set of heuristics for performing a neighborhood search in an attempt to improve this feasible solution is described in Section~\ref{sec:heuristic_neighborhood}. 
In Section~\ref{sec:subgradient_optimization}, the subgradient method used to update the multipliers is introduced. 
The processes of the LM are summarized in Figure~\ref{fig:lagrangean_method_flow}. 
Finally, in Section~\ref{sec:conding_considerations}, we provide implementation considerations for implementing the proposed LM. 

\begin{figure}
    \centering
    \begin{tikzpicture}[
        squarednode/.style={
            rectangle,
            text centered,
            draw=black,
            minimum size=5mm},
        startstop/.style={rectangle, rounded corners, 
            minimum width=1.4cm, 
            minimum height=0.8cm,
            text centered, 
            draw=black,
            fill=gray!10},
        pregunda/.style={diamond, 
            aspect=2.5,
            minimum width=3cm, 
            minimum height=0.5cm,
            text centered, 
            draw=black,
            fill=gray!10}
    ]
        \node[startstop] (start) {Start};
        \node[squarednode] (lagrelax) [below=of start, yshift=0.5cm, align=center] {Compute solution to relaxed problem analytically\\to obtain upper bound $Z_{\rm relax}$ (Section~\ref{sec:upperbound_lagrelax})};
        \node[squarednode] (turnfeasible) [below=of lagrelax, yshift=0.5cm, align=center] {Convert allocation to obtain\\feasible lower bound $Z_{\rm feas}$ (Section~\ref{sec:lowerbound_heuristics})};
        
        \node[squarednode] (neighbor) [below=of turnfeasible, yshift=0.4cm, align=center] {Perform heuristic neighborhood\\search to improve $Z_{\rm feas}$ (Section~\ref{sec:heuristic_neighborhood})};
        
        \node[pregunda] (convergence) [below=of neighbor, yshift=0.5cm] {Converged?};
        \node[startstop] (stop) [right=of convergence, xshift=-0.0cm]  {Stop};
        
        \node[squarednode] (updatemultipliers) [below=of convergence, yshift=0.5cm] {Update Lagrange multipliers (Section~\ref{sec:subgradient_optimization})};
        
        \draw[->] (start.south) -- (lagrelax.north);
        \draw[->] (lagrelax.south) -- (turnfeasible.north);
        \draw[->] (turnfeasible.south) -- (neighbor.north);
        \draw[->] (neighbor.south) -- (convergence.north);
        
        \draw[->] (convergence.south) -- node[anchor=west] {no} (updatemultipliers.north);
        \draw[->] (convergence.east) -- node[anchor=south] {yes} (stop.west);
        \coordinate (leftUB) at ($(lagrelax.west) - (0.5cm,0)$);
        \draw[->]  
            (updatemultipliers.west) -| (leftUB) -- (lagrelax.west);
        
    \end{tikzpicture}
    \caption{Overview of Lagrangian method for the time-expanded $p$-Median problem}
    \label{fig:lagrangean_method_flow}
\end{figure}

\subsection{Upper Bound Solution via Lagrangian Relaxation}
\label{sec:upperbound_lagrelax}
{We first seek to find the decision variables and the corresponding upper-bound solution to the relaxed problem.}
Relaxing constraints~\eqref{eq:TEpMP_constraint_single_dir} and~\eqref{eq:TEpMP_slack_constraint_infnorm} with multipliers $\lambda_{jt} \geq 0$ and $\eta_{tk} \geq 0$ respectively, the relaxed problem is given by
\begin{equation}
    \begin{aligned}
        \maximize{X,Y,\theta} \quad& \sum_{t=1}^{\ell} \sum_{k=1}^q 
        \theta_{tk} 
        - \dfrac{1}{\ell} \sum_{j=1}^n f_j Y_j
        - \sum_{j=1}^n \sum_{t=1}^{\ell} 
        \lambda_{jt} 
        \left(
            \sum_{i=1}^m X_{ijt} - 1
        \right)
        - \sum_{t=1}^{\ell} \sum_{k=1}^q
        \eta_{tk} 
        \left( 
            \theta_{tk} - \sum_{i=1}^{m} \sum_{j=1}^n M^b_{ijtk} X_{ijt} 
        \right)
        \\
        \suchthat \quad
        & \text{\eqref{eq:TEpMP_constraint_pfacility} $\sim$ \eqref{eq:TEpMP_constraint_XYbin},~\eqref{eq:TEpMP_slack_constraint_slackbounds}}
    \end{aligned}
    \label{eq:TEpMP_slack_lagrelax}
\end{equation}
The objective may be reorganized by regrouping variables to
\begin{equation}
    \maximize{X,Y,\theta} \quad
    \sum_{t=1}^{\ell} \sum_{k=1}^q 
        \left(1 - \eta_{tk}\right) \theta_{tk}
    + \sum_{i=1}^m \sum_{j=1}^n \sum_{t=1}^{\ell}
        \left(
            \sum_{k=1}^q 
                \eta_{tk}
                M^b_{ijtk}
            - 
            \lambda_{jt} 
        \right)
        X_{ijt}
    - \dfrac{1}{\ell} \sum_{j=1}^n f_j Y_j
    + \sum_{j=1}^n \sum_{t=1}^{\ell}
        \lambda_{jt}
\end{equation}
The sought optimal decision variables to the relaxed problem~\eqref{eq:TEpMP_slack_lagrelax}, denoted by $\bar{X}$, $\bar{Y}$, and $\bar{\theta}$ generates an upper bound~$Z_{\rm relax}$ to problem~\eqref{eq:TEpMP_slack}. 

We first turn our attention to the slack variables $\theta_{tk}$. Due to the relaxation of constraint~\eqref{eq:TEpMP_slack_constraint_infnorm}, the only remaining constraint on $\theta_{tk}$ is the box bound constraint~\eqref{eq:TEpMP_slack_constraint_slackbounds}. 
Therefore, the optimal value $\bar{\theta}_{tk}$ is simply driven by its coefficient in the objective, and is trivially given by
\begin{equation}    \label{eq:LR_theta}
    \bar{\theta}_{tk}
    = 
    \begin{cases}
        1 & 1 - \eta_{tk} > 0
        \\
        0 & \text{otherwise}
    \end{cases}
\end{equation}

\subsubsection{Decomposition into Subproblems}
We can further recognize that this relaxation decouples the problem on a per-candidate observer location basis.
For each candidate location $j$, the subproblem is given by
\begin{subequations}    \label{eq:TEpMP_slack_lagrelax_subproblem}
\begin{align}
    \maximize{X,Y_j,\theta} \quad& 
    \sum_{i=1}^m \sum_{t=1}^{\ell}
        \left(
            \sum_{k=1}^q 
                \eta_{tk}
                M^b_{ijtk}
            -
            \lambda_{jt} 
        \right)
        X_{ijt}
    - \dfrac{1}{\ell} \left( f_j Y_j \right)
    + \sum_{t=1}^{\ell}
        \lambda_{jt}
    \label{eq:TEpMP_slack_lagrelax_subproblem_objective}
    \\
    \suchthat \quad
    & X_{ijt} \leq Y_j \quad \forall i,t
    \label{eq:TEpMP_slack_lagrelax_subproblem_constraint_existence}
    \\
    & X_{ijt}, Y_j \in \{0,1\} \quad \forall i,t
\end{align}
\end{subequations}
Subproblem~\eqref{eq:TEpMP_slack_lagrelax_subproblem} has as variables $X_{ijt} \in \mathbb{B}^{m \times 1 \times \ell}$, $Y_j \in \mathbb{B}$, and $\theta \in \mathbb{R}^{\ell \times k}$. 
We denote the decision on $X$ and $Y_j$ made with the $j^{\rm th}$ subproblem as $\bar{X}^j$ and $\bar{Y}_j^j$, representing the allocations of the $j^{\rm th}$ observer to pointing directions over time, as well as whether the $j^{\rm th}$ observer should be used, respectively.

We consider separately the cases where $\bar{Y}^{j}_j = 0$ or $\bar{Y}^{j}_j = 1$; {the decision on whether $\bar{Y}^{j}_j = 0$ or $\bar{Y}^{j}_j = 1$ is determined at a later stage, when the relaxed solution is constructed by assembling solutions from each subproblem.}
If $\bar{Y}^{j}_j = 0$, constraint~\eqref{eq:TEpMP_slack_lagrelax_subproblem_constraint_existence} trivially forces $\bar{X}^{j}_{ijt} = 0$ $\forall i,t$. 
If $\bar{Y}^{j}_j = 1$, the choice of each $\bar{X}^{j}_{ijt}$ $\forall i,t$ is driven by its coefficient in the objective, and thus the optimal value $\bar{X}_{ijt}^{j}$ is given by
\begin{equation}
    \bar{X}_{ijt}^{j}
    =
    \begin{cases}
        1 & \sum_{k=1}^q \eta_{tk} M^b_{ijtk} - \lambda_{jt} > 0
        \\
        0 & \text{otherwise}
    \end{cases}
    \label{eq:LR_upperbound_Xijt_j}
\end{equation}
The corresponding objective of the $j^{\mathrm{th}}$ subproblem~\eqref{eq:TEpMP_slack_lagrelax_subproblem_objective} assuming $\bar{Y}_j = 1$, denoted by $Z^{j}_{\mathrm{relax}}$, is given by
\begin{equation}
    Z^{j}_{\mathrm{relax}} = 
    \sum_{i=1}^m \sum_{t=1}^{\ell}
        \left(
            \sum_{k=1}^q 
                \eta_{tk}
                M^b_{ijtk}
            - 
            \lambda_{jt}
        \right)
        \bar{X}_{ijt}^{j}
    - \dfrac{f_j}{\ell} 
    + \sum_{t=1}^{\ell}
        \lambda_{jt}
    \label{eq:subproblem_objective_YBar1}
\end{equation}

\subsubsection{Construction of Lagrangian Relaxed Solution}
To construct the relaxed solution, we first compute $Z^{j}_{\mathrm{relax}}$ for all $j=1,\ldots,n$ subproblems~\eqref{eq:TEpMP_slack_lagrelax_subproblem} according to~\eqref{eq:subproblem_objective_YBar1}. 
We then define the permutation $\pi$ for $Z_{\rm relax}^j$ of the $j = 1,\ldots,n$ candidate observer locations in descending order 
\begin{equation}    \label{eq:LR_Z_sorting}
    Z^{\pi(1)}_{\mathrm{relax}} \geq Z^{\pi(2)}_{\mathrm{relax}} \geq \ldots \geq Z^{\pi(n)}_{\mathrm{relax}}
\end{equation}
To enforce the $p$-median constraint~\eqref{eq:TEpMP_constraint_pfacility}, variables $\bar{Y}$ are chosen such that the $p$ best facilities are chosen
\begin{equation}    \label{eq:LR_Ybar_definition}
    \bar{Y}_{\pi(j)} =
    \begin{cases}
        1 & j \leq p
        \\
        0 & \text{otherwise}
    \end{cases}
\end{equation}
Having chosen the $p$ facility locations to be used, we construct the upper bound solution $\bar{X}_{ijt}$ by making use of subproblem solutions $\bar{X}_{ijt}^j$ from equation~\eqref{eq:LR_upperbound_Xijt_j}, such that
\begin{equation}   \label{eq:LR_Xbar_definition}
    \bar{X}_{i \pi(j) t} 
    = 
    \begin{cases}
        \bar{X}_{i\pi(j)t}^{\pi(j)} & j \leq p
        \\ 
        0 & \text{otherwise}
    \end{cases}
\end{equation}
Finally, the upper bound objective of the Lagrangian relaxed problem~\eqref{eq:TEpMP_slack_lagrelax} $Z_{\mathrm{relax}}$ is given by substituting the decision variables to equation~\eqref{eq:TEpMP_slack_objective},
\begin{equation}    \label{eq:LR_Z_relax}
    Z_{\mathrm{relax}} = \sum_{t=1}^{\ell} \sum_{k=1}^q 
        \left(1 - \eta_{tk}\right) \bar{\theta}_{tk}
    +
    \sum_{s = 1}^p Z^{\pi(s)}_{\mathrm{relax}}
\end{equation}
Algorithm~\ref{alg:lagrangian_relaxation} summarizes the relaxation process. 
Within the algorithm, \verb|argsortDescend| returns the sorted indices of the inputs in descending order.
In effect, the relaxed solution to the original problem~\eqref{eq:TEpMP_slack} is obtained through only simple algebraic evaluations and without any iterative operation. 

\begin{algorithm}
\caption{Lagrangian Relaxation for TE-$p$-MP}
\label{alg:lagrangian_relaxation}
\begin{algorithmic}[1]
    \Require $\lambda_{jt}$ $\forall j,t$, $\eta_{tk}$ $\forall t,k$  
    \State $\bar{Y}, \bar{X} \gets \boldsymbol{0}_{n \times 1}, \boldsymbol{0}_{m \times n \times \ell}$    \Comment{Initialize}
    \State $\bar{\theta}_{tk} \gets \text{eqn.~\eqref{eq:LR_theta} }\forall t,k$
    \For{$j = 1,\ldots,n$}
        \State $\bar{X}_{ijt}^j \gets \text{eqn.~\eqref{eq:LR_upperbound_Xijt_j} }\forall i,t$
            \Comment{Compute allocation decision variables for subproblem $j$}
        \State $Z_{\rm relax}^j \gets \text{eqn.~\eqref{eq:subproblem_objective_YBar1}}$
            \Comment{Compute objective for subproblem $j$}
    \EndFor
    \State $\pi \gets \verb|argsortDescend|(Z_{\rm relax}^1, \ldots, Z_{\rm relax}^n)$
        \Comment{Sort in descending order according to~\eqref{eq:LR_Z_sorting}}
    \For{$j = \pi(1),\ldots,\pi(p)$}
        \State $\bar{Y}_j \gets 1$
            \Comment{Use facility location $j$}
        \State $\bar{X}_{ijt} \gets \bar{X}_{ijt}^j$
            \Comment{Assign targets according to decision made for subproblem $j$}
    \EndFor 
    \State $Z_{\rm relax} \gets \text{eqn.~\eqref{eq:LR_Z_relax}}$
    \State \Return $\bar{X}$, $\bar{Y}$, $\bar{\theta}$, $Z_{\rm relax}$
\end{algorithmic}
\end{algorithm}

\subsection{Lagrangian Heuristics Reallocation}
\label{sec:lowerbound_heuristics}
If the Lagrangian relaxed solution $(\bar{X}, \bar{Y}, \bar{\theta})$ is feasible, then it is the optimal solution $(X^*, Y^*, \theta^*)$ to problem~\eqref{eq:TEpMP_slack}; in most cases, $(\bar{X}, \bar{Y}, \bar{\theta})$ is infeasible as constraints~\eqref{eq:TEpMP_constraint_single_dir} and~\eqref{eq:TEpMP_slack_constraint_infnorm} have been relaxed. 
Thus, we seek to construct a feasible \textit{lower} bound solution $(\underline{X}, \underline{Y}, \underline{\theta})$ using the relaxed solution. 
The overall procedure is summarized in Algorithm~\ref{alg:feasible_allocation}. 

\begin{algorithm}
\caption{Heuristic Feasible Allocation Construction}
\label{alg:feasible_allocation}
\begin{algorithmic}[1]
    \Require $\bar{X}$, $\mathcal{J}^{\bar{Y}}$, $M^b$  
    \State $\underline{X} \gets \bar{X}$  \Comment{Initialize allocation variable $X$ of feasible solution}
    \State $\underline{\theta} \gets \boldsymbol{0}_{\ell \times q}$   \Comment{Initialize slack variable $\theta$ of feasible solution}
    
    \For{$t = 1, \ldots, \ell$}
    
        \State $\mathcal{K}^0_t \gets \mathcal{K}$  \Comment{Initialize set of unobserved targets at time $t$}
        \State $\mathcal{J}^0_t \gets \varnothing$  \Comment{Initialize set of observers with infeasible allocations}
    
        \For{$j \in \mathcal{J}^{\bar{Y}}$}
            \If{$\sum_{i=0}^m \underline{X}_{ijt} = 1$}
                \State $\mathcal{K}^0_t \gets \mathcal{K}^0_t \setminus k$  \Comment{Remove target $k$ from $\mathcal{K}^0_t$}
            \Else
                \State $\mathcal{J}^0_t \gets \mathcal{J}^0_t \cup \{j\}$  \Comment{Append violating facility index $j$ to $\mathcal{J}^0_t$}
                \State $\underline{X}_{:jt} \gets 0$    \Comment{Turn all allocations of facility $j$ at time $t$ off}
            \EndIf
        \EndFor
        
        \State $\underline{X} \gets $ \verb|GreedyAllocation|$(\underline{X}, M^b, \mathcal{J}^0_t, \mathcal{K}^0_t)$ 
        \Comment{Perform allocations of facilities in $\mathcal{J}^0_t$}

        \For{$t=1,\ldots,\ell, \, k=1,\ldots,q$}
            \If{$\sum_{i=1}^m \sum_{j=1}^n M^b_{ijtk} \underline{X}_{ijt} \geq 1$}
                \State $\underline{\theta}_{tk} \gets 1$   \Comment{If at least one observer is seeing target $k$ at time $t$, turn $\underline{\theta}_{tk}$ on}
            \EndIf  
        \EndFor
    \EndFor
    \State \Return $\underline{X}$, $\underline{\theta}$
\end{algorithmic}
\end{algorithm}

The first step is to fix $\underline{Y} = \bar{Y}$ and initialize $\underline{X} = \bar{X}$ and $\underline{\theta} = \bar{\theta}$. Furthermore, let $\mathcal{J}^{\bar{Y}}$ denote the set of indices $j$ that is active in $\bar{Y}$, 
\begin{equation}
    \mathcal{J}^{\bar{Y}} = 
    \left\{
        j \in 1, \ldots, n \mid Y_j = 1
    \right\}
\end{equation}
Then, we seek to satisfy constraint~\eqref{eq:TEpMP_constraint_single_dir}. 
We begin by fixing allocations of observer location $j$ to direction $i$ for each time $t$ that satisfies the constraint with an equality,
\begin{equation}
    \underline{X}_{ijt} = \begin{cases}
        1 & \sum_{i=1}^m \bar{X}_{ijt} = 1
        \\
        0 & \text{otherwise}
    \end{cases}
\end{equation}
and we construct the set of all targets that are unobserved at time $t$, denoted by  $\mathcal{K}^0_t$
\begin{equation}
    \mathcal{K}^0_t = 
    \left\{
        k \in \mathcal{K} 
        \mid \sum_{j=1}^n \sum_{i=1}^m \underline{X}_{ijt} < 1
    \right\}
\end{equation}
We also record facility indices $j$ for which the constraint at time $t$ is violated into set $\mathcal{J}_t^0$, given by
\begin{equation}
    \mathcal{J}_t^0 =
    \left\{
        j \in \mathcal{J}^{\bar{Y}} \mid \sum_{i=1}^m X_{ijt} \neq 1
    \right\}
\end{equation}
Equivalently, this set may be understood as the set of $j$ that still requires allocation to an appropriate pointing direction $i$. 
The next step consists of adding back allocations for $j \in \mathcal{J}_t^0$ for each time-step $t$, which corresponds to the \verb|GreedyAllocation()| functions call in line 14 of Algorithm~\ref{alg:feasible_allocation}. 
{We propose two approaches, coined the \textit{greedy allocation} and \textit{full-factorial allocation}, respectively. The choice between these two algorithms is made based on the number of observers $p$ of a given TE-$p$-MP instance and the relative computational time the user chooses to dedicate to the heuristic search.
The two allocation algorithms are introduced in further detail in subsections~\ref{sec:greedy_allocation} and~\ref{sec:fullfactorial_allocation}.}
Finally, we fix the solution to satisfy constraint~\eqref{eq:TEpMP_slack_constraint_infnorm} by iterating through $t$ and setting $\underline{\theta}_{tk} = 0$ if 
\begin{equation}
    \underline{\theta}_{tk}
    =
    \begin{cases}
        0 & \sum_{i=1}^m \sum_{j=1}^n M^b_{ijtk} X_{ijt} = 0
        \\
        1 & \text{otherwise}
    \end{cases}
    \label{eq:lagrel_fix_theta_tk}
\end{equation}
Finally, the lower bound feasible solution $Z_{\mathrm{feas}}$ is given by
\begin{equation}
    Z_{\mathrm{feas}} = \sum_{t=1}^{\ell} \sum_{k=1}^q 
        \left(1 - \eta_{tk}\right) \underline{\theta}_{tk}
    + \sum_{j=1}^n f_j \underline{Y}_j
\end{equation}

\subsubsection{Greedy Allocation}
\label{sec:greedy_allocation}
{The first allocation strategy we propose is a greedy allocation.}
With the greedy allocation, the facility $j^{\star}$ with pointing direction $i^{\star}$ that results in the largest additional number of targets observed is selected via the expression
\begin{equation}
    (i^{\star}, j^{\star}) = \underset{\substack{i = 1,\ldots,m
        \\j \in \mathcal{J}_t^0}
    }{\operatorname{argmax}}\sum_{k \in \mathcal{K}^0_t} M_{ijtk}
    \label{eq:lagrel_greedy_ijstar}
\end{equation}
This process is repeated, each time setting $X_{i^{\star}j^{\star}t} = 1$ and removing $j^{\star}$ from $\mathcal{J}_t^0$, until $\mathcal{J}_t^0 = \varnothing$. 
{
Pseudo-code for the greedy allocation scheme is provided in Algorithm~\ref{alg:greedy_allocation}. 
Since the algorithm essentially involves summation across $k$ and sorting across $i$ and $j$, its computational time is insignificant with dimensions of typical TE-$p$-MP instances.}

\begin{algorithm}
\caption{Greedy Allocation}
\label{alg:greedy_allocation}
\begin{algorithmic}[1]
    \Require $t, \underline{X}$, $M^b$, $\mathcal{J}^0_t$, $\mathcal{K}^0_t$
    
    \While{$\mathcal{J}^0_t \neq \varnothing \wedge \mathcal{K}^0_t \neq \varnothing$}
        \State $i^{\star}, j^{\star} \gets $ eqn.~\eqref{eq:lagrel_greedy_ijstar}
        \State $\underline{X}_{i^{\star} j^{\star} t} \gets 1$    \Comment{Turn allocation on}
        \State $\mathcal{J}^0_t \gets \mathcal{J}^0_t \setminus j^{\star}$ \Comment{Remove $j^{\star}$ from set of observers with infeasible allocations}
        \State $\mathcal{K}^0_t \gets \mathcal{K}^0_t \setminus \{k \mid M^b_{i^{\star}j^{\star}tk} = 1\}$ \Comment{Remove newly targets observed from set of unobserved targets}
    \EndWhile
    
    \State \Return $\underline{X}$
\end{algorithmic}
\end{algorithm}

\subsubsection{Full-Factorial Allocation}
\label{sec:fullfactorial_allocation}
Typically, the number of observers $p$ in the context of the CSSA problem is relatively small, on the orders of 1s to 10s. 
Provided that $p$ is small enough, it may be computationally reasonable to spend additional time sequentially allocating remaining observers $j \in \mathcal{J}^0_t$ to remaining targets $\mathcal{K}^0_t$. 
We refer to this process as the \textit{full-factorial} allocation; Algorithm~\ref{alg:fullfactorial_allocation} shows the pseudo-code for its procedure. 
The while-loop from the greedy Algorithm~\ref{alg:greedy_allocation} is now replaced by the outer for-loop in line 2, which iterates through all possible permutations of indices of the set $\mathcal{J}^0_t$, denoted by $\Pi_{\mathcal{J}}$. 
For each permutation $\pi_{\mathcal{J}} \in \Pi_{\mathcal{J}}$, we go through indices $j$ in the order of $\pi_{\mathcal{J}}$, each time allocating $j$ to pointing direction $i$ in a greedy manner, such that 
\begin{equation}
    i^{\star} = \underset{\substack{i = 1,\ldots,m}
    }{\operatorname{argmax}}\sum_{k \in \mathcal{K}^0_t} M^b_{ijtk}
    \label{eq:lagrel_fullfactorialgreedy_istar}
\end{equation}
and targets observed by the allocation $(i^{\star},j)$ are removed from $\mathcal{K}^0_t$. 
Once the successive greedy allocation has been tried for all permutations in $\Pi_{\mathcal{J}}$, the allocations that resulted in the largest additional targets observation is returned. 

\begin{algorithm}
\caption{Full-Factorial Allocation}
\label{alg:fullfactorial_allocation}
\begin{algorithmic}[1]
    \Require $t, \underline{X}$, $M^b$, $\mathcal{J}^0_t$, $\mathcal{K}^0_t$

    \State $\Pi_\mathcal{J} \gets $\verb|permutations|$(\mathcal{J}^0_t)$
        \Comment{Compute permutations of $j$ in $\mathcal{J}^0_t$}
        
    \For{$\pi_\mathcal{J} \in \Pi_{\mathcal{J}}$}    \Comment{Iterate through each permutation}
        
        \State $\mathcal{K}^{0 \pi}_t \gets \mathcal{K}^0_t$  \Comment{Make a copy of $\mathcal{K}^0_t$ for current permutation}
        \State $\underline{X}^{\pi} \gets \underline{X}$  \Comment{Make a copy of $\underline{X}$ for current permutation}
        
        \For{$j \in \pi_{\mathcal{J}}$}
            \If{$\mathcal{K}^{0 \pi}_t = \varnothing$}
                \State break  \Comment{If all targets are observed, no more allocations are necessary}
            \EndIf 

            \State $i^{\star} \gets$ eqn.~\eqref{eq:lagrel_fullfactorialgreedy_istar}
            \State $\underline{X}^{\pi}_{i^{\star}jt} \gets 1$

            \For{$k = 1,\ldots,q$}
                \If{$M^b_{i^{\star}jtk}== 1$}
                    \State $\mathcal{K}^{0 \pi}_t \gets \mathcal{K}^0_t \setminus \{k \mid M^b_{i^{\star}j^{\star}tk} = 1\}$ \Comment{Remove newly targets observed from set of unobserved targets}
                \EndIf
            \EndFor
        \EndFor
        
        \State $\Delta^{\pi} \gets |\mathcal{K}^0_t| - |\mathcal{K}^{0 \pi}_t|$  \Comment{Number of additional targets observed by current permutation}
    \EndFor
    
    \State $\pi^* \gets \min_{\pi} \Delta^{\pi}$  \Comment{Choose allocation best permutation}
    \State $\underline{X} \gets \underline{X}^{\pi^*}$  \Comment{Take $\underline{X}$ corresponding to best permutation}
    
    \State \Return $\underline{X}$
\end{algorithmic}
\end{algorithm}

\subsection{Heuristics for Neighborhood Facility Swap}
\label{sec:heuristic_neighborhood}
The major limitation of the relaxation procedure is the lack of convergence guarantee; depending on, among other things, the choice of constraints that are relaxed and the construction of the heuristic lower bound solution, it is possible that the procedure is no longer able to reduce the optimality gap. For the problem in question, preliminary experiments with the greedy/full-factorial allocation heuristics alone have revealed that while the algorithm performs well in most cases, there are instances where it gets stuck with intuitively suboptimal choice of observer locations. We thus develop additional heuristics to mutate the {chosen set of facilities.}

Consider a list of sets $\mathcal{C}$ of length $n$, where each entry is a set of neighboring facility indices to each corresponding facility $j$. 
For a given facility $j$, we denote with $\xi_1, \ldots, \xi_{w_j}$ its $w_j$ neighboring facilities; then, we may express $\mathcal{C}_j$ as
\begin{equation}
    \mathcal{C}_j = \{ \text{neighboring facility indices of $j$} \} 
    = \left\{
        \xi_1, \ldots, \xi_{w_j}
        \mid
        \xi \in \texttt{Neighborhood}(j)
    \right\}
    \label{eq:neighborhood_set_definition}
\end{equation}
where $\texttt{Neighborhood}(j)$ is a set of slot indices that are in the ``neighborhood'' of slot $j$. 
We conceive two kinds of ``neighborhoods'', one within the same LPO as $j$ but with shifted phasing, and another with similar phasing but on different LPOs. These are defined in further detail in the subsequent subsections. 

Using $\mathcal{C}_j$, the neighborhood swap heuristic which seeks to improve an existing solution by exploring neighboring candidate facility slots is conceived. 
The neighborhood swap involves creating a new decision vector $Y_{\mathrm{new}}$ by mutating one of the active slot in the existing best feasible $Y_{\mathrm{best}}$, denoted by $j_{\mathrm{out}}$, to one of its neighboring slot in $\mathcal{C}_{j_{\mathrm{out}}}$, and evaluating the resulting feasible allocations $X_{\mathrm{new}}$ and $\theta_{\mathrm{new}}$ using the heuristic feasible allocation Algorithm~\ref{alg:feasible_allocation}. 
This process is repeated for each candidate neighboring slots $\mathcal{C}_j$ of each active slot in $Y_{\mathrm{best}}$. 
{The neighborhood swap heuristic procedure is} summarized in Algorithm~\ref{alg:neighborhood_swap}. 

\begin{algorithm}
\caption{Neighborhood Swap}
\label{alg:neighborhood_swap}
\begin{algorithmic}[1]
    \Require $C_j$, $X_{\mathrm{best}}$, $Y_{\mathrm{best}}$, $Z_{\mathrm{best}}$, $\theta_{\mathrm{best}}$

    \State $\mathcal{J}^{Y} \gets \left\{ j \in 1, \ldots, n \mid {Y_{\mathrm{best}}}_j = 1 \right\}$
    
    \For{$j_{\mathrm{out}} \in \mathcal{J}^{Y}$}
        \For{$j_{\mathrm{in}} \in \mathcal{C}_{j_{\mathrm{out}}}$}
            \State $\mathcal{J}^{Y_{\mathrm{new}}} \gets \mathcal{J}^{Y} \cup \{j_{\mathrm{in}}\} \setminus j_{\mathrm{out}}$
                \Comment{Create new set of facilities by removing $j_{\mathrm{out}}$ and adding $j_{\mathrm{in}}$}
                
            \State $Y_{\mathrm{new}} \gets$ \verb|SetToVector|$(\mathcal{J}^{Y_{\mathrm{new}}})$
                \Comment{Convert set to vector}
                
            \State $X_{\mathrm{new}}, \theta_{\mathrm{new}} \gets$ \verb|HeuristicFeasibleAllocation|$(Y_{\mathrm{new}})$
                \Comment{Heuristically compute feasible allocation}
            
            \State $Z_{\mathrm{new}} \gets $ eqn.~\eqref{eq:TEpMP_slack_objective}
                \Comment{Evaluate {linear objective function}}
            
            \If{$Z_{\mathrm{new}} > Z_{\mathrm{best}}$}
                    \Comment{If swap results in improvement, overwrite best solution}
                \State $X_{\mathrm{best}} \gets X_{\mathrm{new}}$
                \State $Y_{\mathrm{best}} \gets Y_{\mathrm{new}}$
                \State $Z_{\mathrm{best}} \gets Z_{\mathrm{new}}$
                \State $\theta_{\mathrm{best}} \gets \theta_{\mathrm{new}}$
            \EndIf 
        \EndFor
    \EndFor
    \State \Return $X_{\mathrm{best}}$, $Y_{\mathrm{best}}$, $Z_{\mathrm{best}}$, $\theta_{\mathrm{best}}$
\end{algorithmic}
\end{algorithm}

\begin{figure}
    \centering
    \begin{subfigure}[b]{0.49\textwidth}
         \centering
         \includegraphics[width=\textwidth]{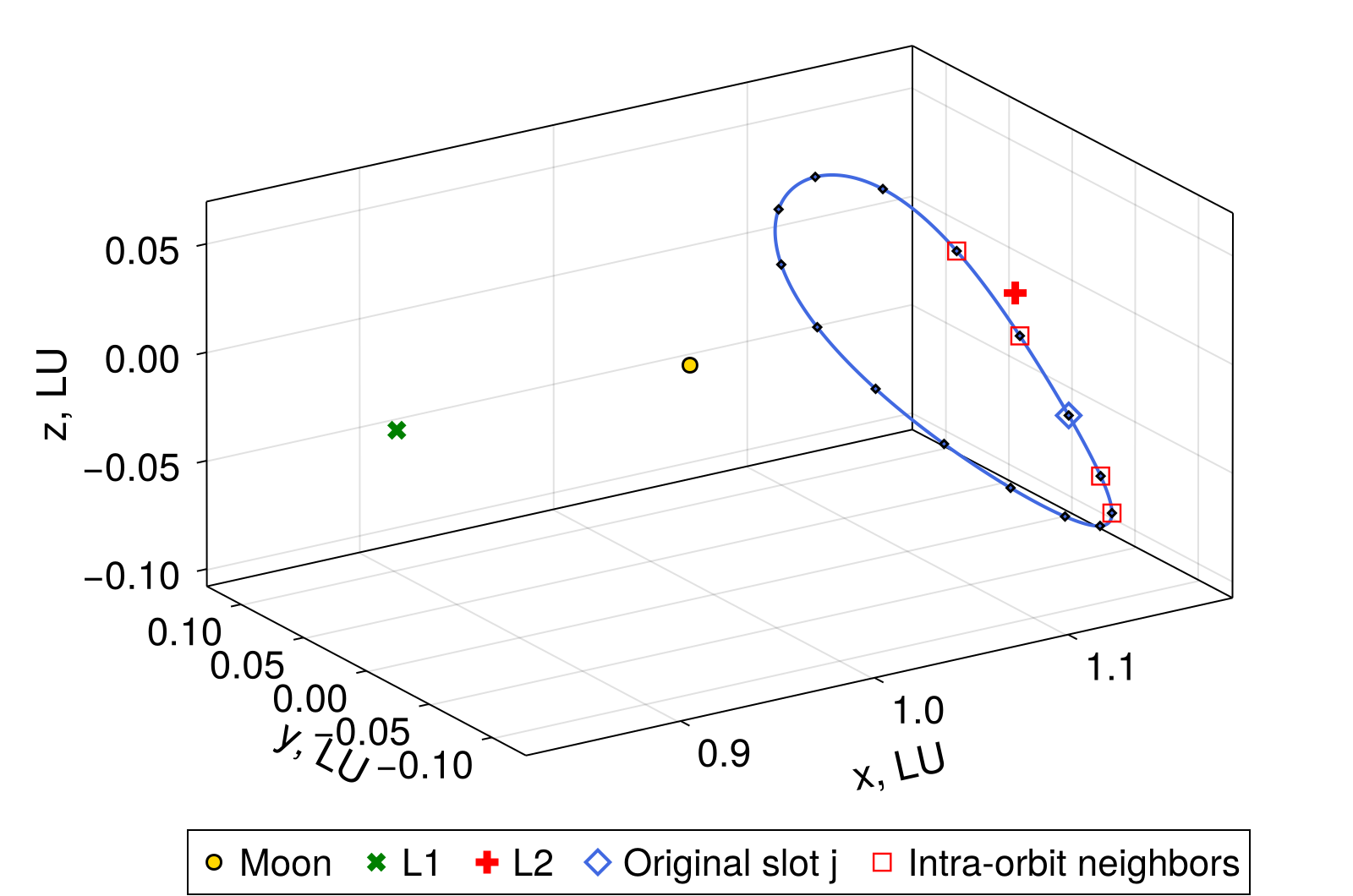}
         \caption{Intra-orbit neighborhood}
         \label{fig:neighborhood_intra}
    \end{subfigure}
    \hfill
    \begin{subfigure}[b]{0.49\textwidth}
         \centering
         \includegraphics[width=\textwidth]{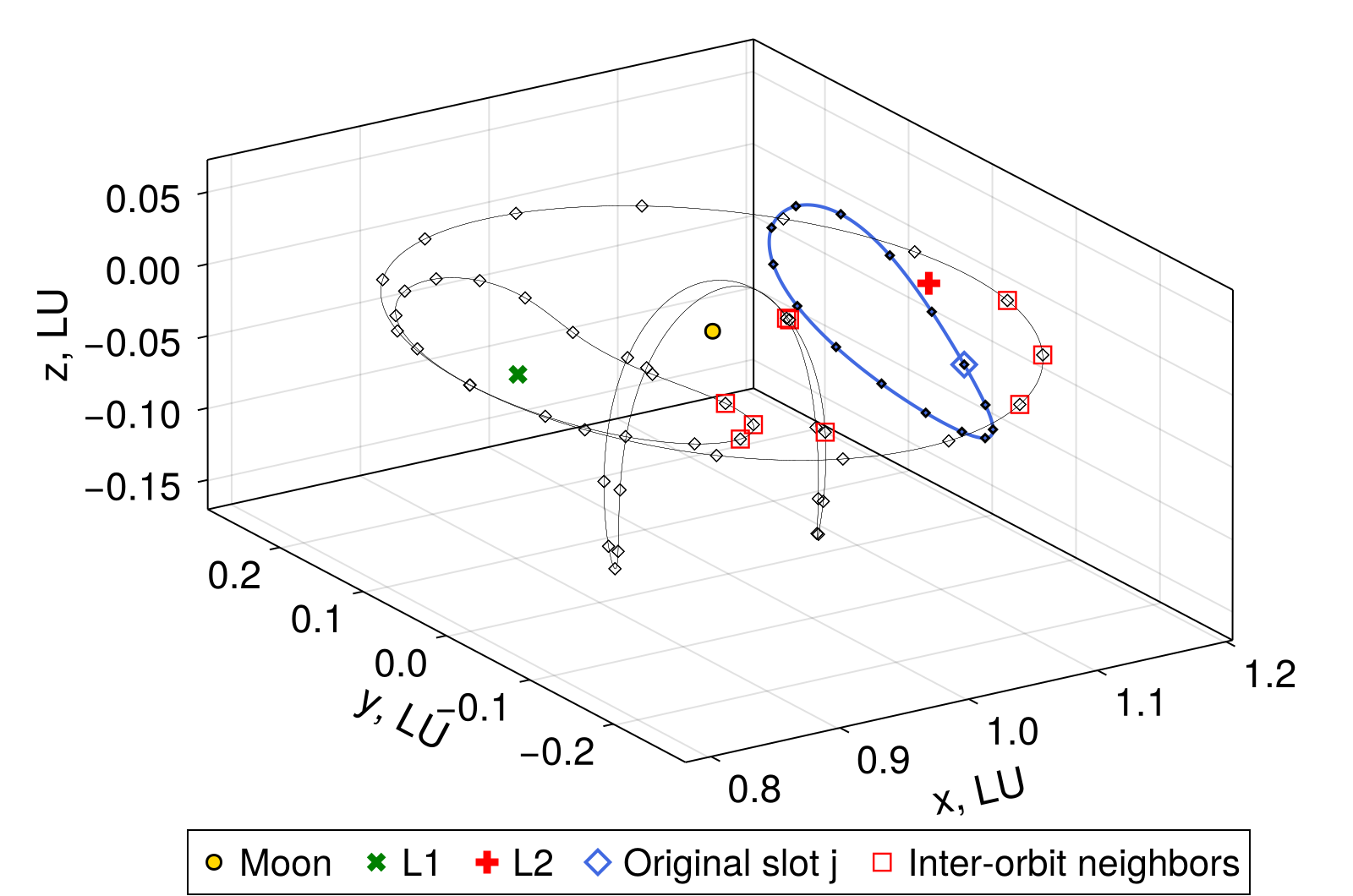}
         \caption{Inter-orbit neighborhood (subset)}
         \label{fig:neighborhood_inter}
    \end{subfigure}
    \caption{Example of intra- and inter-orbit neighborhoods of a slot on the 2:1 resonant L2 Southern Halo Orbit}
    \label{fig:slot_neighborhood}
\end{figure}

\subsubsection{Intra-Orbit Neighborhood}
The intra-orbit neighborhood is defined based on proximity in angular location along a given LPO, as shown in Figure~\ref{fig:neighborhood_intra}. 
Specifically, let $\Xi_j \subset \mathcal{J}$ denote all other facility indices that are on the same LPO as $j$; for a given $j$, its LPO is discretized into $b_j$ slots according to equation~\eqref{eq:discretize_lpo_slots}, then $\Xi_j$ has $b_j - 1$ elements, excluding $j$ itself. 
We further define $\Delta \varphi_{j \xi}$ as the angular separation between facility locations $j$ and $\xi \in \Xi_j$.  
Then, we define the permutation of $\xi$, denoted by $\pi_{\rm intra}$, for which the absolute value of $\Delta \phi_{j \xi}$ is in ascending order, 
\begin{equation}
    | \Delta \varphi_{j \pi_{\rm intra}(1)} | \leq | \Delta \varphi_{j \pi_{\rm intra}(2)} | 
    \leq \ldots 
    \leq | \Delta \varphi_{j \pi_{\rm intra}(b_j - 1)} | 
\end{equation}
We finally construct the intra-orbit neighborhood $\mathcal{C}^{\rm intra}_j$ of size $c^{\rm intra} \leq b_j - 1$ by taking the first $c^{\rm intra}$ elements in the permutation $\pi_{\rm intra}$, 
\begin{equation}
    \mathcal{C}_j^{\rm intra} = \left\{ \xi_{\kappa} \, \mid \,
    \kappa = \pi_{\rm intra}(1), \ldots, \pi_{\rm intra}(c^{\rm intra})
    \right\}
\end{equation}
Note that $c^{\rm intra}$ is a tuning parameter for the heuristics that dictates the angular separation to be considered for swapping; while a larger $c^{\rm intra}$ results in a more thorough evaluation of the neighborhood, the size of $\mathcal{C}_j^{\rm intra}$ directly impacts the number of iterations of the inner for-loop in Algorithm~\ref{alg:neighborhood_swap}, and must thus be kept reasonably low. It is also recommended to choose $c^{\rm intra}$ as a multiple of 2 to evenly explore candidate locations $\xi$ that either lead or lag in angular location relative to the originally chosen $j$.

\subsubsection{Inter-Orbit Neighborhood}
The inter-orbit neighborhood of $j$ is defined as the set of candidate facility locations that are (1) on different LPOs with the same $M$:$N$ resonance, and (2) exhibit a similar solar phase angle $\phi_S$ with respect to a reference target at a reference epoch.
The inter-orbit neighborhood for the same original slot as in Figure~\ref{fig:neighborhood_intra} is shown in Figure~\ref{fig:neighborhood_inter}.
The latter condition is required to isolate a single slot from the slots along an LPO with the same resonance as the LPO along which slot $j$ is located. 

To arrive at an expression for the inter-orbit neighborhood of $j$, we begin by defining by $\mathcal{P}_j$ the set of LPOs that have the same $M$:$N$ resonance as the LPO of $j$.
We also define a reference target position $\rbold_{k_{\rm ref}}$; in this work, we choose the mean location across all targets, $\rbold_{k_{\rm ref}} = \sum_k \rbold_k$. 
Then, for each candidate facility location $\xi$ along each LPO $\rho \in \mathcal{P}_j$, we determine the slot $\xi^*_{\rho}$ that achieves the closest solar phase angle with respect to $\rbold_{k_{\rm ref}}$ as the original slot $j$, given by
\begin{equation}
    \xi^*_{\rho} = \underset{\xi}{\operatorname{argmin}} 
    \left( \boldsymbol{l}_{j k_{\rm ref}}^T \boldsymbol{l}_{S k_{\rm ref}} 
        - \boldsymbol{l}_{\xi k_{\rm ref}}^T \boldsymbol{l}_{S k_{\rm ref}} 
    \right)
\end{equation}
The choice of $\boldsymbol{l}_{Sk}$ is a direct consequence of choosing a reference epoch. 
In this work, we choose the initial epoch with Earth-Moon-Sun alignment as the reference, and define the reference target position as the weighted sum of all considered targets across all times,
\begin{equation}
    \rbold_{k_{\rm ref}} = \dfrac{1}{\ell q} \sum_{t = 1}^{\ell} \sum_{k=1}^q \rbold_k(t)
\end{equation}
Finally, the inter-orbit neighborhood $\mathcal{C}_j^{\rm inter}$ is given by
\begin{equation}
    \mathcal{C}_j^{\rm inter}
    =
    \left\{
        \xi^*_{\rho} \, \mid \, \forall \rho \in \mathcal{P}_j
    \right\}
\end{equation}
The inter-orbit neighborhood swap is only attempted if the relaxed solution and the corresponding feasible solution from Algorithm~\ref{alg:feasible_allocation} fail to reduce the optimality gap over a pre-defined number of iterations. 
We place this restriction due to two considerations: firstly, we avoid prematurely swapping to different LPO slots before the nearby slots in phase along a given LPO has been sufficiently explored. 
Secondly, the inter-orbit neighborhood Algorithm~\ref{alg:neighborhood_swap} requires running the intra-orbit neighborhood Algorithm~\ref{alg:feasible_allocation} for each candidate solution, thus making it more computationally intensive.

\subsection{Lagrange Multipliers Update with Subgradient Method}
\label{sec:subgradient_optimization}
The Lagrange multipliers are updated according to the subgradient method {according to~\eqref{eq:LM_subgradient_multiplier_update}; for the TE-$p$-MP, the multiplier update formulae are}
\begin{align}
    \lambda_{jt}^{(h+1)} &= \lambda_{jt}^{(h)} + s^{(h)} \left(
        \sum_{i=1}^m M_{ijtk} \bar{X}_{ijt} - 1
    \right)
    \quad \forall j, t
    \\
    \eta_{tk}^{(h+1)} &= \lambda_{jt}^{(h)} + s^{(h)} \left(
        \bar{\theta}_{tk} - \sum_{i=1}^m \sum_{j=1}^n M_{ijtk} \bar{X}_{ijt}
    \right)
    \quad \forall t,k
\end{align}
where the expression for the step-size~\eqref{eq:step_definition} becomes
\begin{equation}
    s^{(h)} = \dfrac{\mu^{(h)} (Z_{\mathrm{relax}} - Z_{\mathrm{feas}})}{
    \left[ \sum_{j=1}^n \sum_{t=1}^{\ell}
        \max \left(0,\sum_{i=1}^m M_{ijtk} \bar{X}_{ijt} - 1
    \right) \right]^2
    \left[ \sum_{t=1}^{\ell} \sum_{k=1}^{q}
        \max \left(0, 
            \bar{\theta}_{tk} - \sum_{i=1}^m \sum_{j=1}^n M_{ijtk} \bar{X}_{ijt}
        \right)
    \right]^2
    }
    \label{eq:step_definition_TEpMP}
\end{equation}

\subsection{Considerations for Implementing the Lagrangian Method}
\label{sec:conding_considerations}
The LM is implemented using the Julia language. We make note of three particular design choices that are made to improve the performance of the algorithm, namely on the initialization of the Lagrange multipliers, parallelism and function memorization. 

\subsubsection{Initialization of Lagrange Multipliers}
The choice of the initial values of the Lagrange multipliers requires some care, as the LM does not have a guarantee of convergence. One approach commonly adopted is to initialize them using the dual of the corresponding constraints for the linear program (LP)-relaxed problem, but our initial investigation has shown that this does not yield a particular advantage, despite the additional cost of having to solve the relaxed problem.

\subsubsection{Parallelism}
For a given choice of facilities $Y$, the feasible allocation construction Algorithm~\ref{alg:feasible_allocation} makes allocation decisions $X_{ijt}$ for each time step, each of which is independent of allocation decisions at other time steps. 
This makes the algorithm well-suited for parallelism; we opt for the multiple-threading paradigm through Julia's \verb|Threads.@threads| decorator to the for-loop in line 3 through 20, as the operation within each loop is not too expensive, and we gain more benefit from a shared memory. 
Note that both the intra- and inter-orbit neighborhood swap schemes require recomputing the best heuristic allocations $X_{ijt}$ for each new hypothesized $Y$, thus a speed-up of Algorithm~\ref{alg:feasible_allocation} is particularly beneficial to speed up the overall Lagrangian method. 

\subsubsection{Memoization}
During the heuristic process of the Lagrangian method, it is possible that the Lagrangian relaxed solution $\bar{Y}$ or the result of either an intra- and inter-orbit neighborhood swap result in a choice of facilities that has already been considered by the algorithm. 
We make use of function memoization on the feasible allocation construction Algorithm~\ref{alg:feasible_allocation} by storing the outputs $\underline{X}$ and $\underline{\theta}$ that result from a given set of active facility indices $\mathcal{J}^{\bar{Y}}$.

\section{Numerical Results}
\label{sec:numerical_results}
The TE-$p$-MP is solved for various choices of observation demand as well as observer parameters. 
{
Each instance is solved using the LM and B\&B via Gurobi.
Gurobi is a state-of-the-art optimization suite with a particular focus on B\&B mixed-integer linear programs and serves as a point of comparison to study the efficacy of the proposed LM.}
This Section initially focuses on comparing the performance of the two methods, emphasizing on the quality of solutions they yield. 
Then, closer attention is drawn to the obtained architectures for a select number of instances, in order to provide insights on the CSSA problem. 

The two sets of demand presented in Section~\ref{sec:demand_definition} are considered.
The corresponding problem dimensions are summarized in Table~\ref{tab:problem_dimensions}. 
For each set of demand introduced in Section~\ref{sec:demand_definition}, we consider LPOs given in Appendix~\ref{sec:appendix_lpo_conditions}, with an approximate time-spacing parameter $\Delta t_b = 12$ hours. Using equation~\eqref{eq:discretize_lpo_slots}, this results in $n = 1212$ candidate facility slots in total, also presented in Appendix~\ref{sec:appendix_lpo_conditions}. 
{Since the set of candidate observer LPOs is assumed to be resonant, any architecture involving any combination of LPOs will return to its original configuration after $N_{\rm syn}$ synodic period, as introduced in Section~\ref{sec:LPO}. 
Leveraging this feature, we consider a TE-$p$-MP with a time horizon of $N_{\rm syn} = 4$ synodic months. 
Each synodic month is discretized with a time-step of 30 steps per synodic month, corresponding to roughly 24 hours, thus resulting in $\ell = 120$ time-steps.  
Note that the $N_{\rm syn}$-synodic period-repeating feature of the candidate LPOs is not necessitated by the TE-$p$-MP formulation itself; one may choose to formulate the problem with non-repeating candidate observer locations, and set a longer time-horizon as well.
}

All problem instances are solved on a desktop PC running on 20 i7-12700 CPU cores at 2.10GHz base speed. 
Table~\ref{tab:options_lagrangean_method} gives the hyperparameters for the LM.
For B\&B, we employ Gurobi version 11.0.1, using the barrier method to solve the initial root relaxation as it was found through preliminary experiments to improve the algorithm's performance to TE-$p$-MP instances. 
The optimality gap tolerance is set to $1\%$. 
All other hyperparameters are left unchanged from their default values.

\begin{table}[]
\centering
\caption{Problem dimensions for each demand and observer parameters}
\begin{tabular}{@{}lllll@{}}
\toprule
Demand                &           &            & Cone of Shame & LET transit \\ \midrule
Number of pointing directions $m$ & &           & 14  & 14  \\
Number of facility locations $n$  & &           & 1212 & 1212 \\
Number of time steps $\ell$ &        &          & 120 & 120 \\
Number of targets $q$ &           &             & 304           & 675         \\ \midrule
\multirow{6}{*}{Fraction of entries in $M^b$}
&                     & $\bar{m}_{\mathrm{crit}} = 15$  & 0.0047 & 0.0007\\
& FOV $=60^{\circ}$    & $\bar{m}_{\mathrm{crit}} = 18$  & 0.0460 & 0.0043 \\
&                     & $\bar{m}_{\mathrm{crit}} = 20$  & 0.0581 & 0.0049 \\
\cmidrule{2-5}%
&                     & $\bar{m}_{\mathrm{crit}} = 15$  & 0.0180 & 0.0026 \\
& FOV $=120^{\circ}$ & $\bar{m}_{\mathrm{crit}} = 18$  & 0.1711 & 0.0176 \\
&                     & $\bar{m}_{\mathrm{crit}} = 20$  & 0.2177 & 0.0204 \\
\bottomrule
\end{tabular}
\label{tab:problem_dimensions}
\end{table}

\subsection{Comparison of Lagrangian Method against Branch-and-Bound}
We first compare the proposed Lagrangian method against Gurobi's B\&B, for TE-$p$-MP instances with static and dynamic demands. 
For the sake of a fair comparison, we set a time limit of 500 seconds on both approaches to assess their performance. 
In addition, the solution obtained by Gurobi with a time limit of 3600 seconds is reported. Note that while extending the time limit may allow B\&B to find better solutions, this is not necessarily desirable, especially when various problem parameters, such as the number of satellites $p$, cut of magnitude, or FOV, must be explored and optimized. 
Nevertheless, the solution from Gurobi with a 3600-second time limit serves as a proxy of a near-optimal solution obtained with a \textit{reasonably} long time limit in most instances.

Tables~\ref{tab:SSA_pmedian_gurboi_vs_LR_cone1_fov60}, and~\ref{tab:SSA_pmedian_gurboi_vs_LR_let1_fov60} report the sum of targets observed across all time steps, $\Theta$, and the actual CPU solve time for each of the two observation demands, respectively. 
In both cases, the solutions obtained with Gurobi using time limits of 500 and 3600 seconds and with the LM using a time limit of 500 seconds are reported. 
In addition, the percentage difference of the Lagrangian method solution with respect to each B\&B solution is reported. 
Note that the actual CPU solve time may exceed 500 seconds for both Gurobi and the Lagrangian method; in Gurobi, this is due to the computation of attributes associated with the terminated solution, while with the Lagrangian method, the solve time is checked between the generation of feasible lower bound solution, which may start before the solve time limit, but exceed it before another solve time check happens. 

Across both demands, we first note that LM finds a superior solution to B\&B with the 500 seconds time limit (i.e., under a fair condition) in all cases with $\bar{m}_{\mathrm{crit}} = 20$, and in all except two cases with $\bar{m}_{\mathrm{crit}} = 18$, the latter of which corresponding to $p=2$ with the Cone of Shame and LET demands. For cases with $\bar{m}_{\mathrm{crit}} = 15$, B\&B finds competitive solutions for the Cone of Shame demand, while it is found to struggle with the LET demand for $p\geq 3$. 
Letting Gurobi run for up to 3600 seconds improves the solution in some instances, but the behavior is hard to predict. 
For example, in the Cone of Shame demand, cases with $\bar{m}_{\mathrm{crit}} = 18$ and $p=3$ has a solution that improves from $\Theta_{500} = 0$ to $\Theta_{3600} = 0.6355$, while with $\bar{m}_{\mathrm{crit}} = 18$ and $p=4$, the solution improves only from $\Theta_{500} = 0.6556$ to $\Theta_{3600} = 0.6975$.

The performance of B\&B, both after 500 and 3600 seconds, may be seen as a measure of the difficulty of the TE-$p$-MP instance. 
The problem appears to be easy when the overall visibility of targets is low, corresponding to instances with fewer (low $p$) observers and/or are equipped with worse sensors (lower $\bar{m}_{\mathrm{crit}}$), or when the observation is easy, with many observers (high $p$) and/or high-performance sensors (higher $\bar{m}_{\mathrm{crit}}$), resulting in an instance where the observation task is ``easy''. 
Examples of the former include the cases with $\bar{m}_{\rm crit} = 15,18$ and $p=2,3$ for the LET demand, while examples of the latter include the cases with $\bar{m}_{\rm crit} = 15$ and $p=4,5$, also for the LET demand. 
Instances of intermediate observation difficulty, such as the cases with $\bar{m}_{\rm crit} = 18,20$ and $p=2,3,4,5$ for the Cone of Shame demand, are found to result in the poorest performances from B\&B. 

The performance of the LM is consistent in terms of observation condition: increasing $p$ always increases $\Theta_{\rm LH}$, while increasing $\bar{m}_{\mathrm{crit}}$ also increases $\Theta_{\rm LH}$ except for the two cases with the Cone of Shame demand with $p=2,3$ going from $\bar{m}_{\mathrm{crit}} = 18$ to $20$, albeit with a small reduction, from $0.7334$ to $0.7297$ with $p=2$, and from $0.8192$ to $0.8080$ with $p=3$. 
Note that these two instances correspond to some of the cases where B\&B struggles the most. 
Among the 24 instances, 6 instances resulted in an LM solution that was inferior to B\&B after the 500-second time limit; within them, only 4 cases resulted in a solution that was worse by over 10\%. 
Overall, the consistency of the LM offers an attractive reliable approach for assessing the achievable observation from a given TE-$p$-MP instance. 
Alternatively, one may consider using Gurobi with a relatively short solve time limit such as 500 seconds, and use the LM as a consistent fallback solver to recover decent solutions in case B\&B under-performs.

\begin{table}[]
\centering
\caption{Lagrangian method hyperparameters}
\begin{tabular}{@{}ll@{}}
\toprule
Option & Value 
\\
\midrule
Max number of iterations & 30 \\
Gap tolerance & 0.01 \\
Max number of iterations without improvements & 10 \\
Number of iterations with no improvement to reduce step-size parameter $\mu^{(h)}$ & 5 \\
Number of slots considered for intra-orbit neighborhood swap & 4 \\
Number of iterations with no improvement to try inter-orbit neighborhood swap & 4 \\
Relocation strategy & \verb|FullFactorialGreedy| \\
\bottomrule
\end{tabular}
\label{tab:options_lagrangean_method}
\end{table}

\begin{table}[h]
    \centering
    \small
    \caption{TE-$p$-MP instance with Cone of Shame demand and FOV = $60^{\circ}$ solved using Branch-and-Bound versus Lagrangian method. All reported times are in seconds. }
    \begin{subtable}[h]{\textwidth}
        \centering
        \caption{$\bar{m}_{\mathrm{crit}} = 15$}
        \begin{tabular}{@{}lllllllllllll@{}}
        \toprule
        $p$ & \multicolumn{3}{c}{B\&B (\texttt{TimeLimit} $= 500$)} & & \multicolumn{3}{c}{B\&B (\texttt{TimeLimit} $= 3600$)} & &\multicolumn{4}{c}{Lagrangian Method} \\
        \cline{2-4} \cline{6-8} \cline{10-13} 
          & $\Theta_{\mathrm{500}}$ & Time & Status &
          & $\Theta_{\mathrm{3600}}$ & Time & Status &
          & $\Theta_{\mathrm{LH}}$ & Time 
          & \% diff. to $\Theta_{\mathrm{500}}$
          & \% diff. to $\Theta_{\mathrm{3600}}$
          \\
        \midrule
2 & 0.0419 & 101 & \verb|OPTIMAL| & & - & - & - & & 0.0402 & 620 & $-3.99$ & $-3.99$ \\
3 & 0.0608 & 501 & \verb|TIME_LIMIT| & & 0.0608 & 681 & \verb|OPTIMAL| & & 0.0533 & 575 & $-12.22$ & $-12.22$ \\
4 & 0.0778 & 500 & \verb|TIME_LIMIT| & & 0.0783 & 3601 & \verb|TIME_LIMIT| & & 0.0717 & 614 & $-7.89$ & $-8.54$ \\
5 & 0.0939 & 500 & \verb|TIME_LIMIT| & & 0.0939 & 3601 & \verb|TIME_LIMIT| & & 0.0824 & 603 & $-12.29$ & $-12.29$ \\
        \bottomrule
        \end{tabular}
        \label{tab:SSA_pmedian_gurboi_vs_LR_cone1_m15_fov60}
    \end{subtable}
    \\ \par\bigskip
    \begin{subtable}[h]{\textwidth}
        \centering
        \caption{$\bar{m}_{\mathrm{crit}} = 18$}
        \begin{tabular}{@{}lllllllllllll@{}}
        \toprule
        $p$ & \multicolumn{3}{c}{B\&B (\texttt{TimeLimit} $= 500$)} & & \multicolumn{3}{c}{B\&B (\texttt{TimeLimit} $= 3600$)} & &\multicolumn{4}{c}{Lagrangian Method} \\
        \cline{2-4} \cline{6-8} \cline{10-13} 
          & $\Theta_{\mathrm{500}}$ & Time & Status &
          & $\Theta_{\mathrm{3600}}$ & Time & Status &
          & $\Theta_{\mathrm{LH}}$ & Time 
          & \% diff. to $\Theta_{\mathrm{500}}$
          & \% diff. to $\Theta_{\mathrm{3600}}$
          \\
        \midrule
2 & 0.6058 & 501 & \verb|TIME_LIMIT| & & 0.7259 & 3601 & \verb|TIME_LIMIT| & & 0.6341 & 542 & $+4.67$ & $-12.65$ \\
3 & 0.0000 & 501 & \verb|TIME_LIMIT| & & 0.6355 & 3601 & \verb|TIME_LIMIT| & & 0.7786 & 584 & $+Inf$ & $+22.52$ \\
4 & 0.6756 & 501 & \verb|TIME_LIMIT| & & 0.6975 & 3601 & \verb|TIME_LIMIT| & & 0.8670 & 634 & $+28.32$ & $+24.31$ \\
5 & 0.0000 & 500 & \verb|TIME_LIMIT| & & 0.7576 & 3602 & \verb|TIME_LIMIT| & & 0.9142 & 703 & $+Inf$ & $+20.68$ \\
        \bottomrule
        \end{tabular}
        \label{tab:SSA_pmedian_gurboi_vs_LR_cone1_m18_fov60}
    \end{subtable}
    \\ \par\bigskip
    \begin{subtable}[h]{\textwidth}
        \centering
        \caption{$\bar{m}_{\mathrm{crit}} = 20$}
        \begin{tabular}{@{}lllllllllllll@{}}
        \toprule
        $p$ & \multicolumn{3}{c}{B\&B (\texttt{TimeLimit} $= 500$)} & & \multicolumn{3}{c}{B\&B (\texttt{TimeLimit} $= 3600$)} & &\multicolumn{4}{c}{Lagrangian Method} \\
        \cline{2-4} \cline{6-8} \cline{10-13} 
          & $\Theta_{\mathrm{500}}$ & Time & Status &
          & $\Theta_{\mathrm{3600}}$ & Time & Status &
          & $\Theta_{\mathrm{LH}}$ & Time 
          & \% diff. to $\Theta_{\mathrm{500}}$
          & \% diff. to $\Theta_{\mathrm{3600}}$
          \\
        \midrule
2 & 0.0000 & 501 & \verb|TIME_LIMIT| & & 0.6180 & 3602 & \verb|TIME_LIMIT| & & 0.6328 & 536 & $+Inf$ & $+2.39$ \\
3 & 0.0000 & 501 & \verb|TIME_LIMIT| & & 0.8577 & 3601 & \verb|TIME_LIMIT| & & 0.7223 & 593 & $+Inf$ & $-15.79$ \\
4 & 0.0000 & 501 & \verb|TIME_LIMIT| & & 0.7854 & 3602 & \verb|TIME_LIMIT| & & 0.9271 & 597 & $+Inf$ & $+18.04$ \\
5 & 0.0000 & 501 & \verb|TIME_LIMIT| & & 0.9558 & 3603 & \verb|TIME_LIMIT| & & 0.9599 & 584 & $+Inf$ & $+0.43$ \\
        \bottomrule
        \end{tabular}
        \label{tab:SSA_pmedian_gurboi_vs_LR_cone1_m20_fov60}
    \end{subtable}
     \label{tab:SSA_pmedian_gurboi_vs_LR_cone1_fov60}
\end{table}

\begin{table}[h]
    \centering
    \small
    \caption{TE-$p$-MP instance with LET demand and FOV = $60^{\circ}$ solved using Branch-and-Bound versus Lagrangian method. All reported times are in seconds. }
    \begin{subtable}[h]{\textwidth}
        \centering
        \caption{$\bar{m}_{\mathrm{crit}} = 15$}
        \begin{tabular}{@{}lllllllllllll@{}}
        \toprule
        $p$ & \multicolumn{3}{c}{B\&B (\texttt{TimeLimit} $= 500$)} & & \multicolumn{3}{c}{B\&B (\texttt{TimeLimit} $= 3600$)} & &\multicolumn{4}{c}{Lagrangian Method} \\
        \cline{2-4} \cline{6-8} \cline{10-13} 
          & $\Theta_{\mathrm{500}}$ & Time & Status &
          & $\Theta_{\mathrm{3600}}$ & Time & Status &
          & $\Theta_{\mathrm{LH}}$ & Time 
          & \% diff. to $\Theta_{\mathrm{500}}$
          & \% diff. to $\Theta_{\mathrm{3600}}$
          \\
        \midrule
2 & 0.3737 & 66 & \verb|OPTIMAL| & & - & - & - & & 0.3144 & 523 & $-15.86$ & $-15.86$ \\
3 & 0.3246 & 501 & \verb|TIME_LIMIT| & & 0.5078 & 1426 & \verb|OPTIMAL| & & 0.3967 & 644 & $+22.20$ & $-21.88$ \\
4 & 0.1703 & 501 & \verb|TIME_LIMIT| & & 0.5973 & 1846 & \verb|OPTIMAL| & & 0.4681 & 519 & $+174.92$ & $-21.63$ \\
5 & 0.1845 & 500 & \verb|TIME_LIMIT| & & 0.6670 & 3601 & \verb|TIME_LIMIT| & & 0.5471 & 537 & $+196.50$ & $-17.97$ \\
        \bottomrule
        \end{tabular}
        \label{tab:SSA_pmedian_gurboi_vs_LR_let1_m15_fov60}
    \end{subtable}
    \\ \par\bigskip
    \begin{subtable}[h]{\textwidth}
        \centering
        \caption{$\bar{m}_{\mathrm{crit}} = 18$}
        \begin{tabular}{@{}lllllllllllll@{}}
        \toprule
        $p$ & \multicolumn{3}{c}{B\&B (\texttt{TimeLimit} $= 500$)} & & \multicolumn{3}{c}{B\&B (\texttt{TimeLimit} $= 3600$)} & &\multicolumn{4}{c}{Lagrangian Method} \\
        \cline{2-4} \cline{6-8} \cline{10-13} 
          & $\Theta_{\mathrm{500}}$ & Time & Status &
          & $\Theta_{\mathrm{3600}}$ & Time & Status &
          & $\Theta_{\mathrm{LH}}$ & Time 
          & \% diff. to $\Theta_{\mathrm{500}}$
          & \% diff. to $\Theta_{\mathrm{3600}}$
          \\
        \midrule
2 & 0.9309 & 500 & \verb|TIME_LIMIT| & & 0.9398 & 2263 & \verb|OPTIMAL| & & 0.8339 & 638 & $-10.42$ & $-11.27$ \\
3 & 0.0000 & 501 & \verb|TIME_LIMIT| & & 0.0000 & 3600 & \verb|TIME_LIMIT| & & 0.9056 & 634 & $+Inf$ & $+Inf$ \\
4 & 0.0000 & 500 & \verb|TIME_LIMIT| & & 0.0000 & 3600 & \verb|TIME_LIMIT| & & 0.9355 & 524 & $+Inf$ & $+Inf$ \\
5 & 0.0000 & 500 & \verb|TIME_LIMIT| & & 0.0000 & 3600 & \verb|TIME_LIMIT| & & 0.9461 & 569 & $+Inf$ & $+Inf$ \\
        \bottomrule
        \end{tabular}
        \label{tab:SSA_pmedian_gurboi_vs_LR_let1_m18_fov60}
    \end{subtable}
    \\ \par\bigskip
    \begin{subtable}[h]{\textwidth}
        \centering
        \caption{$\bar{m}_{\mathrm{crit}} = 20$}
        \begin{tabular}{@{}lllllllllllll@{}}
        \toprule
        $p$ & \multicolumn{3}{c}{B\&B (\texttt{TimeLimit} $= 500$)} & & \multicolumn{3}{c}{B\&B (\texttt{TimeLimit} $= 3600$)} & &\multicolumn{4}{c}{Lagrangian Method} \\
        \cline{2-4} \cline{6-8} \cline{10-13} 
          & $\Theta_{\mathrm{500}}$ & Time & Status &
          & $\Theta_{\mathrm{3600}}$ & Time & Status &
          & $\Theta_{\mathrm{LH}}$ & Time 
          & \% diff. to $\Theta_{\mathrm{500}}$
          & \% diff. to $\Theta_{\mathrm{3600}}$
          \\
        \midrule
2 & 0.4778 & 500 & \verb|TIME_LIMIT| & & 0.9701 & 3601 & \verb|TIME_LIMIT| & & 0.8681 & 526 & $+81.70$ & $-10.51$ \\
3 & 0.4780 & 501 & \verb|TIME_LIMIT| & & 0.4780 & 3600 & \verb|TIME_LIMIT| & & 0.9461 & 537 & $+97.91$ & $+97.91$ \\
4 & 0.4794 & 500 & \verb|TIME_LIMIT| & & 0.9957 & 819 & \verb|OPTIMAL| & & 0.9868 & 547 & $+105.84$ & $-0.89$ \\
5 & 0.4794 & 500 & \verb|TIME_LIMIT| & & 0.9993 & 551 & \verb|OPTIMAL| & & 0.9774 & 560 & $+103.88$ & $-2.20$ \\
        \bottomrule
        \end{tabular}
        \label{tab:SSA_pmedian_gurboi_vs_LR_let1_m20_fov60}
    \end{subtable}
     \label{tab:SSA_pmedian_gurboi_vs_LR_let1_fov60}
\end{table}

\subsection{Analysis of CSSA Coverage}
We now look in detail at the trade space for designing CSSA constellations for the cone of shame, and LET transit window demands. 
For consistency across all cases, all results are based on solutions obtained using the LM, with a time limit of 1000 seconds. 
We provide insight on a select number of architectures across a few combinations of parameters. The full list of LPOs used by each problem from Gurobi and LM are provided in Tables in Appendix~\ref{appendix:appendix_lpos_used_by_solutions}. 

\subsubsection{Cone of Shame}
\begin{figure}
    \centering
    \begin{subfigure}[b]{0.45\textwidth}
        \centering
        \includegraphics[width=\textwidth]{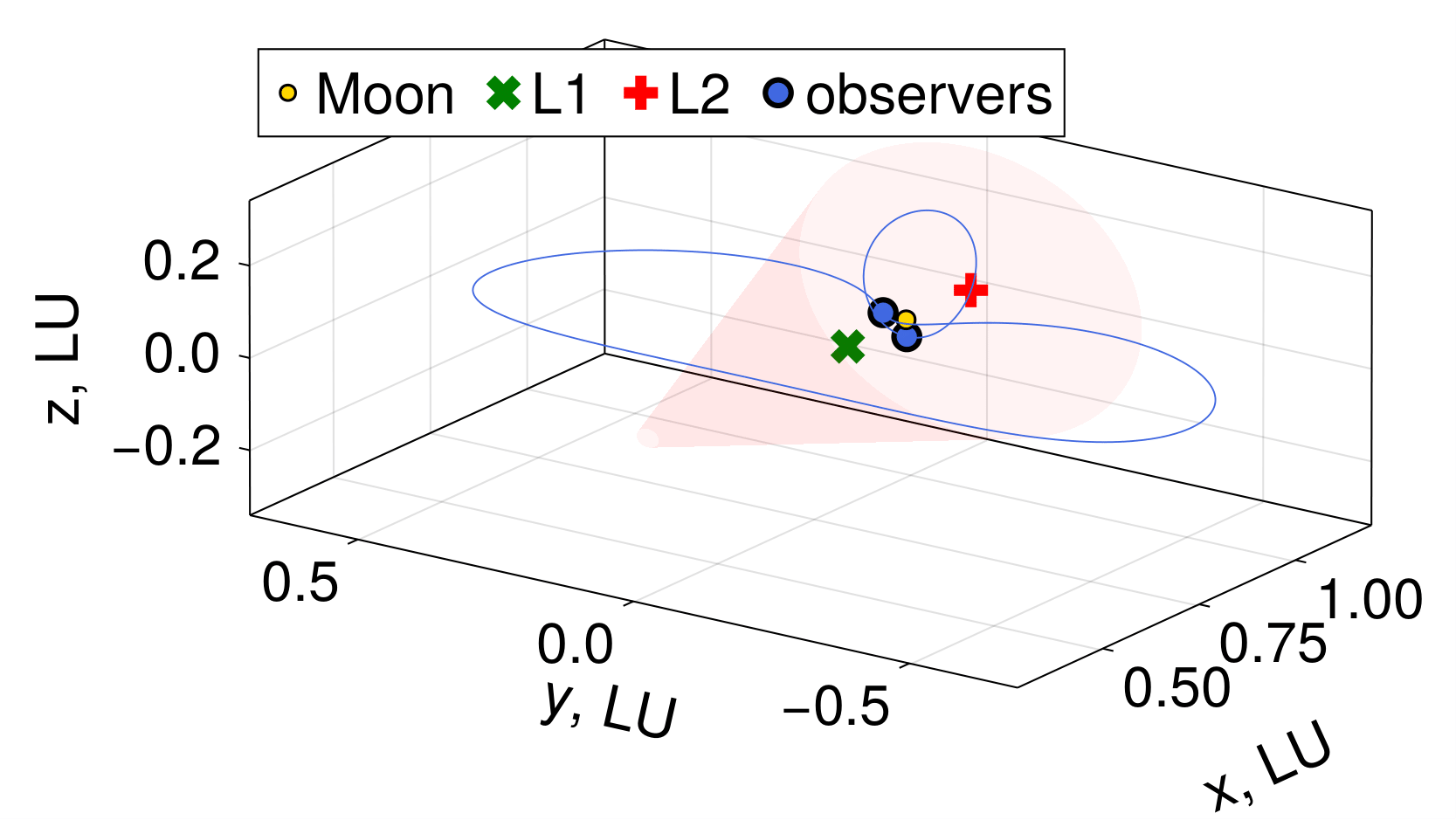}
        \caption{Lagrangian method solution with $\bar{m}_{\mathrm{crit}}=20$, $p=2$}
        \label{fig:arch_traj_cone1_Mcut20_p2_lm}
    \end{subfigure}
    \hfill
    \begin{subfigure}[b]{0.45\textwidth}
        \centering
        \includegraphics[width=\textwidth]{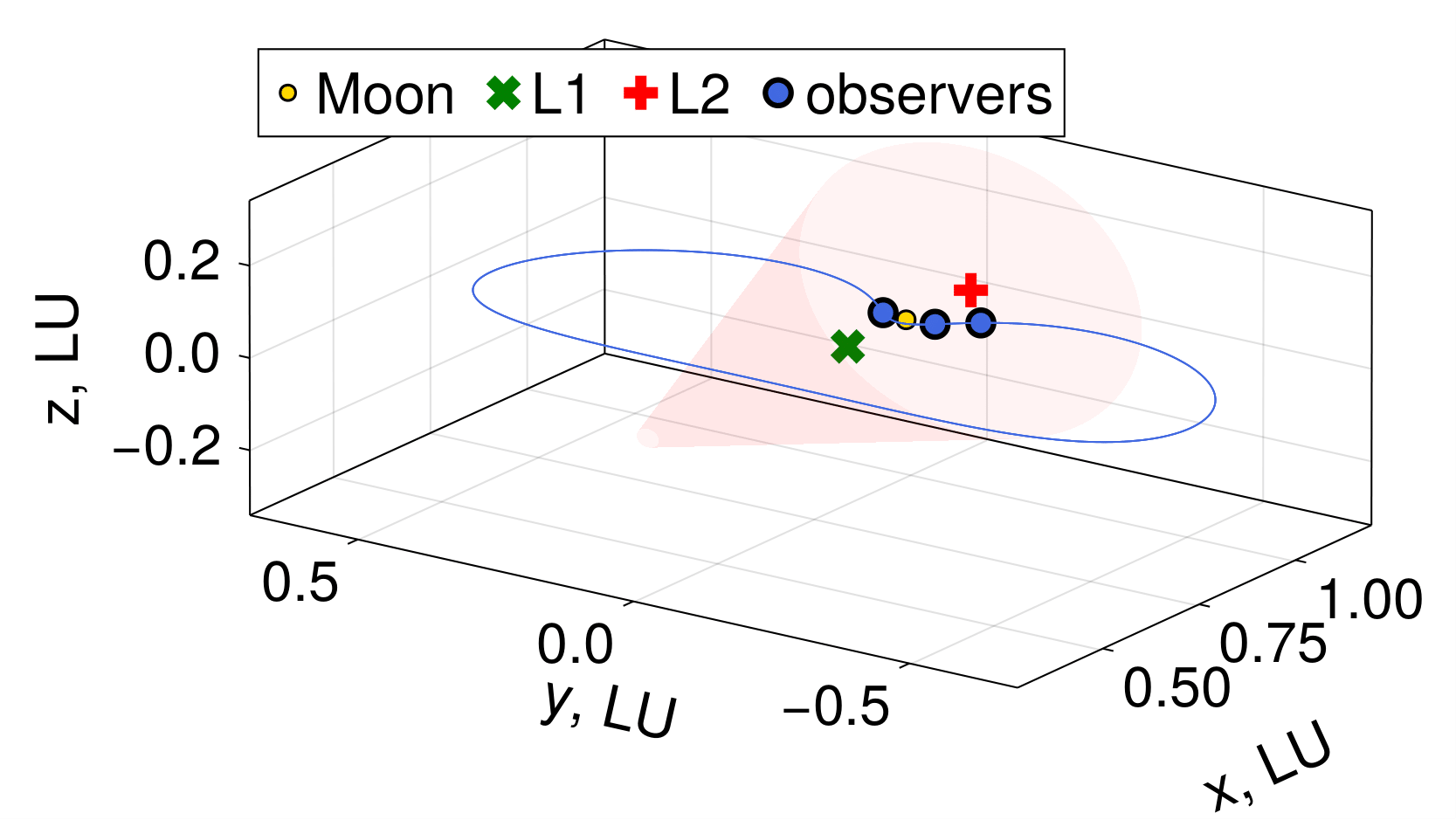}
        \caption{Lagrangian method solution with $\bar{m}_{\mathrm{crit}}=20$, $p=3$}
        \label{fig:arch_traj_cone1_Mcut20_p3_lm}
    \end{subfigure}
    \\
    \begin{subfigure}[b]{0.45\textwidth}
        \centering
        \includegraphics[width=\textwidth]{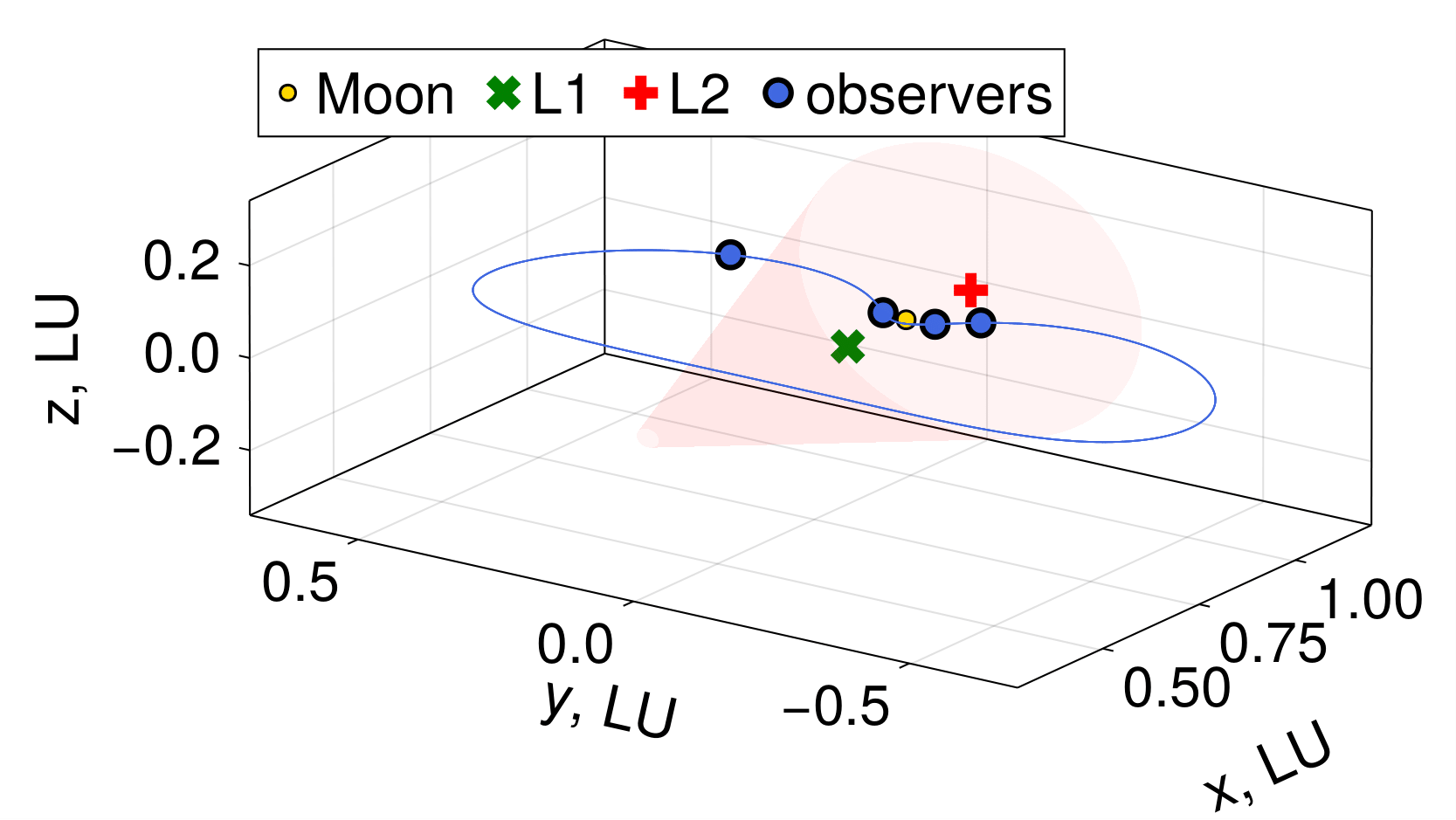}
        \caption{Lagrangian method solution with $\bar{m}_{\mathrm{crit}}=20$, $p=4$}
        \label{fig:arch_traj_cone1_Mcut20_p4_lm}
    \end{subfigure}
    \hfill
    \begin{subfigure}[b]{0.45\textwidth}
        \centering
        \includegraphics[width=\textwidth]{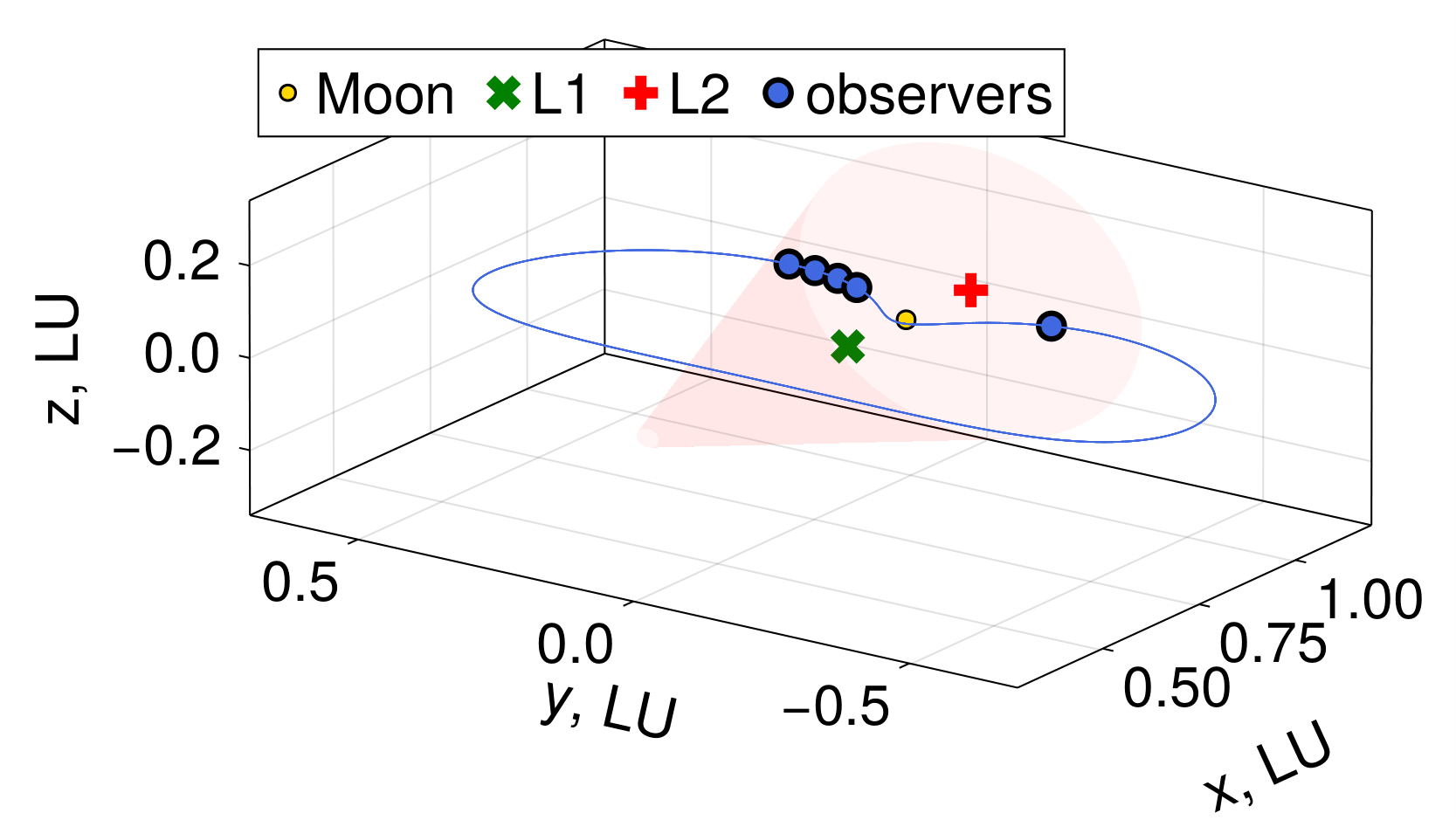}
        \caption{Lagrangian method solution with $\bar{m}_{\mathrm{crit}}=20$, $p=5$}
        \label{fig:arch_traj_cone1_Mcut20_p5_lm}
    \end{subfigure}
    \caption{Optimized observer locations for Cone of Shame demand at time-step $t=1$ (at a new Moon)}
    \label{fig:arch_traj_cone1}
\end{figure}

Figure~\ref{fig:arch_traj_cone1} shows the LM-based solutions for $p=2,3,4,5$ and $\bar{m}_{\mathrm{crit}}=20$. 
In the Cone of Shame demand, the 1:1 resonant L1 Lyapunov orbit is found to be used in many competitive instances. 
As an exception, when only two observers are placed, an alternative constellation configuration involving also a 3:1 resonant halo orbit is found. 
In these three solutions, as the number of observers $p$ is increased, an additional observer is added to the same orbit, increasing the coverage fraction as well. 

Clock-wise rotating LPOs with 1:1 resonance are particularly well suited for optical observations since the solar phase angle $\phi_S$ in the direction of a static target remains nearly constant~\cite{Vendl2021,Visonneau2023}.
Even though the 1:1 L1 Lyapunov has a varying distance from the set of targets, choosing the appropriate phase along the LPO provides consistently low $\bar{m}_{\mathrm{target}}$, due to the strong dependence of $\bar{m}_{\mathrm{target}}$ on $\phi_S$ over the distance, as shown in Figure~\ref{fig:visibility_magnitude_contour}. 
To benefit from favorable constellations, we find solutions across values of $p$ and $\bar{m}_{\mathrm{crit}}$ where observers are thus placed on the 1:1 resonant Lyapunov with a leading/trailing configuration rather than being placed at equal anomaly separation. 
The persistence of this train-like configuration across multiple problem instances highlights the efficacy of this constellation design strategy, particularly for monitoring targets in the Cone of Shame. 
We note that these constellations are unlike typical Earth-bound constellations that are equally spaced along a given orbit; rather than being equally distributed along track, they form a leading/trailing formation. 

\begin{figure}
    \centering
    \begin{subfigure}[b]{0.99\textwidth}
        \centering
        \includegraphics[width=\textwidth]{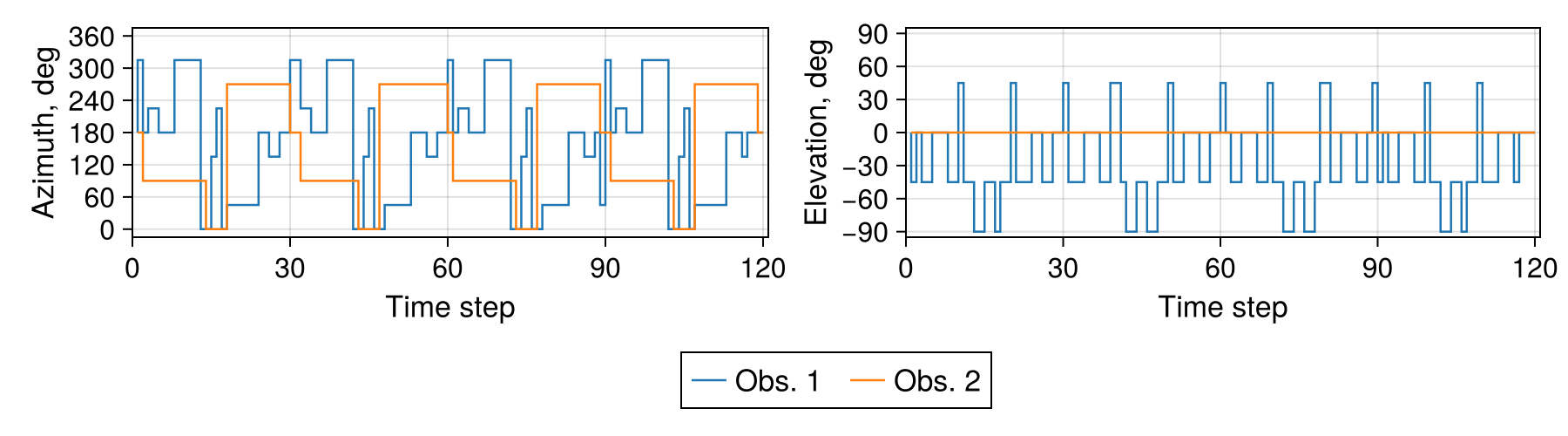}
        \caption{Lagrangian method solution with $\bar{m}_{\mathrm{crit}}=20$, $p=2$}
        \label{fig:lagrangean_targetcone1_p2_Mcut20_fov60}
    \end{subfigure}
    \\
    \begin{subfigure}[b]{0.99\textwidth}
        \centering
        \includegraphics[width=\textwidth]{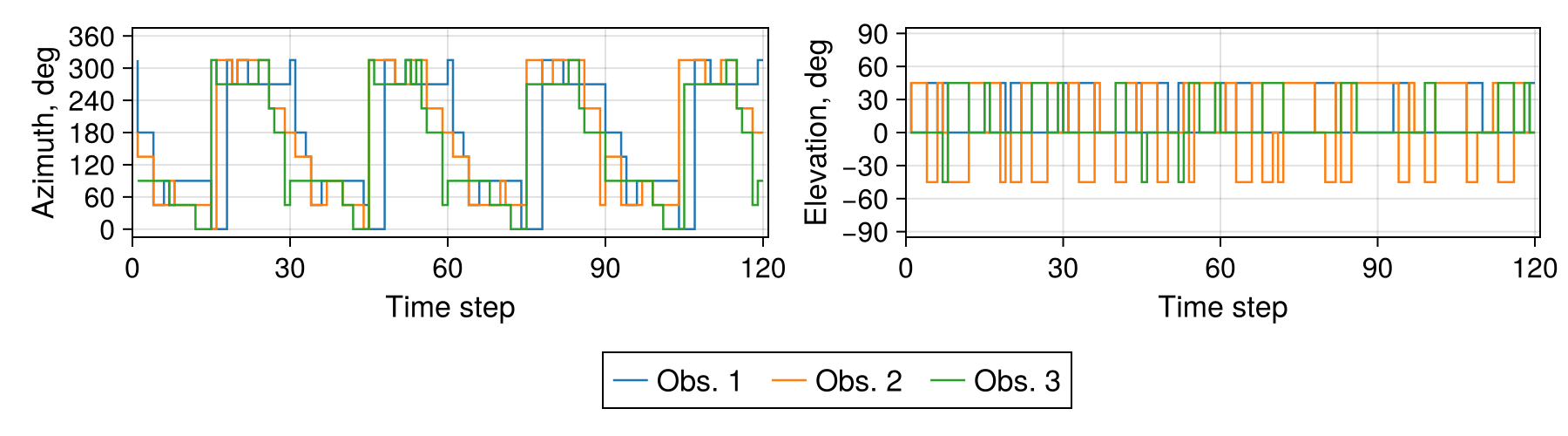}
        \caption{Lagrangian method solution with $\bar{m}_{\mathrm{crit}}=20$, $p=3$}
        \label{fig:lagrangean_targetcone1_p3_Mcut20_fov60}
    \end{subfigure}
    \\
    \begin{subfigure}[b]{0.99\textwidth}
        \centering
        \includegraphics[width=\textwidth]{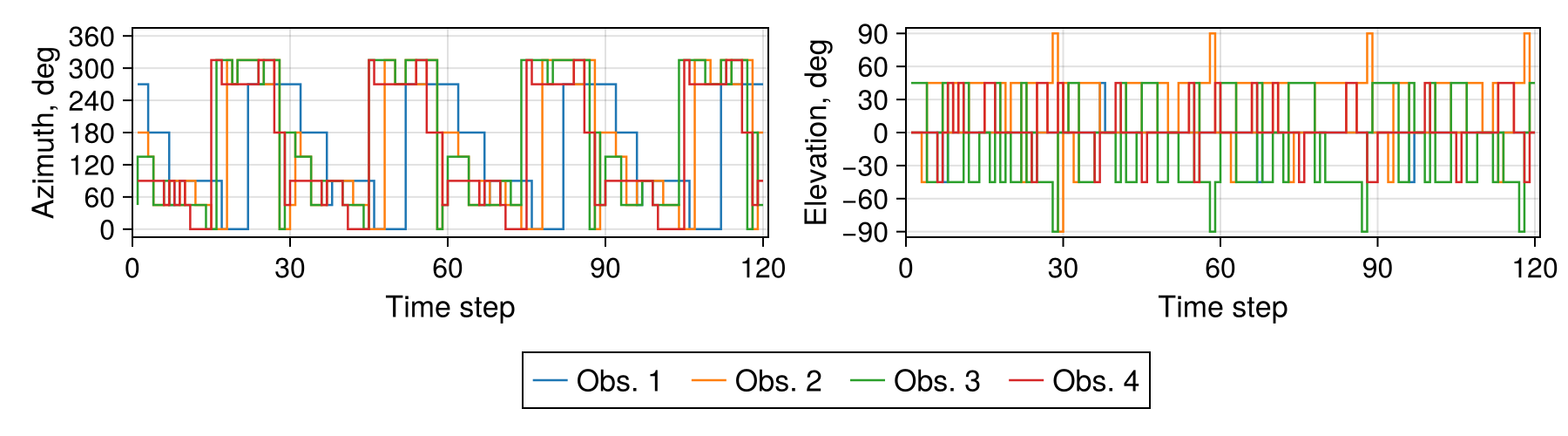}
        \caption{Lagrangian method solution with $\bar{m}_{\mathrm{crit}}=20$, $p=4$}
        \label{fig:lagrangean_targetcone1_p4_Mcut20_fov60}
    \end{subfigure}
    \\
    \begin{subfigure}[b]{0.99\textwidth}
        \centering
        \includegraphics[width=\textwidth]{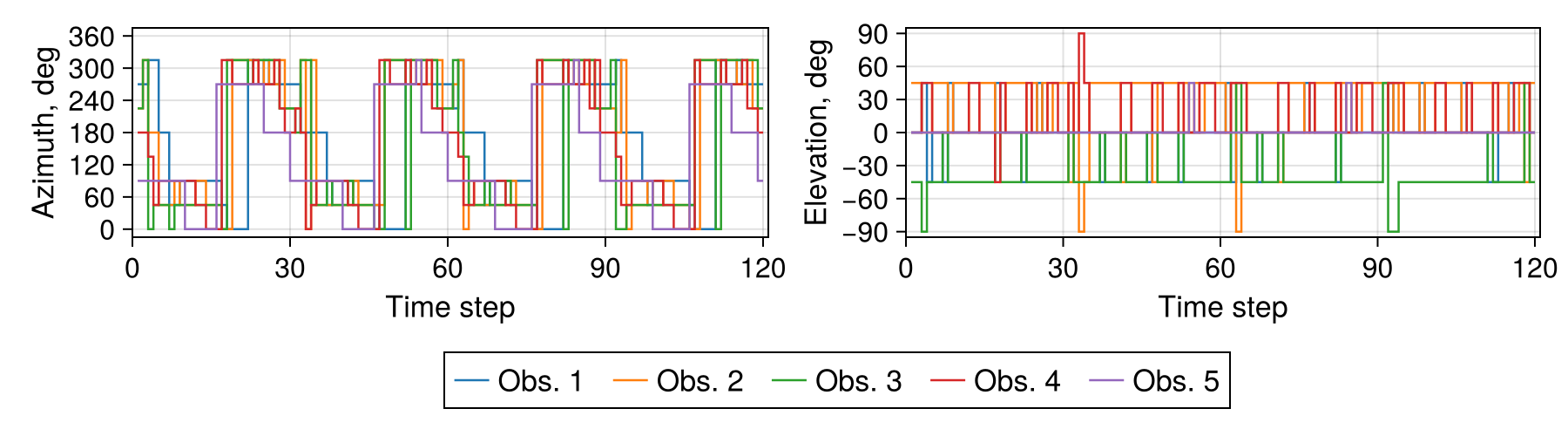}
        \caption{Lagrangian method solution with $\bar{m}_{\mathrm{crit}}=20$, $p=5$}
        \label{fig:lagrangean_targetcone1_p5_Mcut20_fov60}
    \end{subfigure}
    \caption{Optimized observer pointing direction histories for Cone of Shame demand}
    \label{fig:pointingdir_traj_cone1}
\end{figure}

\begin{figure}
    \centering
    \begin{subfigure}[b]{0.45\textwidth}
        \centering
        \includegraphics[width=\textwidth]{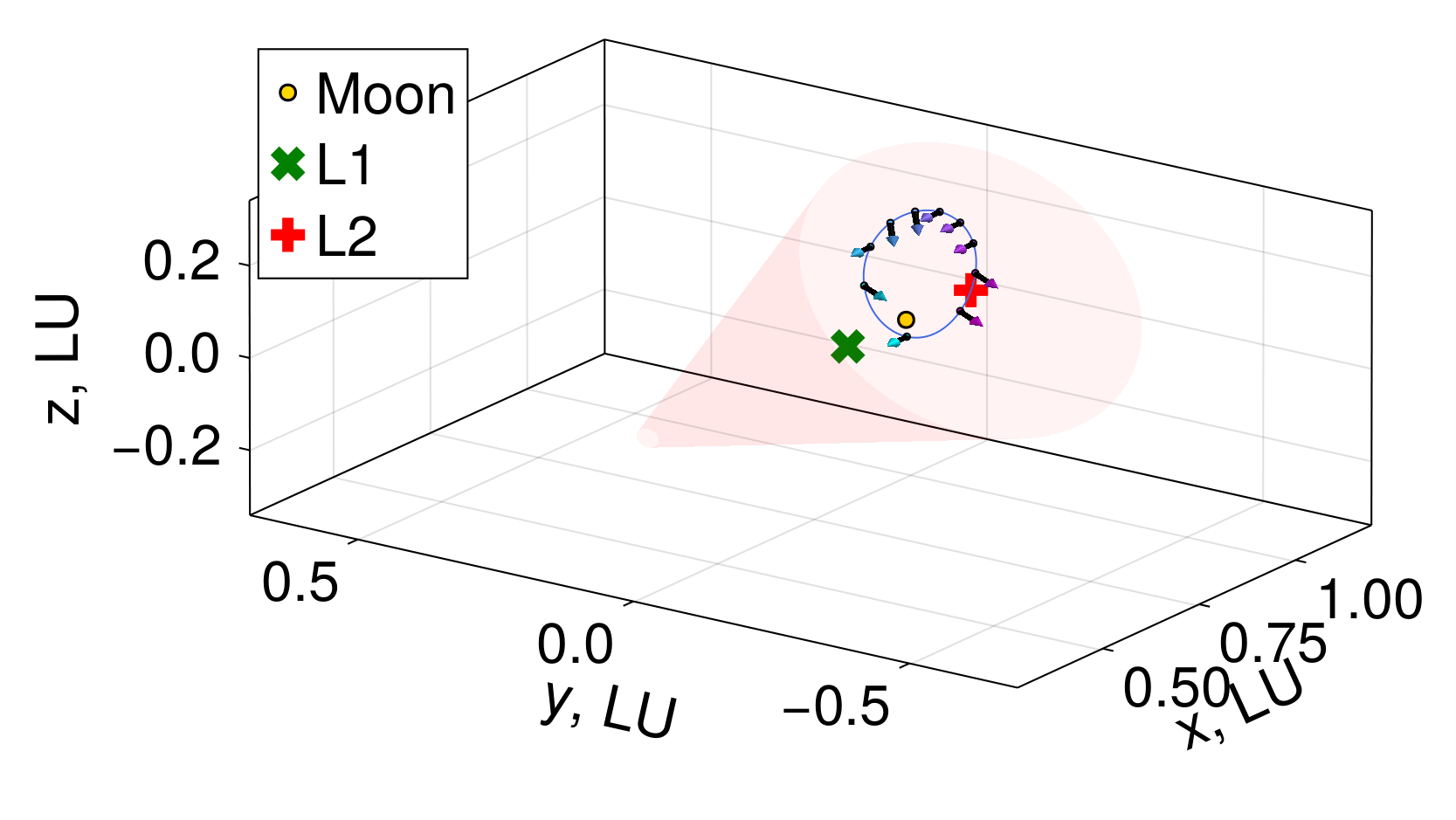}
        \caption{Observer 1 on 3:1 resonant Northern L2 Halo}
        \label{fig:steer_3d_lagrangean_targetcone1_p2_Mcut20_fov60_observer1}
    \end{subfigure}
    \hfill 
    \begin{subfigure}[b]{0.45\textwidth}
        \centering
        \includegraphics[width=\textwidth]{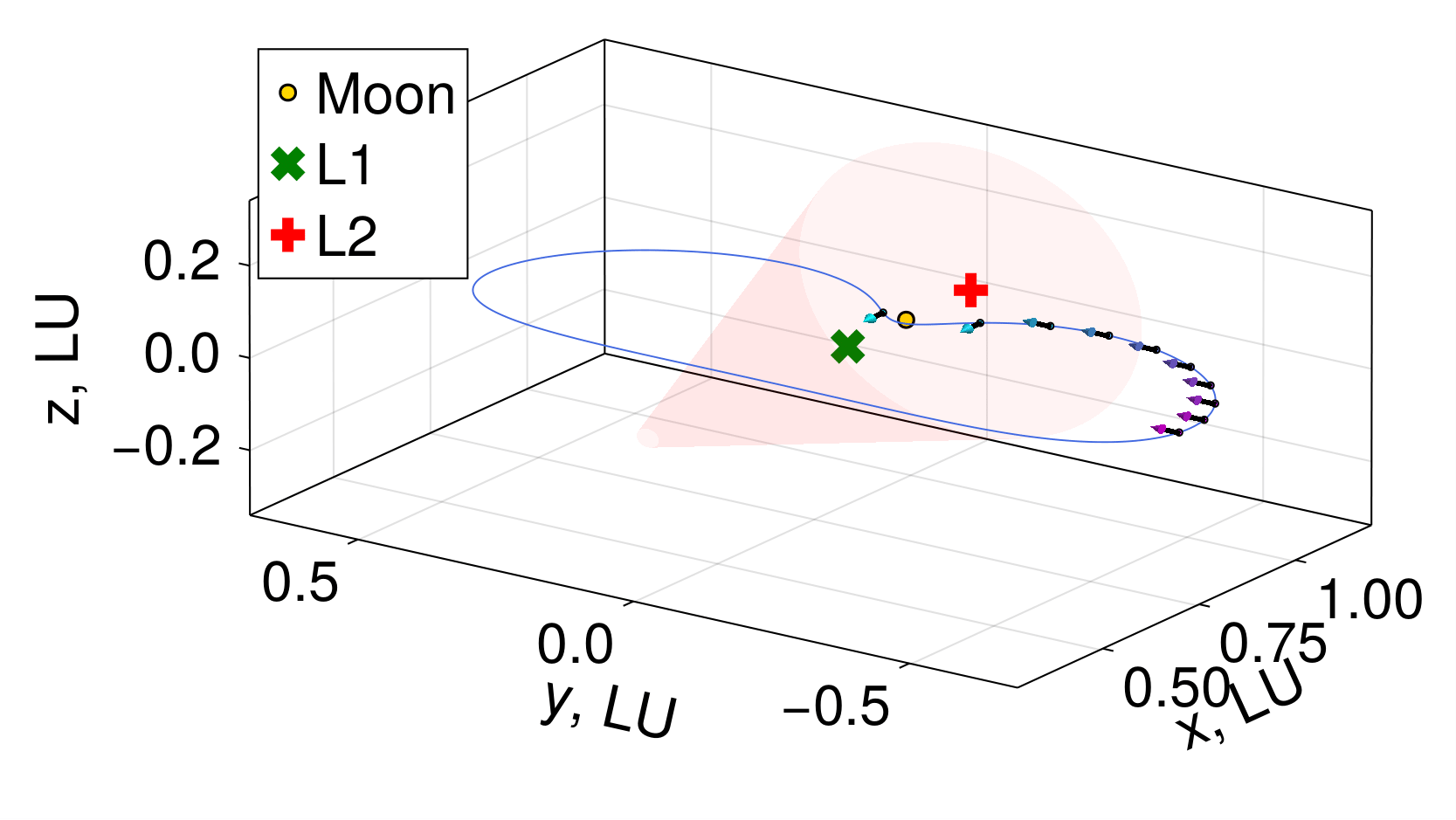}
        \caption{Observer 2 on 1:1 resonant L1 Lyapunov}
        \label{fig:steer_3d_lagrangean_targetcone1_p2_Mcut20_fov60_observer2}
    \end{subfigure}
    \caption{Sensor pointing directions over 10 time-steps with $\bar{m}_{\rm crit} = 20$ $p=2$ for Cone of Shame demand}
    \label{fig:steer_history_cone}
\end{figure}

The steering history for the selected architectures in Figure~\ref{fig:arch_traj_cone1} are shown in Figure~\ref{fig:pointingdir_traj_cone1}. 
For $p=2$, as shown in Figure~\ref{fig:lagrangean_targetcone1_p2_Mcut20_fov60}, observer 1 located on the 3:1 resonant Halo steers its sensor both along the azimuth and elevation angles, while observer 2 has its sensor lied on the $xy$-plane. 
The azimuth steering angle of the two observers exhibits periodicities according to their orbital periods. 
We note that observer 1 is not disadvantaged by the fact that it is not on a 1:1 resonance as its sensor is pointed out of plane with a non-zero elevation angle for the majority of the time.
To provide physical intuition, we plot the steering history of the $p=2$ architecture in position space in Figure~\ref{fig:steer_history_cone}. 
For the sake of clarity, we only show the first 10 steps, corresponding to one revolution along observer 1 located on a 3:1 resonant LPO. 
With architectures for $p=3,4,5$, we are able to see steering histories of observers on the same, 1:1 resonant L1 Lyapunov.
The azimuth angle history for each observer follows a similar track shifted according to their location along the LPO. 
The appropriate azimuth angle to steer the sensor on the 1:1 resonant LPO is largely dictated by the Sun's illumination direction, causing the observers to follow a similar in-plane steering strategy. 
Meanwhile, the elevation steering history is distinct for each observer; the out-of-plane sensor pointing direction is varied in order to achieve a more comprehensive coverage of the volume of interest.

The sensor steering history provides insights as to why in the $p=2$ case, the solution consists of one observer on the L1 Lyapunov orbit and another on a different LPO, while in cases with $p = 3,4,5$ all observers placed on the L1 Lyapunov orbit. 
With at least $p=3$ observers, each can be tasked to observe at elevation angles of $0^{\circ}$, $-30^{\circ}$, and $30^{\circ}$ - thereby complementing the observation to provide holistic coverage. 
Meanwhile, with only $p=2$ observers, the 3:1 resonant Halo, with higher cadence steering across azimuth and elevation angles, is found to be able to better complement the observation of a single L1 Lyapunov orbit observer. 
We note that the elevation angle spread for $p=3$ as seen in Figure~\ref{fig:arch_traj_cone1_Mcut20_p3_lm} is not as thorough as with cases with $p=4$ and $5$, as seen in Figures~\ref{fig:arch_traj_cone1_Mcut20_p4_lm} and~\ref{fig:arch_traj_cone1_Mcut20_p5_lm}. 
This is the reason why the LM solution is $-15.79\%$ worse than the B\&B solution after 3600 seconds, as reported in Table~\ref{tab:SSA_pmedian_gurboi_vs_LR_cone1_m20_fov60}, compared to cases with $p=2,4,5$.

\subsubsection{Low-Energy Transfer Transit Window}
\begin{figure}
    \centering
    \begin{subfigure}[b]{0.45\textwidth}
        \centering
        \includegraphics[width=\textwidth]{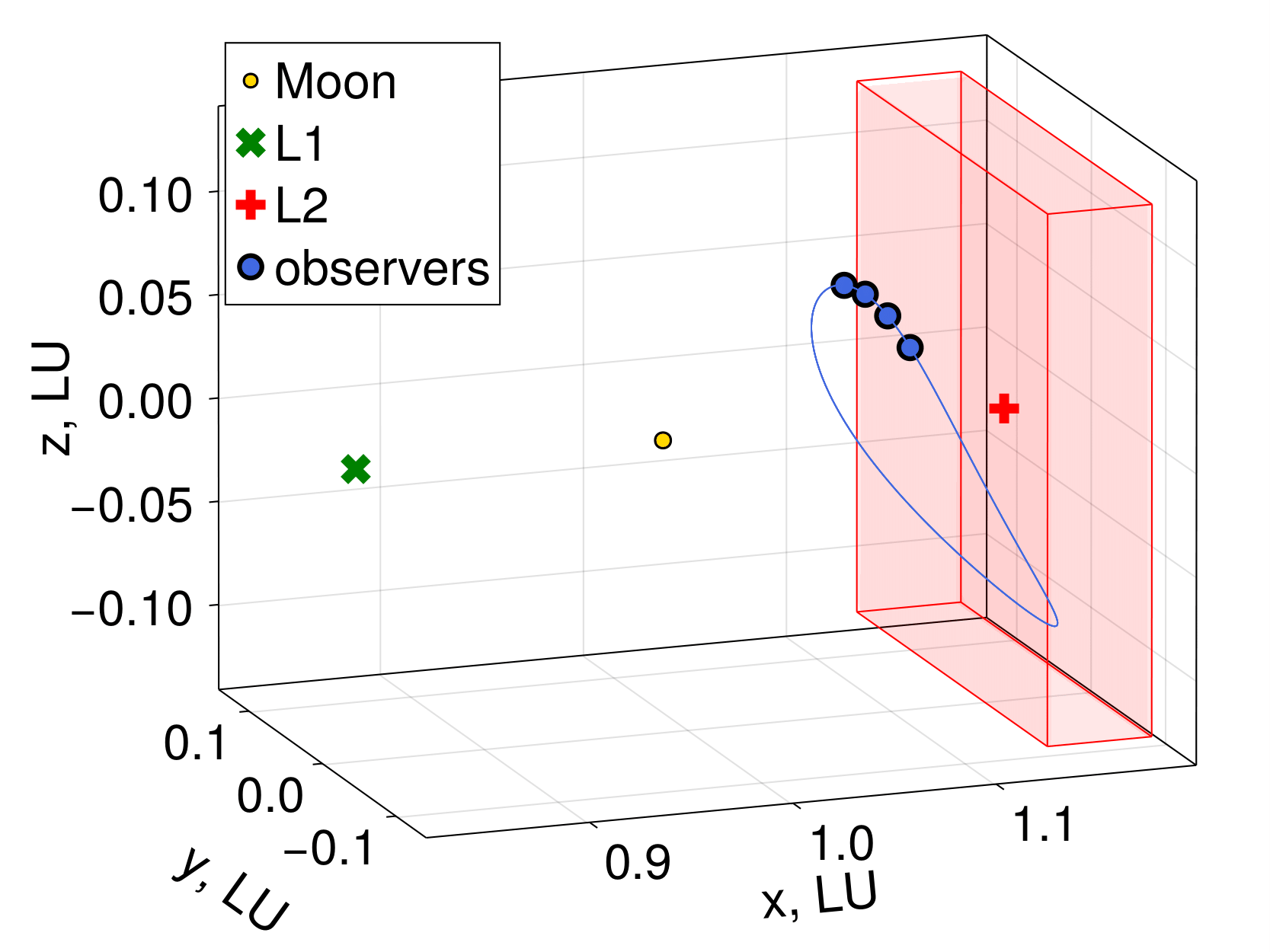}
        \caption{Lagrangian method solution with $\bar{m}_{\mathrm{crit}}=15$, $p=4$}
        \label{fig:arch_traj_let1_Mcut15_p4_lm}
    \end{subfigure}
    \hfill
    \begin{subfigure}[b]{0.45\textwidth}
        \centering
        \includegraphics[width=\textwidth]{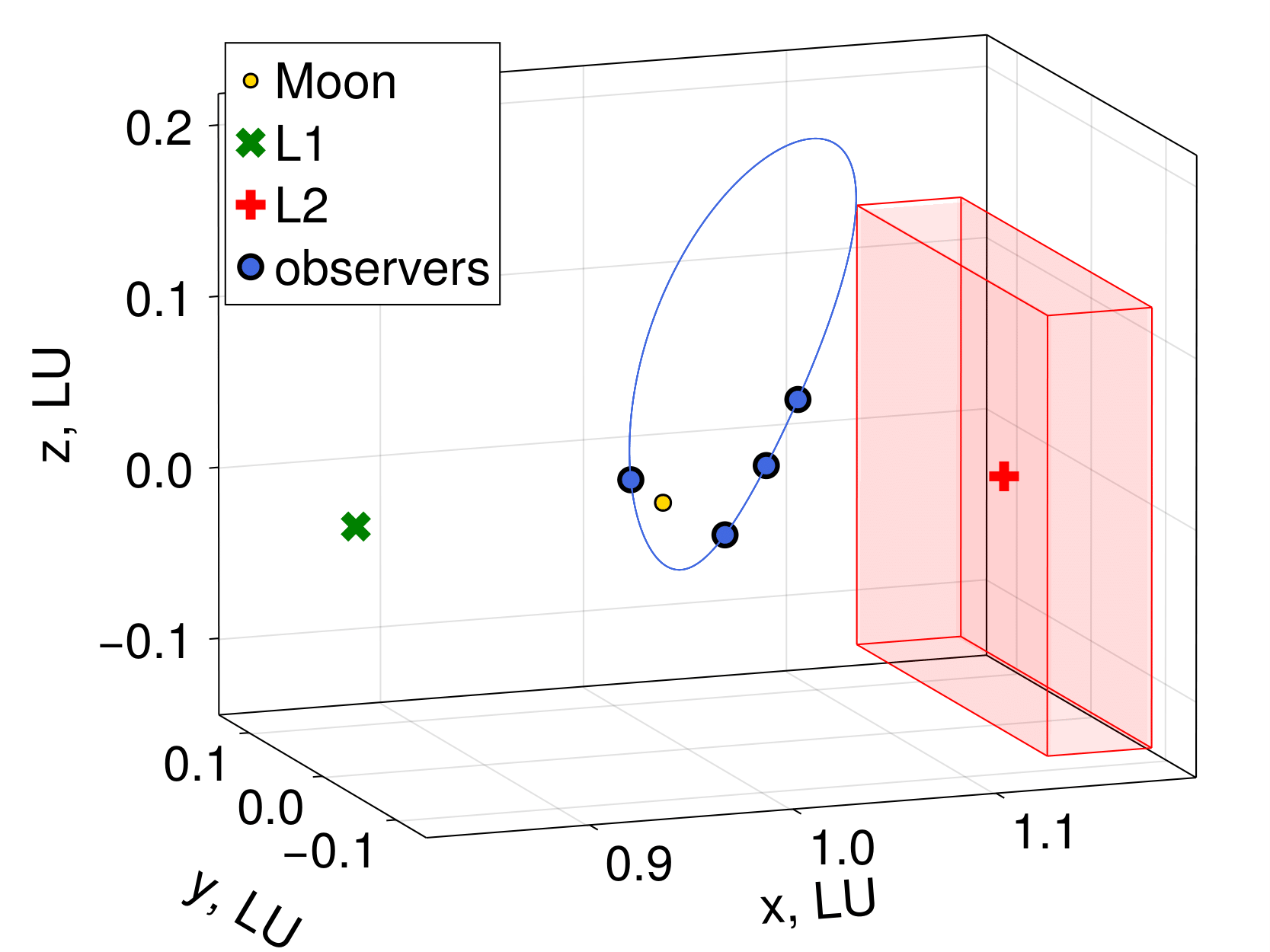}
        \caption{Lagrangian method solution with $\bar{m}_{\mathrm{crit}}=18$, $p=4$}
        \label{fig:arch_traj_let1_Mcut18_p4_lm}
    \end{subfigure}
    \\
    \begin{subfigure}[b]{0.45\textwidth}
        \centering
        \includegraphics[width=\textwidth]{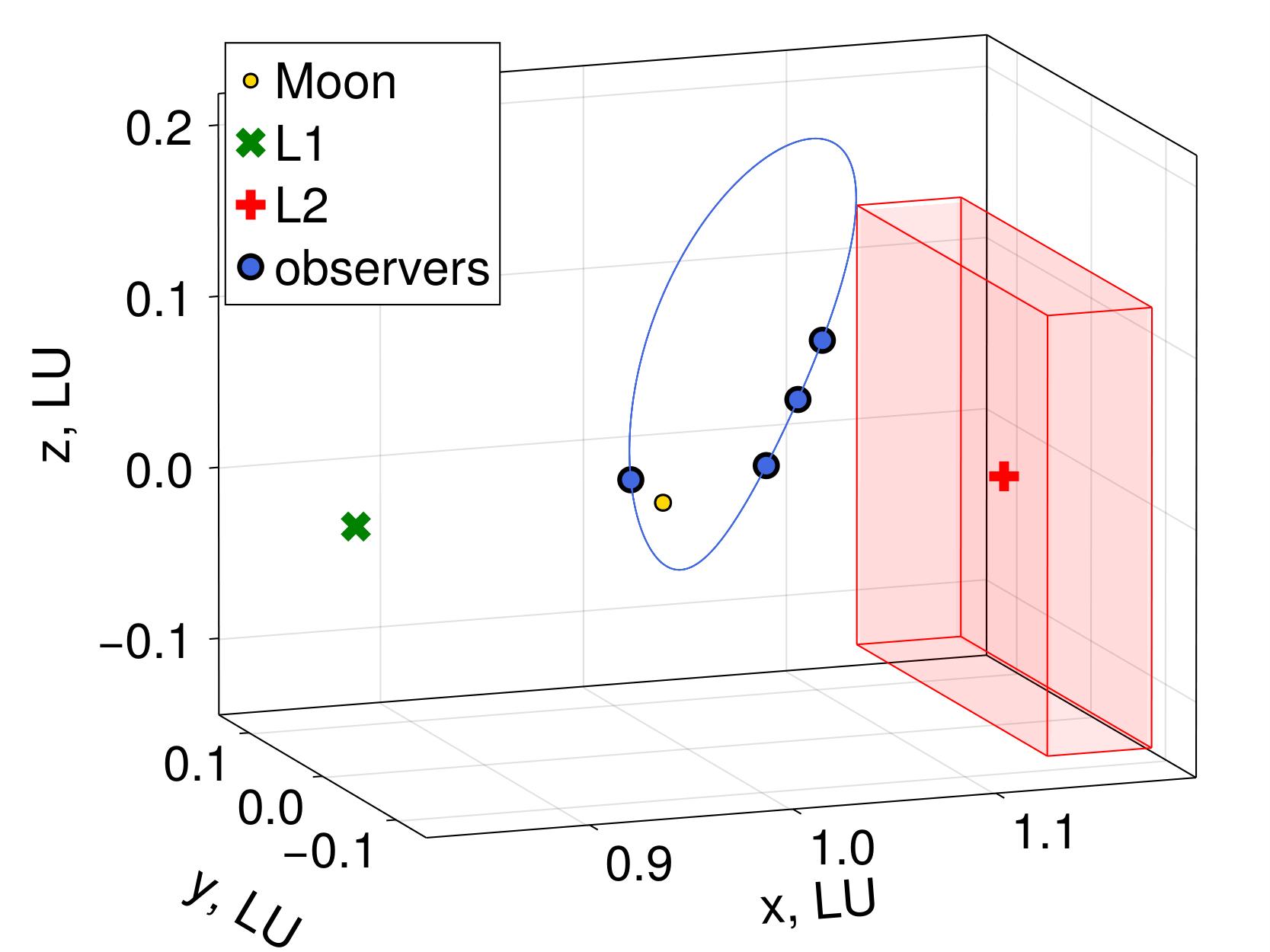}
        \caption{Lagrangian method solution with $\bar{m}_{\mathrm{crit}}=20$, $p=4$}
        \label{fig:arch_traj_let1_Mcut20_p4_lm}
    \end{subfigure}
    \hfill
    \begin{subfigure}[b]{0.45\textwidth}
        \centering
        \includegraphics[width=\textwidth]{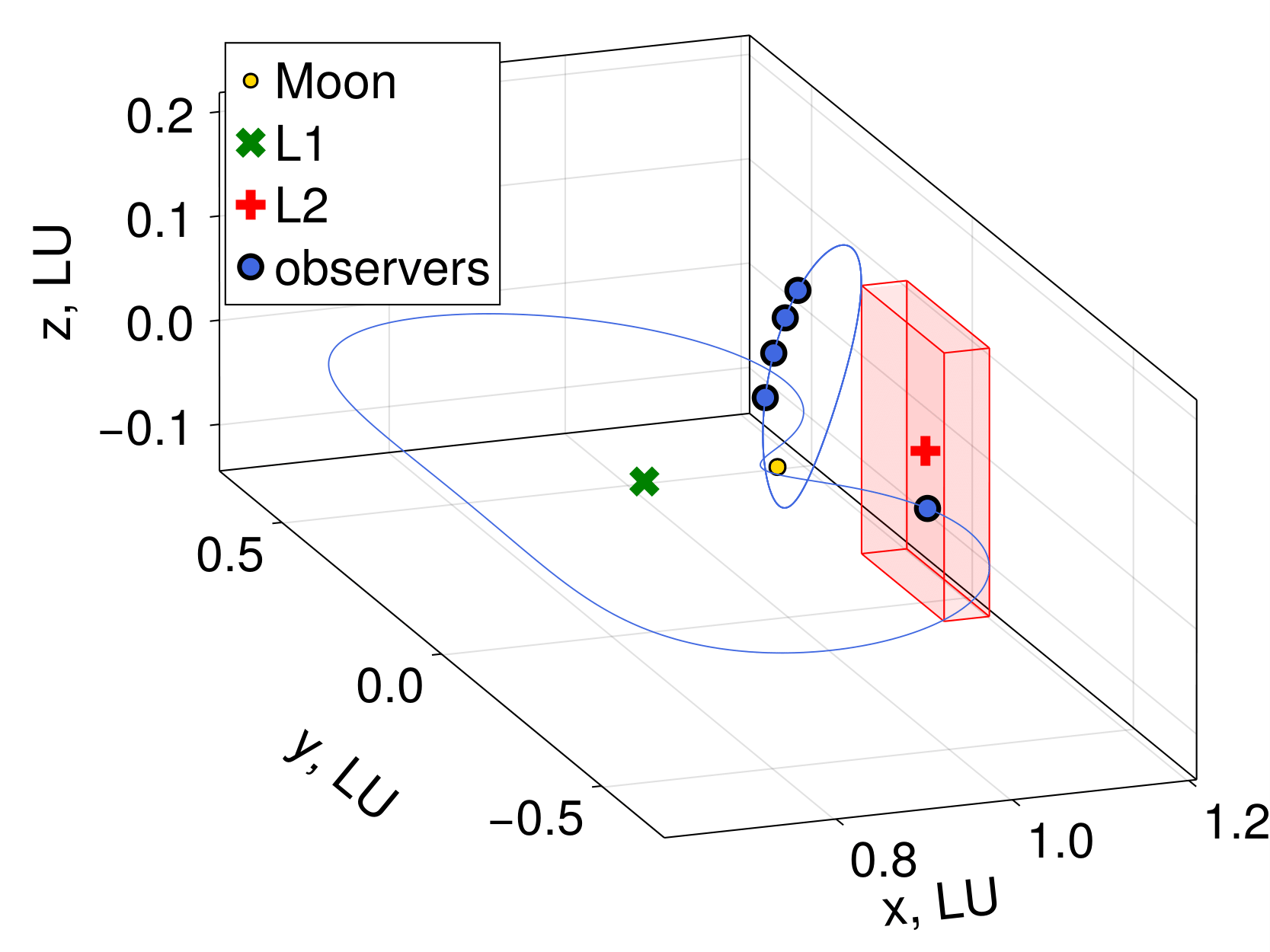}
        \caption{Gurobi solution with $\bar{m}_{\mathrm{crit}}=20$, $p=5$}
        \label{fig:arch_traj_let1_Mcut20_p5_lm}
    \end{subfigure}
    \caption{Optimized observer locations for LET transit demand at time-step $t=1$ (at a new Moon)}
    \label{fig:arch_traj_let1}
\end{figure}

We first compare solutions with increasing cut-off magnitude from 15 to 20, with observers fixed to $p=4$; the corresponding architectures, obtained by the LM, are shown in Figures~\ref{fig:arch_traj_let1_Mcut15_p4_lm},~\ref{fig:arch_traj_let1_Mcut18_p4_lm}, and~\ref{fig:arch_traj_let1_Mcut20_p4_lm}. When the cut-off magnitude is lower, targets become harder to observe, and the observers are placed in a 2:1 resonant L2 Halo orbit that is closer to the set of targets; this comes with the disadvantage that the observer is able to fit fewer targets within its FOV at a given time. 
Meanwhile, the solutions at $\bar{m}_{\mathrm{crit}} = 18,20$ use a 3:1 resonant L2 Halo orbit that is further away, offering larger coverage from a single sensor pointing allocation. 

We also recognize that the train-like constellation, reminiscent of solutions from Figure~\ref{fig:arch_traj_cone1}, is used on the 2:1 resonant Halo for $\bar{m}_{\mathrm{crit}} = 15$, and for three of the four observers on the 3:1 resonant Halo for $\bar{m}_{\mathrm{crit}} = 18,20$. 
The fourth observer in the last two cases is placed with a significant angular separation along the LPO to complement the coverage; note that because the LPO is not planar, the illumination conditions are not necessarily poor for this last observer during the times along the LPO where there is a non-negligible out-of-plane component. 

With $p=5$, the architecture consists of four 3:1 L2 Halo observers and an additional observer on the 1:1 L1 Lyapunov, as shown in Figure~\ref{fig:arch_traj_let1_Mcut20_p5_lm}. 
Now, all four observers on the 3:1 L2 Halo form a train-like configuration, all in adjacent angular locations, and a complementing observer is now placed on a different LPO altogether. Note that the 1:1 resonant L1 Lyapunov is a favorable LPO for placing a single additional observer due to the invariance of the illumination condition; while the larger orbit results in an increased distance from the target set, this does not constitute a penalty in this configuration, where $\bar{m}_{\mathrm{crit}}$ is set high. 

\begin{figure}
    \centering
    \begin{subfigure}[b]{0.99\textwidth}
        \centering
        \includegraphics[width=\textwidth]{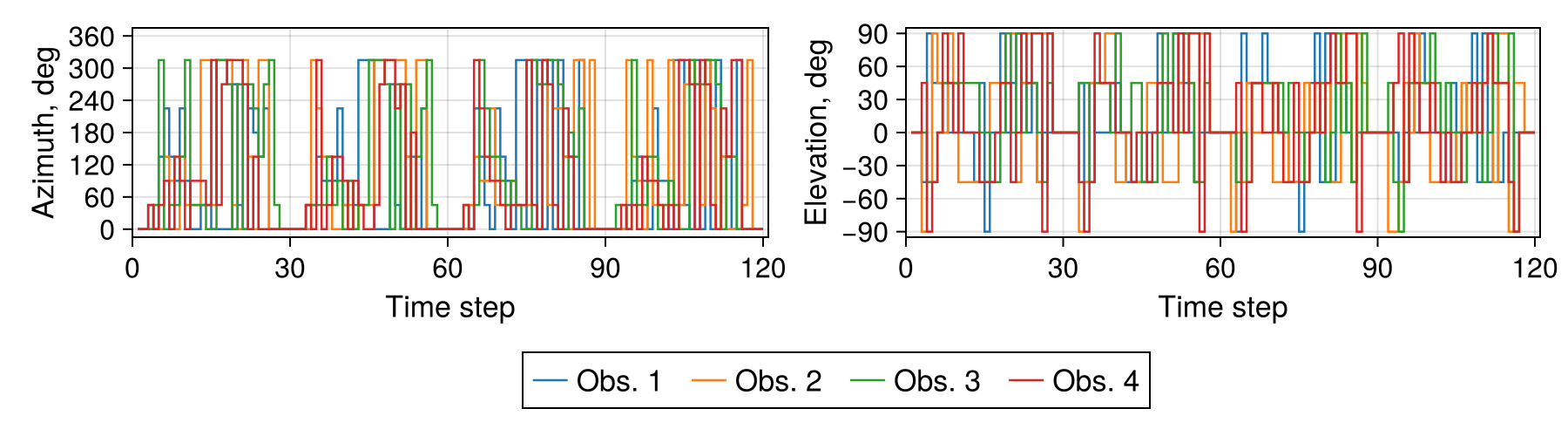}
        \caption{Lagrangian method solution with $\bar{m}_{\mathrm{crit}}=20$, $p=3$}
        \label{fig:lagrangean_targetlet1_p4_Mcut15_fov60}
    \end{subfigure}
    \hfill
    \begin{subfigure}[b]{0.99\textwidth}
        \centering
        \includegraphics[width=\textwidth]{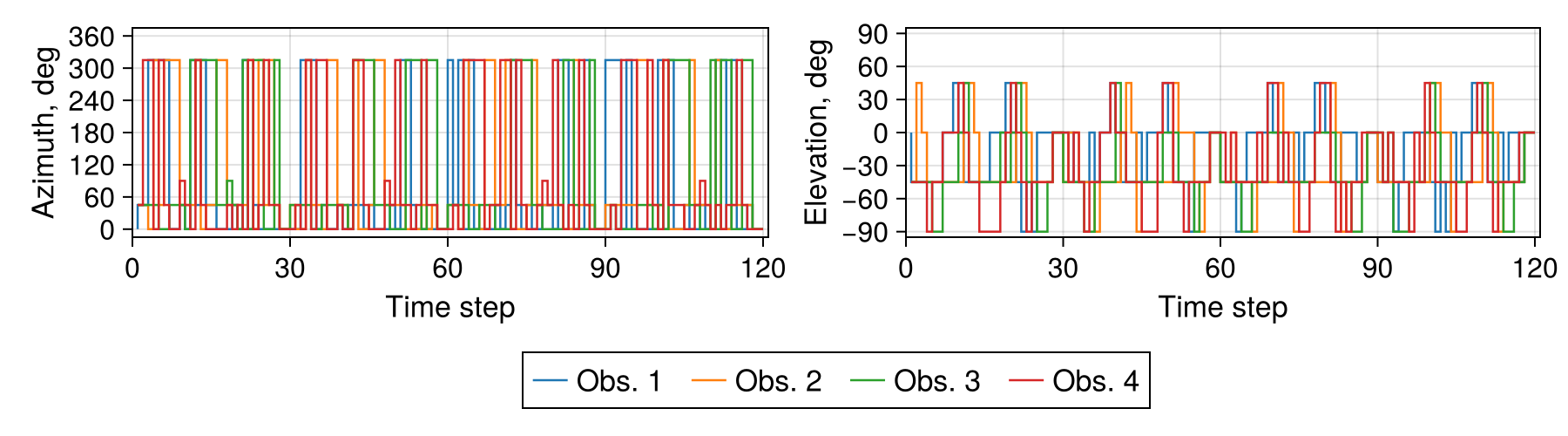}
        \caption{Lagrangian method solution with $\bar{m}_{\mathrm{crit}}=20$, $p=4$}
        \label{fig:lagrangean_targetlet1_p4_Mcut18_fov60}
    \end{subfigure}
    \\
    \begin{subfigure}[b]{0.99\textwidth}
        \centering
        \includegraphics[width=\textwidth]{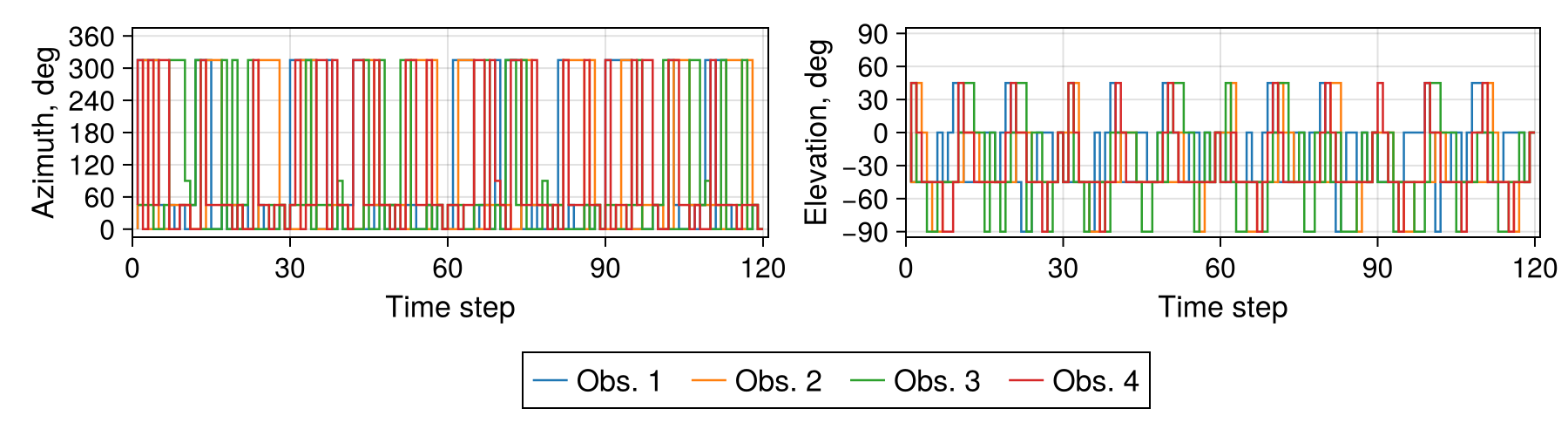}
        \caption{Lagrangian method solution with $\bar{m}_{\mathrm{crit}}=20$, $p=5$}
        \label{fig:lagrangean_targetlet1_p4_Mcut20_fov60}
    \end{subfigure}
    \\
    \begin{subfigure}[b]{0.99\textwidth}
        \centering
        \includegraphics[width=\textwidth]{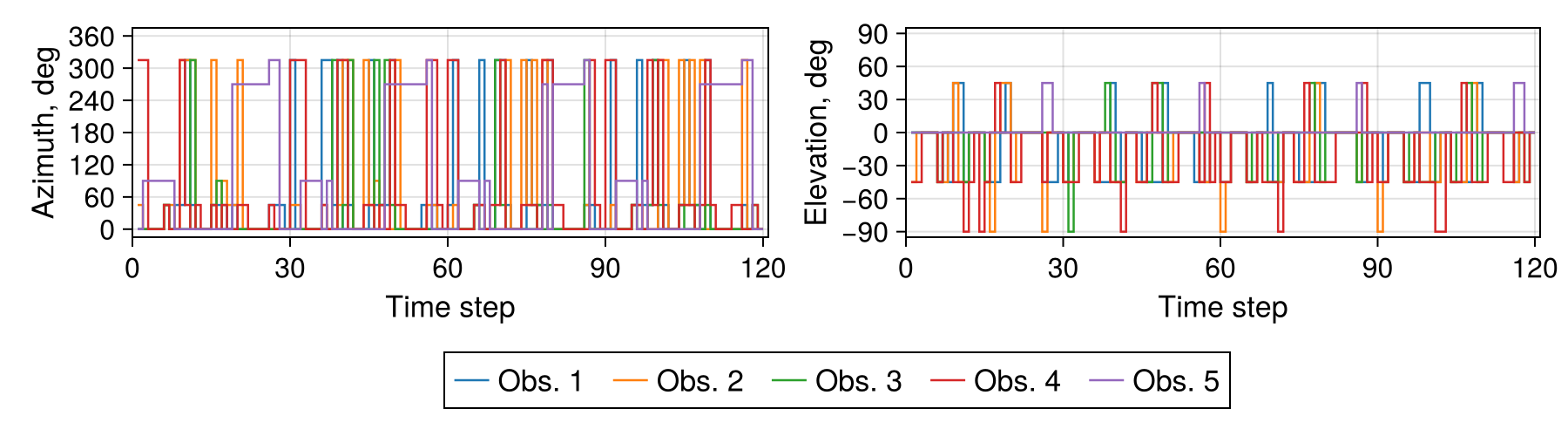}
        \caption{Lagrangian method solution with $\bar{m}_{\mathrm{crit}}=20$, $p=5$}
        \label{fig:lagrangean_targetlet1_p5_Mcut20_fov60}
    \end{subfigure}
    \caption{Optimized observer pointing direction histories for LET transit demand}
    \label{fig:pointingdir_traj_let1}
\end{figure}

\begin{figure}
    \centering
    \begin{subfigure}[b]{0.45\textwidth}
        \centering
        \includegraphics[width=\textwidth]{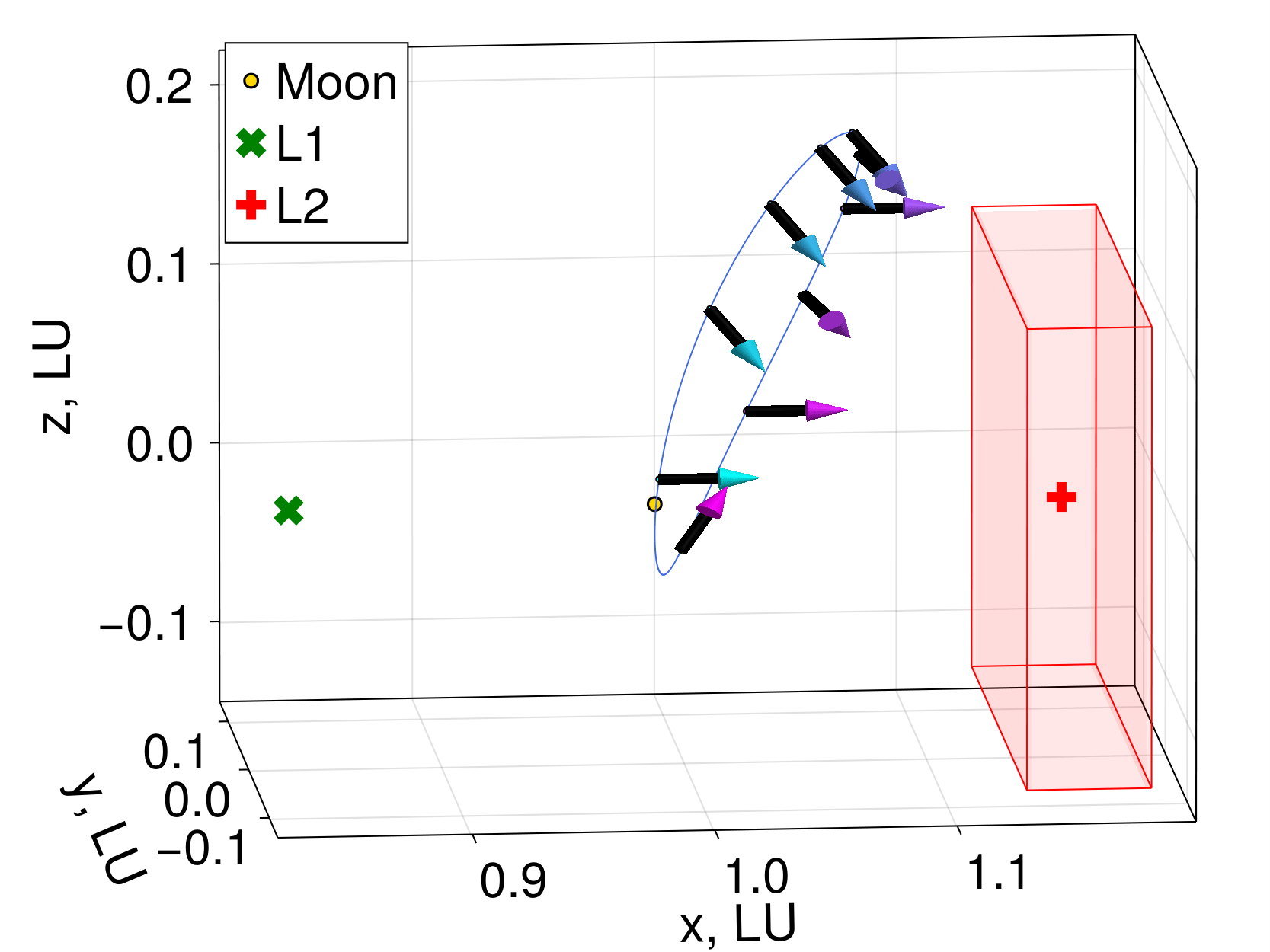}
        \caption{Observer 1 on 3:1 resonant Northern L2 Halo}
        \label{fig:steer_3d_lagrangean_targetlet1_p4_Mcut20_fov60_observer1}
    \end{subfigure}
    \hfill 
    \begin{subfigure}[b]{0.45\textwidth}
        \centering
        \includegraphics[width=\textwidth]{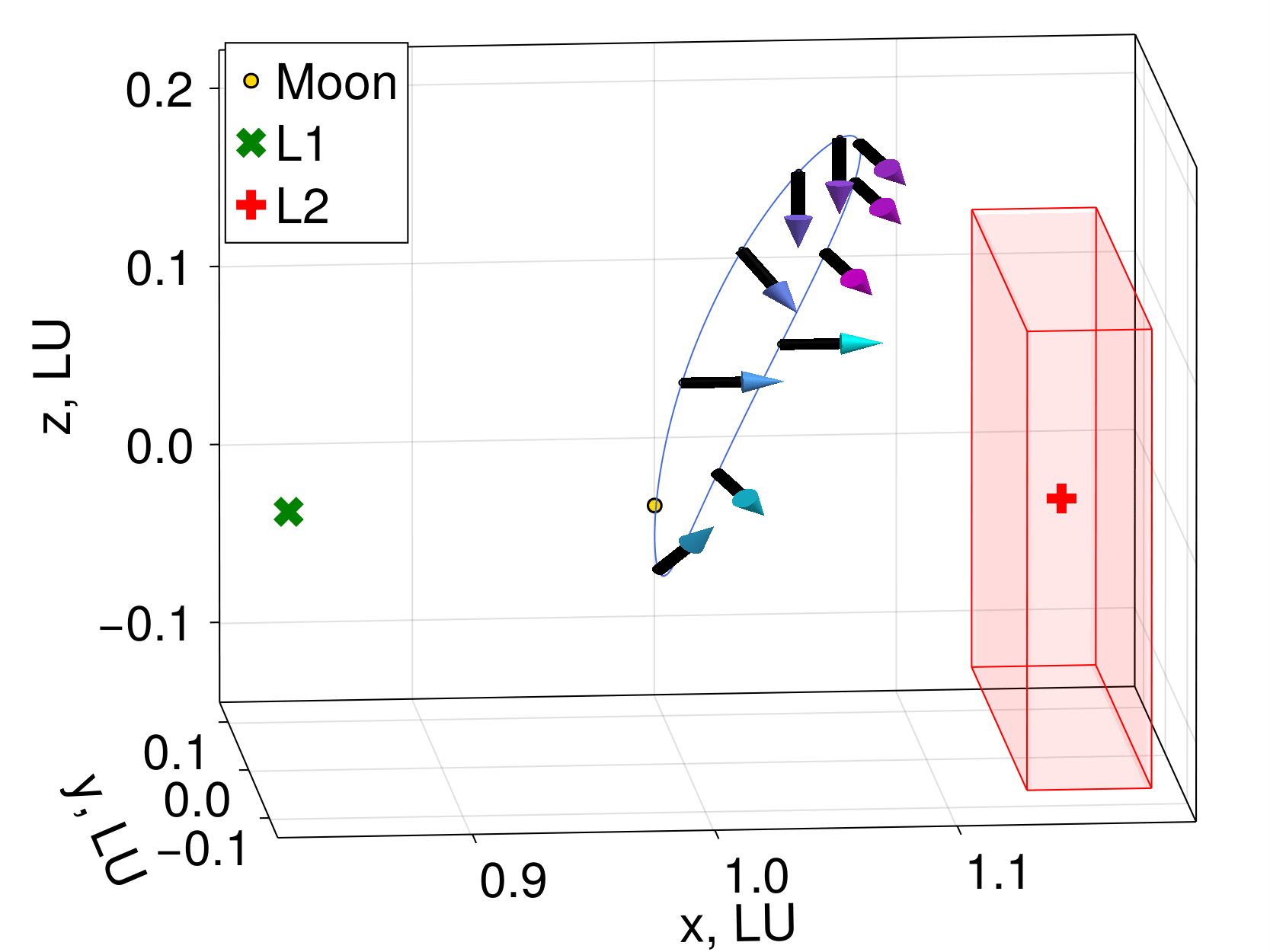}
        \caption{Observer 2 on 3:1 resonant Northern L2 Halo}
        \label{fig:steer_3d_lagrangean_targetlet1_p4_Mcut20_fov60_observer2}
    \end{subfigure}
    \\
    \begin{subfigure}[b]{0.45\textwidth}
        \centering
        \includegraphics[width=\textwidth]{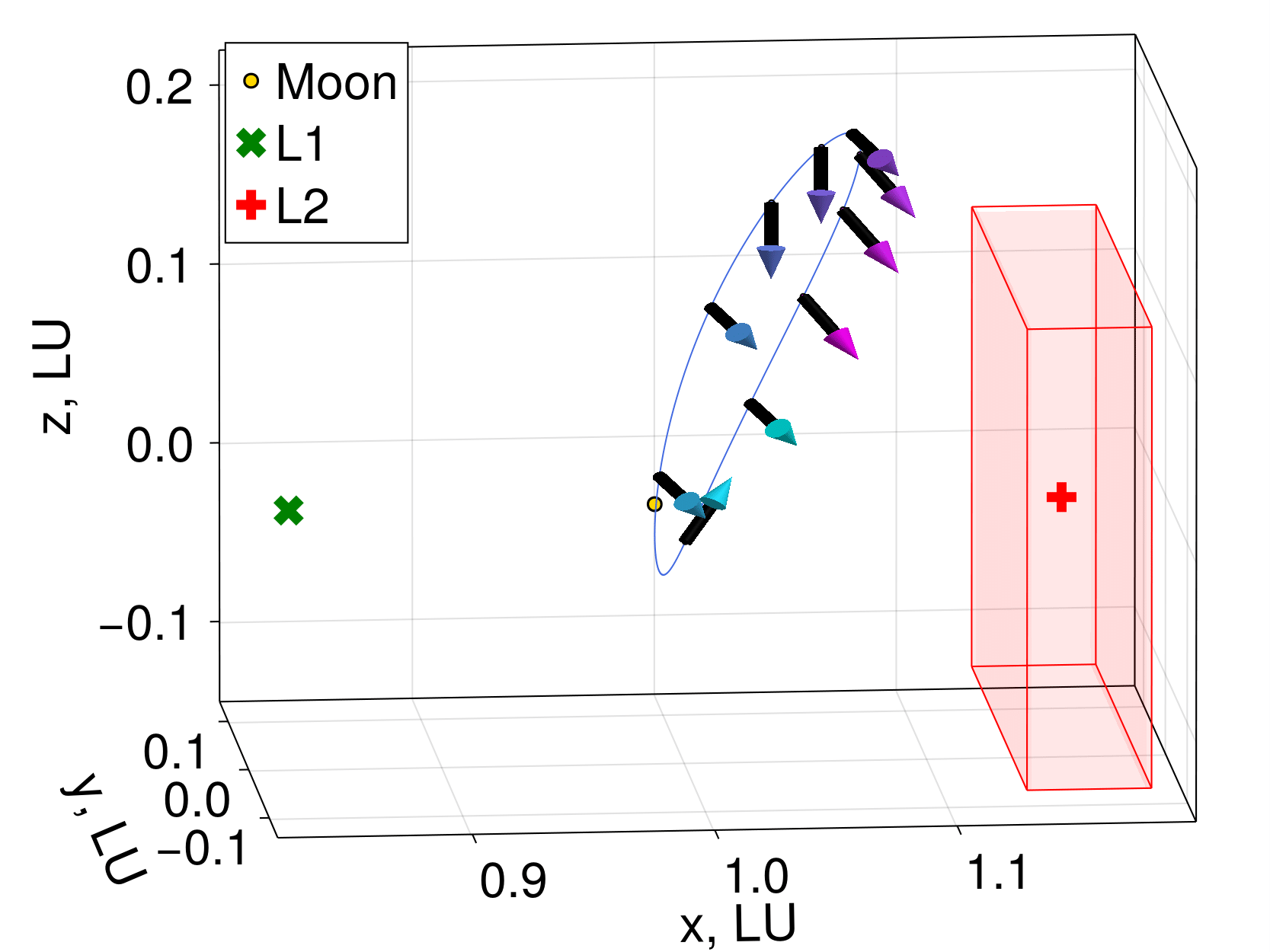}
        \caption{Observer 3 on 3:1 resonant Northern L2 Halo}
        \label{fig:steer_3d_lagrangean_targetlet1_p4_Mcut20_fov60_observer3}
    \end{subfigure}
    \hfill 
    \begin{subfigure}[b]{0.45\textwidth}
        \centering
        \includegraphics[width=\textwidth]{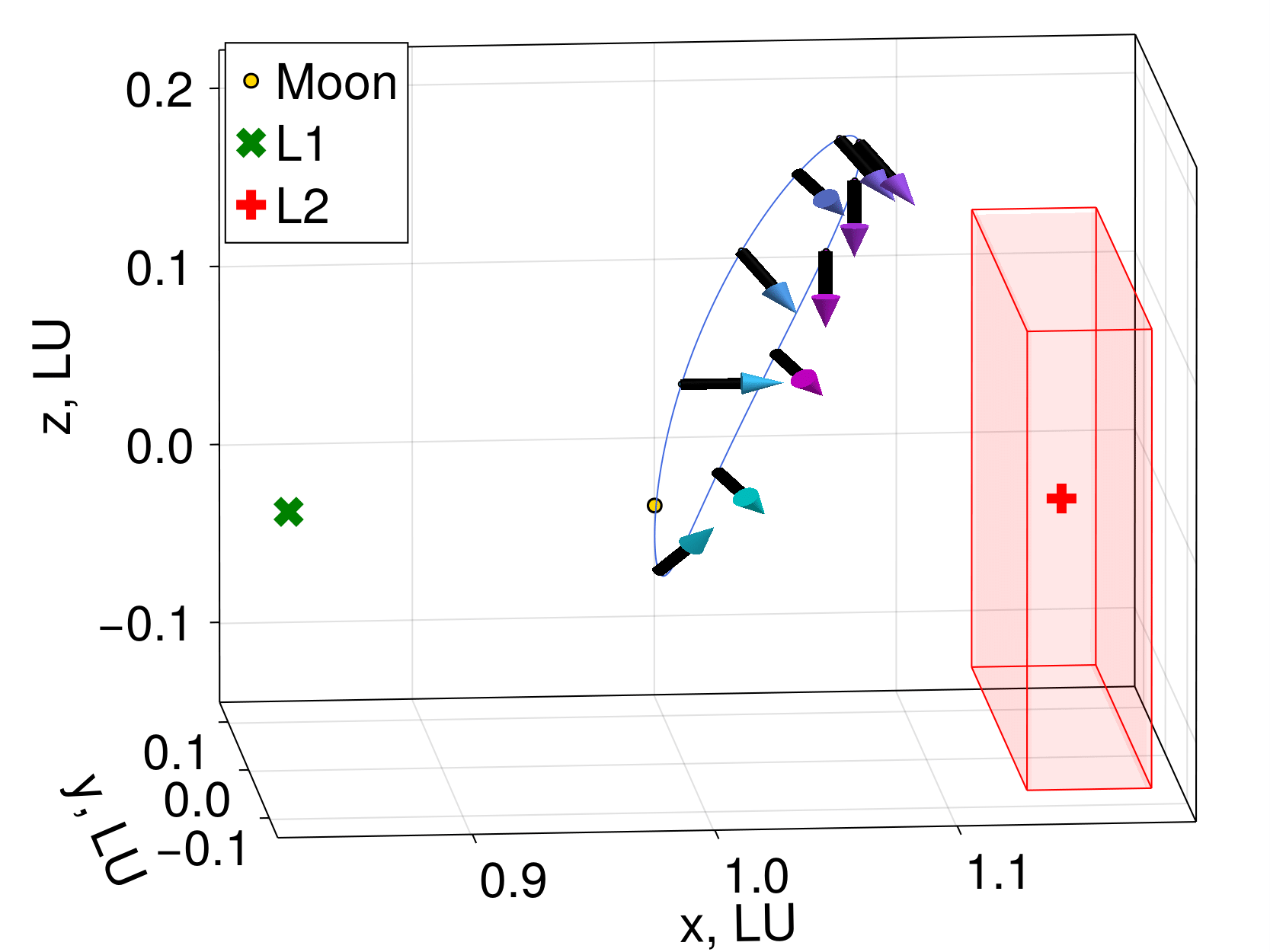}
        \caption{Observer 4 on 3:1 resonant Northern L2 Halo}
        \label{fig:steer_3d_lagrangean_targetlet1_p4_Mcut20_fov60_observer4}
    \end{subfigure}
    \caption{Sensor pointing directions over 10 time-steps with $\bar{m}_{\rm crit} = 20$ $p=4$ for LET transit demand}
    \label{fig:steer_history_let}
\end{figure}

{
Figure~\ref{fig:pointingdir_traj_let1} shows the sensor steering history for architectures in Figure~\ref{fig:arch_traj_let1}. 
In all cases, all observers are located on one side of the region of interest; as a result, the azimuth angle remains mostly between $-30^{\circ}$ and $30^{\circ}$. 
With the $p=2$ case utilizing Southern halos, the elevation angle across observers and across time is skewed on the positive side, while with the $p=3,4,5$ cases utilizing Northern halos, the elevation angle is skewed on the negative side. 
Figure~\ref{fig:steer_history_let} shows the sensor pointing directions for $\bar{m}_{\rm crit} = 20$ and $p = 4$ from Figure~\ref{fig:arch_traj_let1_Mcut20_p4_lm}. 
As the Northern halo has a larger out-of-plane amplitude than the region of interest, sensors are pointed downward around the apolunes along the Northern halo. 
We also note that over the first 10 time steps, the steering directions of each observer have slight variations from one another due to the variation in their respective phase along the LPO. 
}

\subsection{Key Takeaways on the Architecture Design Space}
{Taking the analyses on the obtained architectures into account, the following general remarks can be made: 
\begin{itemize}
    \item For a set of demand spread across L1 and L2, the large 1:1 resonant L1 Lyapunov offers an advantageous observation location; when the sensor is insufficient to observe targets at long ranges, the solution favors an alternative LPO that resides closer to certain regions within the volume of interest where most targets reside.
    \item Due to the solar phase angle largely dictating the visibility of targets, observers on 1:1 resonant LPOs tend to be in a leading/trailing configuration rather than being evenly spaced along the LPO. 
    As expected, the pointing directions are chosen in a manner that each observer complements the volume coverage of one another. 
    \item For a set of demands centered around L2, LPOs covering a smaller spatial region, in particular 2:1 and 3:1 resonant L2 halos, are favored. 
    These LPOs are suitable as the observer fly in the region near but outside the volume of interest; the latter consideration is crucial when considering FOV, as residing \textit{inside} the volume of interest would severely limit the fraction of volume that is visible to a given observer, despite its proximity to targets.
\end{itemize}
}

\section{Conclusions}
\label{sec:conclusions}
In this work, the TE-$p$-MP has been proposed for CSSA. 
The formulation consists of determining the $p$ optimal observer locations along with the orientation of their optical sensors over time to observe a given set of static or dynamic targets. 
{The TE-$p$-MP, inspired by the classical $p$-Median problem, is a BLP discretized in time, candidate locations of observer spacecraft, and sensor pointing directions. 
This formulation simultaneously solves for the constellation design and the sensor-tasking problems.}
{The formulation has been developed together with a} custom LM algorithm, based on the relaxation of complicating constraints and custom heuristics leveraging the problem structure as well as the geometry of the candidate facility locations. 
{A solution to the TE-$p$-MP is a CSSA constellation together with sensor-tasking for each observer that observes a prescribed set of targets.}

The TE-$p$-MP, along with the LM, has been demonstrated to find near-optimal solutions for various CSSA scenarios varying in terms of the distribution of demands or observer parameters in the order of minutes. 
The primary benefit of the LM is the reduction in solve time compared to using generic BLP solvers, which enables exploring the trade space of parameters such as the visible magnitude threshold or the temporal and spatial distribution of observation targets.
The LM has been demonstrated to perform in a largely consistent manner across most TE-$p$-MP instances and to deliver competitive near-optimal constellations, whereas B\&B's performance is found to be far more precarious, depending on the specific parameters of the instance. 
Several noteworthy constellations tailored to monitoring demands based on the Cone of Shame and a dynamic transit window of LETs have been analyzed to demonstrate the proposed approach.

Future works include developing a synergistic scheme to use the LM with B\&B algorithms, thus combining LM's capability of quickly generating feasible, competitive solutions, with the B\&B's guaranteed capability of closing the optimality gap. 
The objective of the TE-$p$-MP may also be modified to account for the cost of changing the sensor's pointing direction.
Finally, while this work specifically explored the use of the TE-$p$-MP for CSSA, a similar formulation may be conceivable for other cislunar constellation-based applications, such as the design of Cislunar on-orbit servicing or GNSS constellations.

\section*{Acknowledgments}
Artificial Intelligence technologies were used for grammar checking.

\section*{Appendix}

\subsection{Parameters of Synodic-Resonant LPOs}
\label{sec:appendix_lpo_conditions}
Parameters of the synodic-resonant LPOs considered as candidate locations are given in Table~\ref{tab:candidate_lpos_info}. For LPOs with Southern and Northern families, only values for the Southern member is provided; the initial state of the corresponding Northern member is obtained by flipping the sign of $x_0$.

\begin{table}[]
\centering
\caption{Parameters of synodic-resonant LPOs}
\begin{tabular}{@{}llllllll@{}}
\toprule
Family & $M$:$N$ & $x_0$ & $z_0$ & $\dot{y}_0$ & Period, TU & Stability & Slots $b$ \\
\midrule
DRO & 9:2 & 0.88976967 & 0 & 0.47183463 & 1.47892343 & 1.00 & 14 \\
DRO & 4:1 & 0.88060589 & 0 & 0.47011146 & 1.66378885 & 1.00 & 15 \\
DRO & 3:1 & 0.85378188 & 0 & 0.47696024 & 2.21838514 & 1.00 & 20 \\
DRO & 9:4 & 0.81807765 & 0 & 0.50559384 & 2.95784685 & 1.00 & 27 \\
DRO & 2:1 & 0.79946085 & 0 & 0.52703349 & 3.32757771 & 1.00 & 30 \\
DRO & 3:2 & 0.73370014 & 0 & 0.62889866 & 4.43677028 & 1.00 & 40 \\
DRO & 5:2 & 0.83249233 & 0 & 0.49184571 & 2.66206217 & 1.00 & 24 \\
L2 Halo (Southern) & 9:2 & 1.01958272 & -0.18036049 & -0.09788185 & 1.47892343 & 1.00 & 14 \\
L2 Halo (Southern) & 4:1 & 1.03352559 & -0.18903385 & -0.12699215 & 1.66378885 & 1.00 & 15 \\
L2 Halo (Southern) & 3:1 & 1.07203837 & -0.20182525 & -0.18853332 & 2.21838514 & 1.00 & 20 \\
L2 Halo (Southern) & 9:4 & 1.12518004 & -0.18195085 & -0.22544142 & 2.95784685 & 28.78 & 27 \\
L2 Halo (Southern) & 2:1 & 1.16846916 & -0.09994291 & -0.19568201 & 3.32757771 & 282.87 & 30 \\
L2 Halo (Southern) & 5:2 & 1.10193101 & -0.19829817 & -0.21702846 & 2.66206217 & 6.93 & 24 \\
L2 Halo (Northern) & 9:2 & 1.01958272 & 0.18036049 & -0.09788185 & 1.47892343 & 1.00 & 14 \\
L2 Halo (Northern) & 4:1 & 1.03352559 & 0.18903385 & -0.12699215 & 1.66378885 & 1.00 & 15 \\
L2 Halo (Northern) & 3:1 & 1.07203837 & 0.20182525 & -0.18853332 & 2.21838514 & 1.00 & 20 \\
L2 Halo (Northern) & 9:4 & 1.12518004 & 0.18195085 & -0.22544142 & 2.95784685 & 28.78 & 27 \\
L2 Halo (Northern) & 2:1 & 1.16846916 & 0.09994291 & -0.19568201 & 3.32757771 & 282.87 & 30 \\
L2 Halo (Northern) & 5:2 & 1.10193101 & 0.19829817 & -0.21702846 & 2.66206217 & 6.93 & 24 \\
DPO & 4:1 & 1.06189575 & 0 & 0.35989734 & 1.66378885 & 2.26 & 15 \\
DPO & 3:1 & 1.06335021 & 0 & 0.38222392 & 2.21838514 & 10.98 & 20 \\
DPO & 9:4 & 1.05547996 & 0 & 0.45941661 & 2.95784685 & 76.76 & 27 \\
DPO & 2:1 & 1.04880058 & 0 & 0.51457559 & 3.32757771 & 159.21 & 30 \\
DPO & 3:2 & 1.02851298 & 0 & 0.71048482 & 4.43677028 & 587.57 & 40 \\
DPO & 5:2 & 1.05978399 & 0 & 0.42240630 & 2.66206217 & 37.71 & 24 \\
DPO & 1:1 & 1.00515914 & 0 & 1.16888350 & 6.65515541 & 1399.19 & 59 \\
L1 Lyapunov & 9:4 & 0.81109465 & 0 & 0.26078428 & 2.95784685 & 746.89 & 27 \\
L1 Lyapunov & 2:1 & 0.79987674 & 0 & 0.35828602 & 3.32757771 & 407.88 & 30 \\
L1 Lyapunov & 3:2 & 0.76511295 & 0 & 0.49115556 & 4.43677028 & 133.00 & 40 \\
L1 Lyapunov & 1:1 & 0.63394833 & 0 & 0.79045684 & 6.65515541 & 53.98 & 59 \\
Butterfly (Northern) & 9:4 & 0.94130132 & -0.16165899 & -0.03565177 & 2.95784685 & 5.79 & 27 \\
Butterfly (Northern) & 2:1 & 0.91204757 & -0.14952514 & -0.02724245 & 3.32757771 & 12.45 & 30 \\
Butterfly (Northern) & 3:2 & 0.91414032 & -0.14492270 & -0.11588220 & 4.43677028 & 1.00 & 40 \\
Butterfly (Northern) & 1:1 & 0.99265217 & -0.17814460 & -0.26312433 & 6.65515541 & 1.00 & 59 \\
Butterfly (Southern) & 9:4 & 0.94130132 & 0.16165899 & -0.03565177 & 2.95784685 & 5.79 & 27 \\
Butterfly (Southern) & 2:1 & 0.91204757 & 0.14952514 & -0.02724245 & 3.32757771 & 12.45 & 30 \\
Butterfly (Southern) & 3:2 & 0.91414032 & 0.14492270 & -0.11588220 & 4.43677028 & 1.00 & 40 \\
Butterfly (Southern) & 1:1 & 0.99265217 & 0.17814460 & -0.26312433 & 6.65515541 & 1.00 & 59 \\
L2 Lyapunov & 3:2 & 1.02557297 & 0 & 0.77068285 & 4.43677028 & 115.15 & 40 \\
L2 Lyapunov & 1:1 & 0.99695262 & 0 & 1.64068576 & 6.65515541 & 49.78 & 59 \\
\bottomrule
\end{tabular}
\bigskip
Number of slots $b$ is based on $\Delta t_b = 12$ \SI{}{hours}, resulting in 1212 slots in total.
\label{tab:candidate_lpos_info}
\end{table}

\subsection{LPOs used in TE-$p$-MP}
\label{appendix:appendix_lpos_used_by_solutions}
Tables ~\ref{tab:lpos_used_cone} and~\ref{tab:lpos_used_let} provide a summary of LPOs used by solutions to the TE-$p$-MP by B\&B with a time limit of 3600 seconds and the LM with a time limit of 500 seconds. 

\begin{table}[h]
    \small
    \caption{LPOs used in TE-$p$-MP instances with Cone of Shame demand from Branch-and-Bound and Lagrangian method solutions}
    \begin{subtable}[h]{0.45\textwidth}
        \centering
        \caption{$p = 2$}
        \begin{tabular}{lll}
            \hline
            $\bar{m}_{\mathrm{crit}}$ & B\&B (3600 \SI{}{sec}) & Lagrangian Method \\ \hline
15   & \begin{tabular}[c]{@{}l@{}} L2 Halo (N) (9:2) $\times$ 1 \\L2 Halo (S) (9:2) $\times$ 1 \end{tabular} & \begin{tabular}[c]{@{}l@{}} Btrfly. (S) (1:1) $\times$ 1 \\Btrfly. (N) (1:1) $\times$ 1 \end{tabular} \\ \midrule
18   & \begin{tabular}[c]{@{}l@{}} L1 Lyap. (1:1) $\times$ 2 \end{tabular} & \begin{tabular}[c]{@{}l@{}} L1 Lyap. (1:1) $\times$ 1 \\L2 Lyap. (1:1) $\times$ 1 \end{tabular} \\ \midrule
20   & \begin{tabular}[c]{@{}l@{}} L2 Lyap. (3:2) $\times$ 1 \\L2 Lyap. (1:1) $\times$ 1 \end{tabular} & \begin{tabular}[c]{@{}l@{}} L1 Lyap. (1:1) $\times$ 1 \\L2 Halo (N) (3:1) $\times$ 1 \end{tabular} \\ \bottomrule
        \end{tabular}
       \label{tab:lpos_used_cone_p2}
    \end{subtable}
    \hfill
    \begin{subtable}[h]{0.45\textwidth}
        \centering
        \caption{$p = 3$}
        \begin{tabular}{lll}
            \hline
            $\bar{m}_{\mathrm{crit}}$ & B\&B (3600 \SI{}{sec}) & Lagrangian Method \\ \hline
15   & \begin{tabular}[c]{@{}l@{}} L2 Halo (N) (9:2) $\times$ 1 \\Btrfly. (S) (1:1) $\times$ 1 \\L2 Halo (S) (9:2) $\times$ 1 \end{tabular} & \begin{tabular}[c]{@{}l@{}} Btrfly. (S) (1:1) $\times$ 1 \\Btrfly. (N) (1:1) $\times$ 2 \end{tabular} \\ \midrule
18   & \begin{tabular}[c]{@{}l@{}} L2 Lyap. (3:2) $\times$ 1 \\Btrfly. (N) (2:1) $\times$ 1 \\L2 Lyap. (1:1) $\times$ 1 \end{tabular} & \begin{tabular}[c]{@{}l@{}} L1 Lyap. (1:1) $\times$ 3 \end{tabular} \\ \midrule
20   & \begin{tabular}[c]{@{}l@{}} L2 Lyap. (1:1) $\times$ 3 \end{tabular} & \begin{tabular}[c]{@{}l@{}} L1 Lyap. (1:1) $\times$ 3 \end{tabular} \\ \bottomrule
        \end{tabular}
        \label{tab:lpos_used_cone_p3}
    \end{subtable}
    \\ \par\bigskip
    \begin{subtable}[h]{0.45\textwidth}
        \centering
        \caption{$p = 4$}
        \begin{tabular}{lll}
            \hline
            $\bar{m}_{\mathrm{crit}}$ & B\&B (3600 \SI{}{sec}) & Lagrangian Method \\ \hline
15   & \begin{tabular}[c]{@{}l@{}} L2 Halo (N) (9:2) $\times$ 1 \\Btrfly. (S) (1:1) $\times$ 1 \\L2 Halo (S) (9:2) $\times$ 1 \\Btrfly. (N) (1:1) $\times$ 1 \end{tabular} & \begin{tabular}[c]{@{}l@{}} Btrfly. (S) (1:1) $\times$ 2 \\Btrfly. (N) (1:1) $\times$ 2 \end{tabular} \\ \midrule
18   & \begin{tabular}[c]{@{}l@{}} L2 Lyap. (3:2) $\times$ 1 \\Btrfly. (N) (2:1) $\times$ 1 \\L2 Lyap. (1:1) $\times$ 1 \\DPO (9:4) $\times$ 1 \end{tabular} & \begin{tabular}[c]{@{}l@{}} L1 Lyap. (1:1) $\times$ 1 \\L2 Lyap. (1:1) $\times$ 3 \end{tabular} \\ \midrule
20   & \begin{tabular}[c]{@{}l@{}} L2 Lyap. (3:2) $\times$ 1 \\Btrfly. (N) (2:1) $\times$ 1 \\L2 Lyap. (1:1) $\times$ 1 \\DPO (9:4) $\times$ 1 \end{tabular} & \begin{tabular}[c]{@{}l@{}} L1 Lyap. (1:1) $\times$ 4 \end{tabular} \\ \bottomrule
        \end{tabular}
       \label{tab:lpos_used_cone_p4}
    \end{subtable}
    \hfill
    \begin{subtable}[h]{0.45\textwidth}
        \centering
        \caption{$p = 5$}
        \begin{tabular}{lll}
            \hline
            $\bar{m}_{\mathrm{crit}}$ & B\&B (3600 \SI{}{sec}) & Lagrangian Method \\ \hline
15   & \begin{tabular}[c]{@{}l@{}} L2 Halo (N) (4:1) $\times$ 1 \\Btrfly. (S) (1:1) $\times$ 1 \\Btrfly. (N) (3:2) $\times$ 1 \\Btrfly. (N) (1:1) $\times$ 2 \end{tabular} & \begin{tabular}[c]{@{}l@{}} L1 Lyap. (9:4) $\times$ 1 \\Btrfly. (N) (1:1) $\times$ 4 \end{tabular} \\ \midrule
18   & \begin{tabular}[c]{@{}l@{}} L2 Halo (N) (9:4) $\times$ 1 \\L2 Lyap. (3:2) $\times$ 1 \\Btrfly. (N) (2:1) $\times$ 1 \\L2 Lyap. (1:1) $\times$ 1 \\DPO (9:4) $\times$ 1 \end{tabular} & \begin{tabular}[c]{@{}l@{}} L1 Lyap. (1:1) $\times$ 3 \\L2 Lyap. (1:1) $\times$ 2 \end{tabular} \\ \midrule
20   & \begin{tabular}[c]{@{}l@{}} L1 Lyap. (1:1) $\times$ 3 \\L1 Lyap. (3:2) $\times$ 1 \\DRO (2:1) $\times$ 1 \end{tabular} & \begin{tabular}[c]{@{}l@{}} L1 Lyap. (1:1) $\times$ 5 \end{tabular} \\ \bottomrule
        \end{tabular}
        \label{tab:lpos_used_cone_p5}
     \end{subtable}
     \label{tab:lpos_used_cone}
\end{table}

\begin{table}[h]
    \small
    \caption{LPOs used in TE-$p$-MP instances with LET demand from Branch-and-Bound and Lagrangian method solutions}
    \begin{subtable}[h]{0.45\textwidth}
        \centering
        \caption{$p = 2$}
        \begin{tabular}{lll}
            \hline
            $\bar{m}_{\mathrm{crit}}$ & B\&B (3600 \SI{}{sec}) & Lagrangian Method \\ \hline
15   & \begin{tabular}[c]{@{}l@{}} L2 Halo (N) (2:1) $\times$ 1 \\L2 Halo (S) (2:1) $\times$ 1 \end{tabular} & \begin{tabular}[c]{@{}l@{}} L2 Halo (S) (2:1) $\times$ 2 \end{tabular} \\ \midrule
18   & \begin{tabular}[c]{@{}l@{}} L1 Lyap. (1:1) $\times$ 2 \end{tabular} & \begin{tabular}[c]{@{}l@{}} L2 Halo (N) (3:1) $\times$ 2 \end{tabular} \\ \midrule
20   & \begin{tabular}[c]{@{}l@{}} L1 Lyap. (1:1) $\times$ 2 \end{tabular} & \begin{tabular}[c]{@{}l@{}} L2 Halo (N) (3:1) $\times$ 2 \end{tabular} \\ \bottomrule
        \end{tabular}
       \label{tab:lpos_used_let_p2}
    \end{subtable}
    \hfill
    \begin{subtable}[h]{0.45\textwidth}
        \centering
        \caption{$p = 3$}
        \begin{tabular}{lll}
            \hline
            $\bar{m}_{\mathrm{crit}}$ & B\&B (3600 \SI{}{sec}) & Lagrangian Method \\ \hline
15   & \begin{tabular}[c]{@{}l@{}} L2 Halo (N) (2:1) $\times$ 1 \\L2 Halo (S) (2:1) $\times$ 2 \end{tabular} & \begin{tabular}[c]{@{}l@{}} L2 Halo (S) (2:1) $\times$ 3 \end{tabular} \\ \midrule
18   & \begin{tabular}[c]{@{}l@{}} L2 Lyap. (3:2) $\times$ 1 \\Btrfly. (N) (2:1) $\times$ 1 \\L2 Lyap. (1:1) $\times$ 1 \end{tabular} & \begin{tabular}[c]{@{}l@{}} L2 Halo (N) (3:1) $\times$ 3 \end{tabular} \\ \midrule
20   & \begin{tabular}[c]{@{}l@{}} L2 Lyap. (3:2) $\times$ 1 \\Btrfly. (N) (2:1) $\times$ 1 \\L2 Lyap. (1:1) $\times$ 1 \end{tabular} & \begin{tabular}[c]{@{}l@{}} L2 Halo (N) (3:1) $\times$ 3 \end{tabular} \\ \bottomrule
        \end{tabular}
        \label{tab:lpos_used_let_p3}
    \end{subtable}
    \\ \par\bigskip
    \begin{subtable}[h]{0.45\textwidth}
        \centering
        \caption{$p = 4$}
        \begin{tabular}{lll}
            \hline
            $\bar{m}_{\mathrm{crit}}$ & B\&B (3600 \SI{}{sec}) & Lagrangian Method \\ \hline
15   & \begin{tabular}[c]{@{}l@{}} L2 Halo (N) (2:1) $\times$ 2 \\L2 Halo (S) (2:1) $\times$ 2 \end{tabular} & \begin{tabular}[c]{@{}l@{}} L2 Halo (S) (2:1) $\times$ 4 \end{tabular} \\ \midrule
18   & \begin{tabular}[c]{@{}l@{}} L2 Lyap. (3:2) $\times$ 1 \\Btrfly. (N) (2:1) $\times$ 1 \\L2 Lyap. (1:1) $\times$ 1 \\DPO (9:4) $\times$ 1 \end{tabular} & \begin{tabular}[c]{@{}l@{}} L2 Halo (N) (3:1) $\times$ 4 \end{tabular} \\ \midrule
20   & \begin{tabular}[c]{@{}l@{}} DRO (3:2) $\times$ 4 \end{tabular} & \begin{tabular}[c]{@{}l@{}} L2 Halo (N) (3:1) $\times$ 4 \end{tabular} \\  \bottomrule
        \end{tabular}
       \label{tab:lpos_used_let_p4}
    \end{subtable}
    \hfill
    \begin{subtable}[h]{0.45\textwidth}
        \centering
        \caption{$p = 5$}
        \begin{tabular}{lll}
            \hline
            $\bar{m}_{\mathrm{crit}}$ & B\&B (3600 \SI{}{sec}) & Lagrangian Method \\ \hline
15   & \begin{tabular}[c]{@{}l@{}} L2 Halo (N) (2:1) $\times$ 2 \\L2 Halo (S) (2:1) $\times$ 3 \end{tabular} & \begin{tabular}[c]{@{}l@{}} L2 Halo (N) (2:1) $\times$ 1 \\L2 Halo (S) (2:1) $\times$ 4 \end{tabular} \\ \midrule
18   & \begin{tabular}[c]{@{}l@{}} L2 Halo (N) (9:4) $\times$ 1 \\L2 Lyap. (3:2) $\times$ 1 \\Btrfly. (N) (2:1) $\times$ 1 \\L2 Lyap. (1:1) $\times$ 1 \\DPO (9:4) $\times$ 1 \end{tabular} & \begin{tabular}[c]{@{}l@{}} L2 Halo (N) (3:1) $\times$ 5 \end{tabular} \\ \midrule
20   & \begin{tabular}[c]{@{}l@{}} L2 Halo (N) (3:1) $\times$ 1 \\DRO (9:4) $\times$ 1 \\DRO (3:2) $\times$ 3 \end{tabular} & \begin{tabular}[c]{@{}l@{}} L1 Lyap. (1:1) $\times$ 1 \\L2 Halo (N) (3:1) $\times$ 4 \end{tabular} \\ \bottomrule
        \end{tabular}
        \label{tab:lpos_used_let_p5}
     \end{subtable}
     \label{tab:lpos_used_let}
\end{table}

\bibliography{references}

\end{document}